%% file: main_CNL.tex
\documentclass[opre,nonblindrev]{informs3wp}

\OneAndAHalfSpacedXI 

\usepackage[pdfpagemode={UseOutlines},bookmarks=true,bookmarksopen=true,
bookmarksopenlevel=0,bookmarksnumbered=true,hypertexnames=false,
colorlinks,linkcolor={blue},citecolor={blue},urlcolor={red},
pdfstartview={FitV},unicode,breaklinks=true]{hyperref}
\hypersetup{urlcolor=blue, colorlinks=true}

\usepackage{endnotes}
\let\footnote=\endnote

\usepackage{natbib}
\bibpunct[, ]{(}{)}{,}{a}{}{,}%
\def\bibfont{\small}%
\TheoremsNumberedThrough  
\ECRepeatTheorems

\EquationsNumberedThrough 

\usepackage[justification=centering]{caption}
\usepackage{subcaption}
\usepackage{float}
\usepackage{graphicx}
\usepackage{tablefootnote}
\usepackage{booktabs}
\usepackage{multirow}
\usepackage{hhline}
\usepackage{footnote}
\usepackage[ruled,vlined]{algorithm2e}
\usepackage{algorithmic}
\usepackage[title]{appendix}
\usepackage{amsmath,tabu}

\renewcommand{\endproof}{\strut\hfill\Halmos\endTrivlist\addvspace{0pt}}

\newcommand{\cX}{{\mathcal{X}}}

\newcommand{\cO}{{\mathcal{O}}}

\newcommand{\bx}{\textbf{x}}
\newcommand{\by}{\textbf{y}}
\newcommand{\bs}{\textbf{s}}

\newcommand{\bA}{\textbf{A}}

\newcommand{\ba}{\textbf{a}}

\newcommand{\bz}{\textbf{z}}
\newcommand{\bb}{\textbf{b}}

\newcommand{\bbR}{\mathbb{R}}

\newtheorem{note}{Note}

\usepackage{mathtools}

\usepackage{framed} 

\DeclarePairedDelimiter\ceil{\lceil}{\rceil}

\newif\ifnotes\notestrue
\def\htien#1{}

\usepackage{xcolor}

\begin{document}
	
	
	\RUNAUTHOR{ Cuong Le and Tien Mai}
	
	\RUNTITLE{Constrained Assortment Optimization under Cross-Nested Logit}
	
	\TITLE{Constrained Assortment Optimization\\ under the Cross-Nested Logit Model}
	
	\ARTICLEAUTHORS{%
		\AUTHOR{Cuong Le and Tien Mai}
		\AFF{School of Computing and Information Systems, Singapore Management University\\
						\EMAIL{cuonglv@smu.edu.sg}, \EMAIL{atmai@smu.edu.sg}}
		\today
	} 
	
	\ABSTRACT{%
		We study the assortment optimization problem under general linear constraints, where the customer choice behavior is captured by the Cross-Nested Logit model. In this problem, there is a set of products organized into multiple subsets (or nests), where each product can belong to more than one nest. The aim is to find an assortment to offer to customers so that the expected revenue is maximized. We show that, under the Cross-Nested Logit model, the assortment problem is NP-hard, even without any constraints. To tackle the assortment optimization problem, we develop a new discretization mechanism to approximate the problem by a linear fractional program
  with a performance guarantee of $\frac{1 - \epsilon}{1+\epsilon}$, for any accuracy level $\epsilon>0$. We then show that optimal solutions to the approximate problem can be obtained by solving mixed-integer linear programs. We further show that our discretization approach  can also be applied to solve a joint assortment optimization and pricing problem, as well as an assortment problem under a mixture of Cross-Nested Logit models to account for multiple classes of customers. Our empirical results on a large number of randomly generated test instances demonstrate that, under a performance guarantee of 90\%, the percentage gaps between the objective values obtained from our approximation methods and the optimal expected revenues are no larger than 1.2\%. 
	}%
	
	
	
	\maketitle

	
	%
	
	\section{Introduction}
	In revenue management, assortment optimization is a fundamental problem in which retailers need to select a subset of products to offer for sale, whether in a physical store, an online marketplace, or any other retail environment. The objective is typically to maximize the expected revenue obtained from the customers while taking into account some business constraints such as available budgets or limitations on shelf space. In reality, the products may serve as substitutes while the customers may make a choice among the products that satisfy their needs. These create a situation where the demand for a specific product depends on what assortment of products is offered. For example, the customer only purchases product A when an equivalent product B in the same category is unavailable. Hence, it is crucial for the retailer to thoroughly examine the underlying relationships between the products and the purchasing patterns of customers before making any decisions.
	
	The first work that directly incorporates customer purchase behaviors into the assortment optimization was conducted by \citet{talluri2004revenue}, in which the authors modeled the choice process of airline consumers by using the multinomial logit model (MNL). This choice model assumes that the unobserved utilities of all products are independent (also known as the Independent from Irrelevant Alternatives, or IIA assumption). Two decades after this study, the MNL model has been widely applied and has achieved great success in the literature of revenue management  \citep[see][for instance]{rusmevichientong2012robust, abdallah2021demand, jagabathula2022personalized}. However, studies have also shown that the IIA assumption could be unrealistic in many practical settings. In such situations, extensions such as the nested logit (NL) or cross-nested logit models are widely preferred as alternatives to account for the violation of the IIA property \citep{davis2014assortment, feldman2015capacity}. In the NL model, products are partitioned into several non-overlapping groups of products (which are often called nests), where the IIA holds within each nest and is relaxed over products from different nests. Despite the nesting structure, the NL model still has limitations resulting from the unambiguous assignment of products to nests, which implies that the cross-elasticities are the same for all pairs of products in each nest, making the NL model unable to correctly capture complex choice processes that arise from intricate  topological relationships.
	
	In this paper, we consider an assortment optimization problem with general linear constraints while the customer's choice behavior is captured by the cross-nested logit (CNL) model. The CNL model allows each product to belong to different nests with different levels of membership. This flexible structure makes CNL powerful for capturing a wide range of complex relationships between products, where the correlations can not be handled by the MNL or NL models.  The CNL model is also a member of the Generalized Extreme Value (GEV) family, and thus it is consistent with the theoretical foundation of the random utility theory \citep{McFa78, bierlaire2006theoretical}.  Compared to other choice models like MNL and NL, the CNL model is more general and more sophisticated. In fact, the CNL is fully flexible as it can approximate any random utility maximization model \citep{fosgerau2013choice}. Nevertheless,  the arisen optimization problem is also significantly harder to handle. As far as we know, there is still an open question regarding how to achieve nearly optimal solutions for assortment optimization under the CNL model. Our work addresses this question.

In the paper, we study, for the first time, a constrained assortment optimization problem under the CNL model. We show that the assortment optimization problem is NP-hard even with no constraint. We then present a discretization mechanism that enables us to approximate the optimization problem by a linear fractional program, which can be further resolved by solving a mixed-integer linear program.
 Additionally, we show that the proposed method can also be applied to a joint assortment and price optimization problem, as well as an extension that considers multiple classes of customers (i.e. assortment optimization under a mixture of customer types). We summarize our main contributions in the following. 
	
	\subsection{Our Contributions}
	
	\noindent\textbf{Problem formulation under the CNL model}. To the best of our knowledge, this is the first study considering the assortment optimization problem under the CNL model. We present a problem formulation with general linear constraints that could be capacity, knapsack, or partition constraints, or a combination of any possible linear constraints. It should be noted that, in order to simplify the optimization process, previous studies in the literature have often examined the assortment problem with a single constraint rather than combining multiple constraints together. The requirement of satisfying multiple constraints arises in many real-world situations. For example, a retailer may need to consider both the total budget allocated for products within each nest as well as the total amount of shelf space that all products across all nests will occupy, or require that some products should or should not be selected together.  We show that the assortment problem under the CNL model is \textit{NP-hard}, even without constraints. The proof is based on a reduction from the \textit{set partition} problem, which is a well-known NP-complete problem \citep{karp2010reducibility}.
	
	\noindent\textbf{Approximation and linear reformulations}. We develop a discretization mechanism to approximate the original problem by a linear-fractional program (LFP), which is easier to handle. Our approximation method relies on a piecewise linear approximation approach, where the number of linear segments required to achieve a desired accuracy level is minimized, leading to an optimal size for the approximate program. We show that solving the approximate program can yield a solution with a worst-case performance guarantee of $\frac{1 - \epsilon}{1 + \epsilon}$ for any given accuracy level $\epsilon>0$. We then develop two methods for solving the approximate LFP. The first method is based on a parametric approach in which a sequence of feasible linear programs is solved to find a near-optimal solution. Our second approach involves converting the approximated LFP into a mixed-integer linear program (MILP) and solving it using a readily available solver, such as CPLEX, to obtain a near-optimal solution.
	
	\noindent\textbf{Joint assortment and price optimization, and mixture of CNL models}. We further consider the problem of jointly optimizing the assortment and prices of offered products, where there is a finite set of possible prices for each product and the customer choice behavior is governed by the CNL model. Our approach is to convert this joint problem into an equivalent assortment problem by treating each product and its associated prices as multiple products, where each product corresponds to a single price. This way allows us to apply the proposed approximation method as well as the linear reformulations to solve the joint assortment-price problem and obtain a solution with the same performance guarantee as for the assortment problem.
	
  We also consider an assortment problem under a mixture of multiple CNL models, which arises from the availability of multiple customer classes. 
  For this problem, we apply the proposed approximation method to each customer class, resulting in a sum-of-ratios LFP problem that can be solved to obtain a $\frac{1 - \epsilon}{1 + \epsilon}$ solution, for any given $\epsilon>0$. We also provide a MILP reformulation for the approximate program, which can be further solved to optimality by an off-the-shelf solver.
	
	\noindent\textbf{Experimental study}. Extensive experiments are conducted on the formulated problems to evaluate the performance of the proposed methods. To accomplish this, we generate a large number of test instances for each problem and run our methods under different settings of performance guarantee. Numerical results show that our proposed methods perform impressively well. In particular, under the performance guarantee of 90\%, the average gaps between the expected revenues yielded by the proposed methods and those obtained from offering optimal assortments are no larger than 1.2\%.

 In summary, we consider constrained assortment optimization under the CNL model or a mixture of CNL models. We propose, for the first time, solution methods that offer  $\frac{1 - \epsilon}{1 + \epsilon}$ approximation solutions. Our methods involve solving  MILPs by an off-the-shelf solver. 
Our solution methods would be of high value in practice; our approach allows for the application of highly complex choice models and works with any linear constraints. Moreover, the use of commercial software
allows for the inclusion of complex business constraints without the need for redesigning algorithms, and is
continually improved with new optimization methods and advanced hardware.

	\subsection{Literature Review}
	In this section, our discussion will primarily center around studies that have utilized discrete choice models, such as the MNL and NL models to address the assortment optimization problem. We will also discuss the advantages and existing applications of the CNL model.
	
	Following the work of \citet{talluri2004revenue}, there is a rich literature on  assortment optimization using the MNL model to predict  customers' product choice probabilities. For example,  \citet{rusmevichientong2010dynamic} study a joint parameter estimation and assortment optimization with a capacity constraint. \citet{rusmevichientong2012robust} investigate robust formulations of the assortment optimization problem under the MNL model with uncertain choice parameters. \citet{bront2009column} and \citet{mendez2014branch} propose mixed-integer linear reformulations and a branch-and-cut algorithm to obtain  optimal solutions for assortment optimization under the  mixed logit model. Furthermore, under the mixed logit model, \citet{rusmevichientong2014assortment} show that the assortment problem is NP-hard even when the choice model is a mixture of only two MNL models. \citet{desir2014near} develop a fully polynomial-time approximation scheme (FPTAS) for the capacitated assortment problem with a runtime that grows exponentially as the number of customer types increases. \citet{sen2018conic} develop a conic quadratic mixed-integer program to handle large instances of the capacitated assortment problem under mixed logit. \cite{abdallah2021demand} study a  joint assortment optimization and customization problem, in which the firm first offers an assortment with consideration over all the customer types and then customizes the assortment based on a specific type of the arrived customer.
	
	Under the nested logit model, \citet{davis2014assortment} show that the unconstrained assortment optimization problem can be solved in polynomial time when the dissimilarity parameters of all the nests are less than one and the customer always makes a purchase after choosing a nest. Relaxing one of these conditions makes the problem NP-hard. Not surprisingly, the majority of previous studies on assortment optimization under the NL model have concentrated on more challenging problems in the latter direction. For an assortment optimization problem defined on $m$ products and $n$ nests without any constraint, \citet{li2014greedy} propose a greedy algorithm with a complexity of $\cO(mn\log n)$ that  removes at most one product from each nest at each iteration. \citet{gallego2014constrained} study an assortment problem under the NL model where the capacity constraint is applied separately in each nest. They proposed an approximation algorithm with a performance guarantee of 0.5 under certain conditions on the products. \citet{rusmevichientong2009ptas} and \citet{feldman2015capacity} develop 4-approximation algorithms to handle the assortment problem under the NL model where the total capacity consumption of the whole assortment is upper bounded.
	
	Our work is also relevant to studies that jointly optimize the assortment and price. \citet{wang2012capacitated} and \citet{chen2020capacitated} study a capacitated assortment and price optimization problem under different logit models where they show that  the joint problems can be  converted into a fixed point problem. \citet{gallego2014constrained} transform the joint assortment optimization and pricing problem under the NL model into an assortment problem by creating multiple virtual copies for each product. \citet{gao2021assortment} study the assortment optimization and pricing problem with impatient customers under the MNL model, where customers incrementally view products in multiple stages. They show that for a given assortment, optimal prices can be computed by a convex program. For the joint assortment-price optimization, they develop an approximation algorithm with a performance guarantee of 87.8\%, in which all products are offered in the first stage with their optimal prices. \citet{miao2021dynamic} study a joint assortment and pricing problem under the MNL model with unknown parameters where customers are assumed to arrive sequentially. The authors then develop a learning algorithm considering the trade-off between parameter learning and revenue extraction.
	
	As mentioned previously, both MNL and NL models  have their own limitations that prevent them from accurately reflecting complex correlations in many real-world scenarios. The IIA assumption in the MNL model implies that the cross-elasticities of all products caused by the change in one product's attributes are identical. This assumption, however, may be too restrictive when it comes to the assortment problem, where the products are often complements or substitutes to others. In the case of the NL,  products are partitioned  into disjoint nests, which relaxes the IIA assumption for products belonging to different nests. Nevertheless, the IIA property still holds within each nest, implying that the correlation levels between all product pairs in each nest must be the same. Besides, the strict partitioning in the NL model implies that the correlation between products must be transitive. For example, if product A is correlated with product B, and B is correlated with product C, the NL model will group all three products into the same nest, regardless of whether A is correlated to C. 
	
	There is a work by \citet{zhang2020assortment} and an extension by \citet{ghuge2022constrained} that handle the above non-transitive correlation in the context of the assortment problem. The authors employ the Paired Combinatorial Logit (PCL) model to capture the choice process, where each pair of products institutes a nest with its own degree of independence. They show that under the PCL model, the unconstrained assortment problem is strongly NP-hard and develop polynomial-time approximation algorithms to solve the problem under different settings. However, under the PCL model, with $m$ products, the uniform allocation of a product to $m-1$ nests limits its correlations with each other to $1/(m-1)$ and thus restricts the magnitude of cross-elasticities of all the pairs. This limitation may become a serious issue if there are more than a few products \citep[see][for more details]{koppelman2000closed, wen2001generalized}. In contrast, the CNL model permits assigning various proportions of each item to multiple nests. This design determines a product's correlation with others based on the allocation fraction to their common nest, thereby offering greater flexibility in the degree of cross-elasticities.
	
	Regarding the CNL model, this model was introduced lately by \cite{small1987discrete, vovsha1997application, ben1999discrete}. \citet{bierlaire2006theoretical} established a theoretical foundation for the CNL model by formally proving that the model belongs to the GEV family \citep{McFa78} and providing a novel estimation procedure based on non-linear programming instead of heuristics like previous studies. This work also indicates that the MNL and NL models are indeed special cases of the CNL model.
 \cite{fosgerau2013choice} show that the CNL is fully flexible as it is able to approximate any random utility maximization (RUM) model.
 With a flexible structure, the CNL model has successfully been applied to model the choice behavior in a wide range of transportation problems. To name a few, the applications of the CNL model can be found in the mode choice problem  \citep{yang2013cross}, departure time choice \citep{ding2015cross}, route choice  \citep{lai2015specification,mai2016method}, air travel management  \citep{drabas2013modelling}, and recently, the location choice in international migration \citep{beine2021new}. Despite the success in transportation management, the potentials of the CNL model in revenue management still remain unexplored. To the best of our knowledge, this is the first study employing the CNL model to capture substitution behavior in  assortment optimization.
	
	\subsection{Organization}
	The rest of this paper is organized as follows. In Section \ref{sec:formulation-assort}, we provide a formulation and an example of the assortment problem where the customers' behavior is captured by the CNL model. In Section \ref{sec:approximation}, we present our proposed discretization mechanism to approximate the formulated assortment optimization problem. To solve the  approximate problem, two solution methods are given in Section \ref{sec:solution_methods}. In Section \ref{sec:formulation-jap}, we consider the joint assortment optimization and pricing problem under the CNL model. Section \ref{sec:mixture} provides a further extension to a mixture of CNL models to account for the availability of multiple customer classes. The experimental results are provided in Section \ref{sec:experiment}, where we evaluate the performance of the proposed methods through extensive experiments on the formulated problems.  
    Section \ref{sec:conclu} concludes the paper.  The appendix contains supplemental analyses and proofs that are not included in the main paper.
	
	\noindent
	\textbf{Notation:}
	Boldface characters represent matrices or vectors or sets, and $a_i$ denotes the $i$-th element of $\ba$ if $\ba$ is indexable. We use $[m]$, for any $m\in \mathbb{N}$, to denote the set $\{1,\ldots,m\}$. 
	
	\section{Assortment Optimization  under the CNL Model}
	\label{sec:formulation-assort}
	In this section, we formulate the constrained assortment problem under the CNL model. We also give an example where the correlation between products cannot be captured by  MNL or NL models.
	
	\subsection{Problem Formulation and NP-hardness}
	Consider a set of $m$ products indexed by $\{1, 2, \ldots, m\}$, where each product $i$ has a predefined revenue $r_i$ and a preference weight (or utility) $v_i$. We assume that $m$ products have been assigned to $n$ subsets (or nests) $S_1, S_2, \ldots, S_N \subseteq [m]$ according to their attributes. Noting that $S_1, S_2, \ldots, S_N$ are not necessary to be disjoint, i.e., a product $i$ can belong to multiple nests. In this case, we use a non-negative quantity $\alpha_{in}$ to capture the level of membership of product $i$ in the nest $S_n$, where $\sum_{n\in[N]}\alpha_{in} = 1\ \forall{i\in[m]}$. Without prejudice to the generality of the formulation, we assume that $S_n = [m]\ \forall{n\in [N]}$, and if a product $i$ is not a member of $S_n$, we can set  $\alpha_{in} = 0$ without affecting the previous assumption. 
	
	Let $\bx = (x_1, x_2, \ldots, x_m) \in \{0, 1\}^m$ be a binary vector representing an assortment  decision where $x_i = 1$ if, and only if, the product $i$ is offered. Under the CNL model, the customer's choice behavior can be decomposed into two stages where each of these takes the form of the standard logit. In the first stage, the customer selects a nest $S_n$ from $N$ nests. Each nest $S_n$ is thus selected with a probability $P(S_n|\bx) = W_n^{\sigma_n}/ \sum_{n'\in[N]}W_{n'}^{\sigma_{n'}}$, where $W_n = v_{0n} + \sum_{i\in S_n}\alpha_{in}x_iv_i$ is the total preference weights of all alternatives in $S_n$, $\sigma_n$ is the dissimilarity parameter of $S_n$, and $v_{0n}$ is the preference weight for leaving nest $S_n$ without purchasing anything. It is typically assumed that the value of $\sigma_n$ varies in the unit interval for all nests to guarantee that the model is consistent with the RUM framework \citep{McFa78,bierlaire2006theoretical}.  In the second stage, the customer decides to leave or purchase a product $i$ from the nest $S_n$ that was selected in the previous stage. In this case,  the conditional probability of choosing product $i$ is given by $P(i|S_n) = \alpha_{in}x_iv_i / W_n$. It is worth noting that the two-stage decomposition of the decision-making process is used to facilitate the representation of the correlation structure between products. Nevertheless, customers do not necessarily follow these stages when making a purchase decision in reality. Given an assortment $\bx$, the CNL model posits that a product $i \in [m]$ is purchased with the probability
	\begin{align*}
		P(i|\bx) = \sum_{n\in [N]}P(S_n|\bx) \cdot P(i|S_n) &= \sum_{n\in[N]}\frac{W_{n}^{\sigma_n}}{ \sum_{n'\in[N]}W_{n'}^{\sigma_{n'}}}\times \frac{\alpha_{in}x_iv_{i}}{W_{n}}\\
		&=\frac{\sum_{n\in [N]}{W_n}^{\sigma_n-1} (\alpha_{in}x_iv_i)}{ \sum_{n\in [N]} W_{n}^{ \sigma_{n}}}.
	\end{align*}
	By the law of total expectation, the expected revenue obtained from the customer by offering the assortment $\bx$ can  be calculated as
	\begin{align}
		F(\bx)  = \frac{\sum_{i\in [m]} \sum_{n\in [N]}{W_n}^{\sigma_n-1} (\alpha_{in}x_ir_iv_i)}{\sum_{n\in [N]} W_{n}^{ \sigma_{n}}}. \nonumber
	\end{align}
	In this study, we consider the assortment optimization problem under general linear constraints. In particular, let $\cX = \{\bx \in \{0, 1\}^m\ |\ \bA\bx \leq \bb\}$ denote the set of all feasible assortments. Then $\bA$ and $\bb$ can be created using various methods to represent diverse forms of business constraints, depending on the requirements of the actual situation. With all of these in mind, the assortment optimization problem with general linear constraints under the CNL model can be formulated as follows:
	\begin{align}
		\max_{\bx \in \cX }\qquad &\left\{F(\bx)  = \frac{\sum_{i\in [m]} \sum_{n\in [N]}{W_n}^{\sigma_n-1} (\alpha_{in}x_ir_iv_i)}{\sum_{n\in [N]} W_{n}^{ \sigma_{n}}}\right\}\label{prob:CNL-assort}\tag{\sf Assort} \\
        \text{subject to} \quad &W_n = v_{0n} + \sum_{i\in S_n}\alpha_{in}x_iv_i, \forall{n \in [N]}. \nonumber
	\end{align}
	In the following theorem, we show that the problem \eqref{prob:CNL-assort} is NP-hard, even without any constraint. 
    \begin{theorem}
    \label{theorem:np_hard}
        The unconstrained assortment optimization problem under the CNL model \eqref{prob:CNL-assort} is NP-hard.
    \end{theorem}
    The proof of Theorem \ref{theorem:np_hard} can be found in Appendix \ref{appendix:proof_np_hard}, where we construct a reduction from the \textit{set partition} problem \citep{karp2010reducibility}. From the perspective of computational complexity, we will not restrict ourselves to polynomial-time algorithms for addressing the constrained assortment problem \eqref{prob:CNL-assort}.
    
    \subsection{Example of Cross-Nested Correlation Structure}
    \begin{figure}
        \centering
        \includegraphics{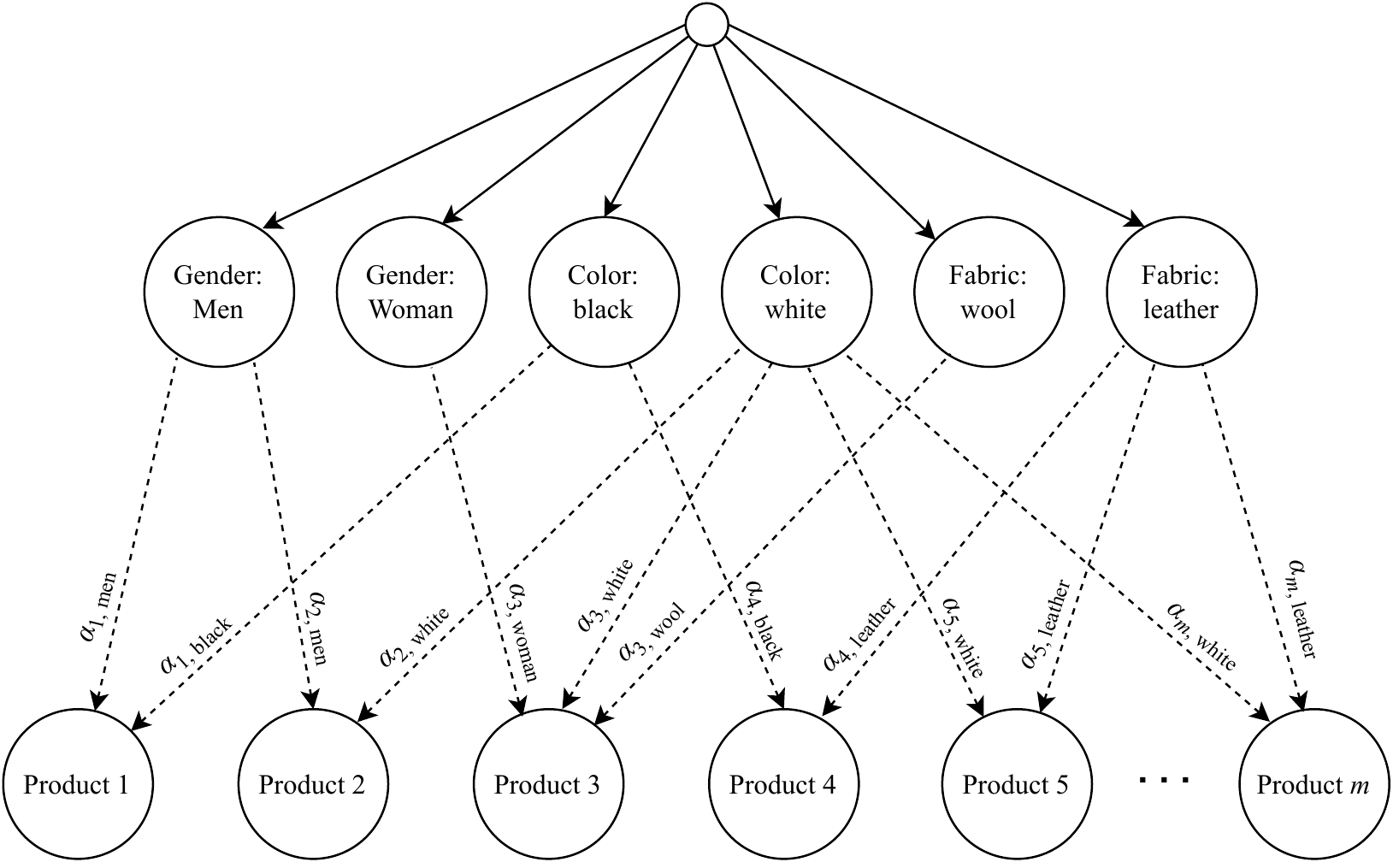}
        \caption{An example of a cross-nested correlation structure.}
        \label{fig:example}
    \end{figure}
    It is not difficult to give examples where the products are correlated to others according to a cross-nested structure. Let's consider the fashion retail industry, where the correlation between products has a significant impact on customer purchase behavior. It is a common strategy in the fashion industry to use correlation to maintain a consistent brand image, facilitate marketing and branding efforts, optimize production, and drive sales. Correlation is often achieved through  observations of similar colors, patterns, fabrics, or design elements, which creates a sense of harmony and unity between different items in a fashion collection. By creating a cohesive and recognizable aesthetic, a brand can encourage customers to purchase multiple items from the same collection. Correlated products can also enhance the customer experience by making it easier for customers to find complementary pieces and create cohesive outfits. Using the correlation, the brand can encourage customers to purchase products that they may not have considered otherwise, thereby expanding the range of products that customers consider purchasing.
    
    In Figure \ref{fig:example}, we provide an example of a fashion shop, where the retailer wants to select several products from a collection of $m$ items to display to the customer. To optimize production and sales, products were designed to be correlated to others. Based on the correlation, products are organized into six nests, including two nests for genders (men and women), two nests for colors (black and white), and two nests for fabrics (wool and leather). The customer's choice process can be seen as taking place in two stages. In the first stage, the customer enters the shop and selects a nest based on their shopping habits or product attributes they are most interested in. For instance, a customer who loves leather material will take the first look at leather products rather than non-leather items. In the second stage, the customer decides to buy a product within the selected nest or leave the shop. Through this example, we will show that a product can correlate to other products at different levels, and the correlations would not always be transitive.
    
    We first demonstrate that the correlation levels can differ between different product pairs by considering the first three items 1, 2 and 3, and four nests of gender and color. As illustrated in Figure \ref{fig:example}, Product 1 (black for men) and Product 2 (white for men) have a correlation as they are both designed for men. Likewise, Product 2 (white for men) and Product 3 (white for women) are correlated as they share the same color (customers who are couples can purchase them in pairs). Let us now consider a scenario where the demand for Product 2 suddenly spikes and the shop owner decides to raise its price. As a result, some customers may rethink buying this product since the value they receive from it is no longer worth the increased cost, i.e., the utility is no longer maximized. In this case, these customers are more likely to switch to Product 1 and accept a different color, but at a lower price, rather than choose Product 3 with the same color as Product 2, but designed for women. This suggests that a change in Product 2's attributes will affect  the demand for Product 1 more than for Product 3, indicating that the correlation between Product 1 and Product 2 is stronger than the correlation between Product 2 and Product 3.

    Secondly, we show that the correlation relationships between products are not always transitive. Let's consider three products in our example, including Product 3, Product 4, and Product 5. As shown in Figure \ref{fig:example}, there is a correlation between Product 3 (white and made from wool) and Product 5 (white and made from leather) because they share the same color. Similarly, Product 4 (black and made from leather) is correlated to Product 5 (white and made from leather) as they are made from the same material. However, Product 3 and Product 4 are uncorrelated since they do not share any common attributes.

     In the above two examples, the correlation structures cannot be captured by the MNL and the NL models due to their IIA assumption and disjoint nested structure. On the other hand, the flexible structure of the CNL model allows each product to belong to multiple nests with different allocation proportions, thereby providing the model the ability to work with complex correlation topologies. Here, we note that the correlations in our example can also be handled by the PCL model, as suggested by \cite{zhang2020assortment} and \cite{ghuge2022constrained}. However, while the CNL model requires at most $6m + 6$ parameters to be estimated (at most 6 allocation parameters for each product and 6 nest dissimilarity parameters), the number of model parameters in the PCL model grows quadratically with the number of products $m$, making it computationally more challenging to estimate, especially when the number of products is large.

	\section{Discretization and Approximation}
	\label{sec:approximation}
	This section presents our proposed method to approximate the problem \eqref{prob:CNL-assort} by a linear fractional program with an accuracy guarantee of $\frac{1 - \epsilon}{1 + \epsilon}$ for any given $\epsilon > 0$. The general idea is to replace the non-linear components, i.e., $W_n^{\sigma_n-1}$ and $W_n^{\sigma_n}$, in the objective function of \eqref{prob:CNL-assort} with piecewise linear approximation functions so that the obtained approximate problem is easier to handle. In Subsection \ref{subsec:approximation} below, we will show how to deal with the non-linear component of the objective function and how we transform the original problem \eqref{prob:CNL-assort} into an LFP program. Since the size of the LFP program is partly decided by the accuracy level of the approximation, specifically, a higher accuracy guarantee requires more linear segments, resulting in a larger model in size, Subsection \ref{subsec:discretization} will present a procedure to discretize the objective function in such a way that the number of linear segments required to achieve the performance guarantee of $\frac{1 - \epsilon}{1 + \epsilon}$ is minimized. Finally, in Subsection \ref{subsec:K_opt_bounds}, we demonstrate the optimality of the proposed procedure and provide the bounds on the number of generated sub-intervals.
	
	\subsection{Approximation Approach}
        \label{subsec:approximation}
	First, let $L_n$ and $U_n$ be the lower bound and upper bound of $W_n$, obtained by solving $\min/\max_{\bx \in \cX} \sum_{i\in S_n}(v_{0n} + \alpha_{in}x_{i}v_{i})$. On $[L_n, U_n]$, let us define $f^n(W_n) = W_n^{\sigma_n - 1}$ and $g^n(W_n) = W_n^{\sigma_n}$. We rewrite the objective function of the problem \eqref{prob:CNL-assort} as
    \begin{equation}
        F(\bx) = \frac{\sum_{i\in [m]} \sum_{n\in [N]}f^n(W_n) (\alpha_{in}x_ir_iv_i)}{\sum_{n\in [N]} g^n(W_n)}.
        \label{eq:obj_assort}
    \end{equation}
    To approximate this objective function, we will replace the non-linear functions $f^n(W_n)$ and $g^n(W_n)$ by piecewise linear inner-approximators $\widehat{f^n}(W_n)$ and $\widehat{g^n}(W_n)$ (we will discuss how this can be done in the next subsection). We then obtain an approximation problem of the following form
	\begin{align}
		\max_{\bx \in \cX } &\;  \left\{\widehat{F}(\bx)  = \frac{\sum_{i\in [m]} \sum_{n\in [N]}\widehat{f^n}(W_n) (\alpha_{in}x_ir_iv_i)}{\sum_{n\in [N]} \widehat{g^n}(W_n)}\right\}\label{prob:approx}\tag{\sf Assort-Approx} \\
        \text{subject to} \quad &W_n = v_{0n} + \sum_{i\in S_n}\alpha_{in}x_iv_i, \forall{n \in [N]}. \nonumber
	\end{align}
	
	\begin{figure}
		\centering
		\includegraphics[scale=0.55]{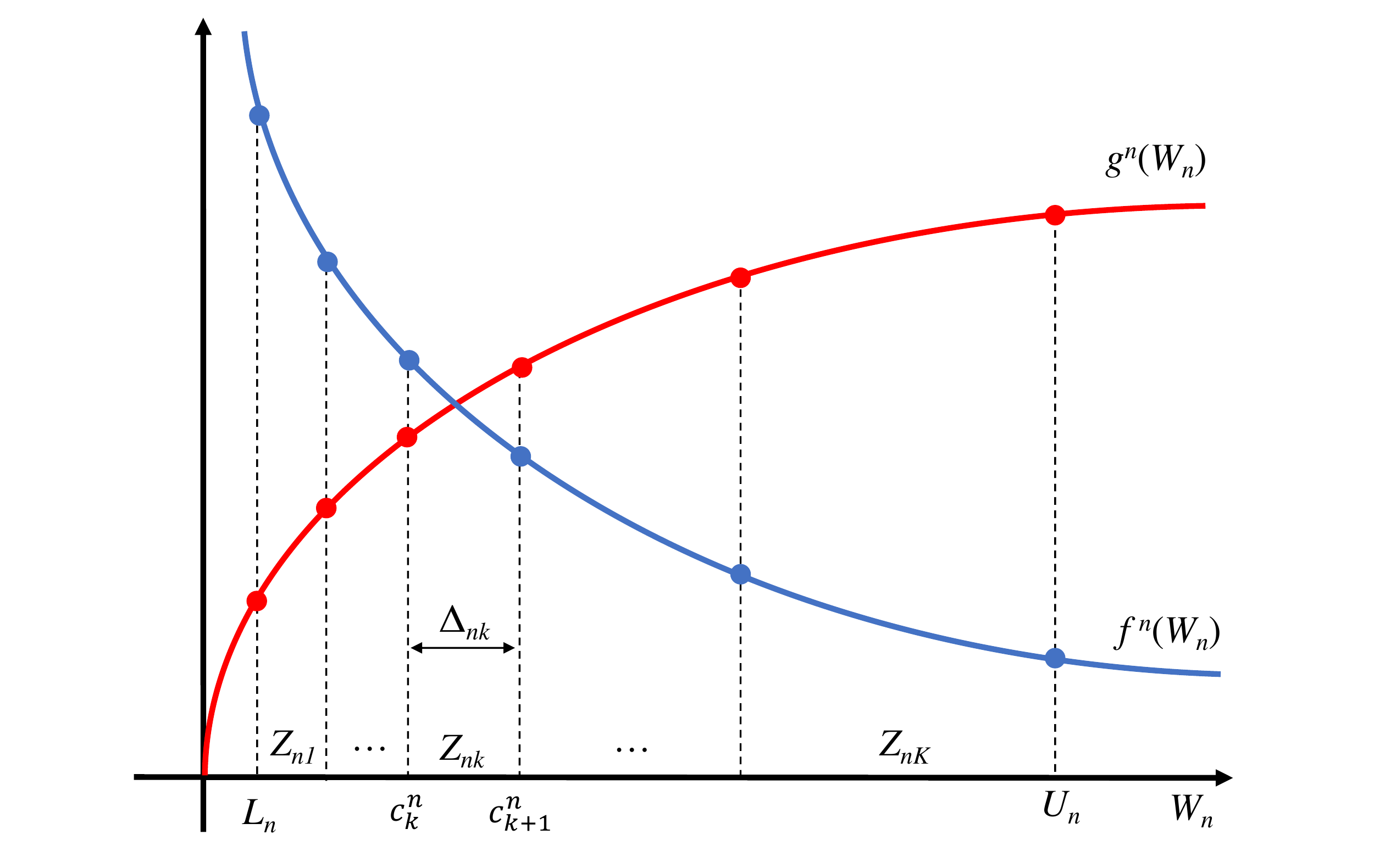}
		\caption{Illustration of the proposed approximation method.}
		\label{fig:approx_illustration}
	\end{figure}

        Recall that $\sigma_n$ is assumed to be in $[0, 1]$ for the consistency of the model with the maximum utility behavior \citep{McFa78}. It follows that the continuous relaxation of $f^n(W_n)$ is convex and decreasing, and this of $g^n(W_n)$ is concave and increasing in $W_n$, implying that $f^n(W_n) \leq \widehat{f^n}(W_n)$ and $g^n(W_n) \geq \widehat{g^n}(W_n)$ for $L_n \leq W_n \leq U_n$. We have the following approximation bounds that  are essential to establish the performance guarantees:
        
	\begin{lemma}
        \label{lemma:piecewise_function}
		Let 
		\begin{align*}
			\epsilon^{fn} = \max_{L_n \leq W_n \leq U_n}\left\{\frac{\widehat{f^n}(W_n) - f^n(W_n)}{f^n(W_n)}\right\} \quad \text{and} \quad \epsilon^{gn} = \max_{L_n \leq W_n \leq U_n}\left\{\frac{g^n(W_n) - \widehat{g^n}(W_n)}{g^n(W_n)}\right\}
		\end{align*}
		be the maximum relative approximation errors of $\widehat{f^n}(W_n)$ and $\widehat{g^n}(W_n)$ for all $n\in[N]$. Let $\epsilon = \max_{n\in [N]}\left\{\epsilon^{fn}, \epsilon^{gn}\right\}$. Then, for any feasible solution $\bx$ of problem \eqref{prob:CNL-assort}, we have
		\begin{equation*}
			F(\bx) \leq \widehat{F}(\bx) \leq \frac{1 + \epsilon}{1 - \epsilon}F(\bx).
		\end{equation*}
	\end{lemma}
	\proof{Proof.}
	The first inequality follows directly from the properties of the piecewise linear functions, i.e., $f^n(W_n) \leq \widehat{f^n}(W_n)$ and $g^n(W_n) \geq \widehat{g^n}(W_n)$ for all $n\in[N]$ and $L_n \leq W_n \leq U_n$.  
	We now prove the second inequality. From the definitions of $\epsilon^{fn}$, $\epsilon^{gn}$, and $\epsilon$, we have $\widehat{f^n}(W_n) \leq (1 + \epsilon) f^n(W_n)$ and $\widehat{g^n}(W_n) \geq (1 - \epsilon) g^n(W_n)$. Substitute these in the objective function of problem \eqref{prob:approx}, we get
	\[
	\widehat{F}(\bx) \leq \frac{\sum_{i\in [m]} \sum_{n\in [N]}(1+\epsilon){f^n}(W_n) (\alpha_{in}x_ir_iv_i)}{\sum_{n\in [N]} (1-\epsilon){g^n}(W_n)} = \frac{1 + \epsilon}{1 - \epsilon}F(\bx),
	\]
	which completes the proof.
	\endproof
	
	Now, suppose that the interval $[L_n, U_n]$ has already been partitioned into $K_n$ closed successive sub-intervals $[c^n_k, c^n_{k+1}]$ for $k \in [K_n-1]$, where $c^n_1 = L_n$ and $c^n_{K_n} = U_n$ (the partition procedure will be discussed in the next sub-section). Let $\gamma^{fn}_k = (f^n(c^n_{k+1}) - f^n(c^n_k))/(c^n_{k+1} - c^n_k)$ be the slope of $\widehat{f^n}(W_n)$ and $\gamma^{gn}_k = (g^n(c^n_{k+1}) - g^n(c^n_k))/(c^n_{k+1} - c^n_k)$ be the slope of $\widehat{g^n}(W_n)$ on the sub-interval $[c^n_k, c^n_{k+1}]$. Then, for every $W_n \in [c^n_k, c^n_{k+1}]$, the forms of $\widehat{f^n}(W_n)$ and $\widehat{g^n}(W_n)$ are given as follows:
	\begin{align}
		\label{eq:approx_f}
		&\widehat{f^n}(W_n) = f^n(c^n_k) + \gamma^{fn}_k(W_n - c^n_k)\\
		\label{eq:approx_g}
		&\widehat{g^n}(W_n) = g^n(c^n_k) + \gamma^{gn}_k(W_n - c^n_k). 
	\end{align}
	To use these approximations for our formulation, let us define real variables $z_{nk} \in [0, 1]$ for all $n\in[N], k\in [K_n]$ such that if $k^*\in [K_n-1]$ such that  $c^n_{k^*} \leq W_n < c^n_{k^*+1}$, we have $z_{nk} = 1\ \forall{k < k^*}$, $z_{nk} = (W_n - c^n_{k^*})/(c^n_{k^*+1} - c^n_{k^*})$ for $k = k^*$, and $z_{nk} = 0\ \forall{k > k^*}$. Then $\widehat{f^n}(W_n)$ and $\widehat{g^n}(W_n)$ can be rewritten as follows
	\begin{align}
		&\widehat{f^n}(W_n) =  f^n(L_n) + \sum_{k\in [K_n]}  \gamma^{fn}_{k}\Delta_{nk}{z_{nk}}\nonumber\\
		&\widehat{g^n}(W_n) =  g^n(L_n) + \sum_{k\in [K_n]}  \gamma^{gn}_{k}
		\Delta_{nk}{z_{nk}}, \nonumber
	\end{align}
	where $\Delta_{nk} = c^n_{k+1} - c^n_k$ is the length of the $k^{th}$ sub-intervals. Figure \ref{fig:approx_illustration} illustrates the functions $f^n(W_n)$ and $g^n(W_n)$ and the way that $[L_n, U_n]$ is partitioned into $K_n$ sub-intervals, where the $k^{th}$ sub-interval is associated with a real variable $z_{nk}$. 
	The objective function of the approximate problem \eqref{prob:approx} can be rewritten as follows
	\begin{align}
		\widehat{F}(\bx)  
		&=\frac{\sum_{i\in [m]}\sum_{n\in [N]}\left(f^n(L_n) + \sum_{k\in [K_n]}  \gamma^{fn}_{k}\Delta_{nk}z_{nk}\right) \left(\alpha_{in}x_ir_iv_i\right)}{\sum_{n\in [N]}\left(g^n(L_n) + \sum_{k\in [K_n]}  \gamma^{gn}_{k}\Delta_{nk}{z_{nk}}\right)} \nonumber \\
		&= \frac{\sum_{i\in[m],n\in [N]} a_{in} x_i + \sum_{i\in [m],n\in [N],k\in [K_n]}b_{ink}x_{i}z_{nk}
		}{c+ \sum_{n\in [N], k\in [K_n]}d_{nk} {z_{nk}}},
		\label{eq:approx_obj}
	\end{align}
	where $a_{in} = \alpha_{in}r_iv_i f^{n}(L_{n})$, $b_{ink} = \alpha_{in}r_iv_i \gamma^{fn}_{k}\Delta_{nk}$, $c =  \sum_{n\in [N]}g^n(L_n)$ and $d_{nk} = \Delta_{nk}\gamma^{gn}_{k}$ for all $i\in [m], n\in [N], k\in [K_n]$.
	Let $s_{ink} = x_iz_{nk}$ for all $i\in [m], n \in [N], k \in [K_n]$. We can linearize this bi-linear term by using the Glover’s linearization scheme \citep{glover1975improved} and reformulate the problem \eqref{prob:approx} as a linear-fractional program as follows
	\begin{align}
		\max_{\bx \in \cX} \ &  \left\{\widehat{F}(\bx)  = \frac{\sum_{i\in[m],n\in[N]} a_{in} x_i + \sum_{i\in [m],n\in [N],k\in [K_n]}b_{ink}s_{ink}
		}{c+ \sum_{n\in [N], k\in [K_n]}d_{nk} {z_{nk}} } \right\}\label{prob:approx1}\tag{\sf Assort-LFP}\\
		\text{subject to}\qquad & v_{0n} + \sum_{i\in S_n} \alpha_{in}x_iv_i = L_n+ \sum_{k\in [K_n]}\Delta_{nk}z_{nk},\;\forall n\in [N]\nonumber \\ 
		&y_{n,k+1} \leq y_{nk} \leq z_{nk},\;\forall n\in [N], k\in [K_n-1] \nonumber \\
		&z_{n,k+1}\leq y_{nk}, \forall n\in [N], k\in [K_n-2] \nonumber \\
		&s_{ink} \leq x_i, s_{ink} \leq z_{nk},\;\forall i\in [m], n\in [N], k \in [K_n] \nonumber\\
		&s_{ink} \geq x_i+z_{nk}-1,\;\forall i\in [m], n\in [N], k \in [K_n] \nonumber\\
		& y_{nk} \in \{0,1\},z_{nk}\in[0,1],  s_{ink}\in [0, 1], \forall i\in [m], n\in [N], k \in [K_n].\nonumber
	\end{align}
	Here, we use auxiliary binary variables $y_{nk}$ to ensure that the real variables $z_{nk}$ get proper values as their definitions. 
Theorem \ref{theo:approx_bound} below provides a performance guarantee when solving the approximate problem \eqref{prob:approx1}.
	\begin{theorem}
		\label{theo:approx_bound}
        For dissimilarity parameters $\sigma_n \in [0, 1], \forall n\in [N]$, suppose that $\widehat{f^n}(W_n)$ and $\widehat{g^n}(W_n)$ are constructed satisfying a given relative approximation error $\epsilon$.
		Let $\bx^*$ be an optimal solution to the original problem \eqref{prob:CNL-assort} and $\widehat{\bx}$ be an optimal solution to the approximate problem \eqref{prob:approx1}. Then we have
		\begin{align*}
			F(\bx^*) \leq \widehat{F}(\widehat{\bx}) \leq \frac{1 + \epsilon}{1 - \epsilon}F(\bx^*) \quad \text{and}
			\quad \frac{1 - \epsilon}{1 + \epsilon}F(\bx^*) \leq F(\widehat{\bx}) \leq F(\bx^*).
		\end{align*}
	\end{theorem}
	\proof{Proof.}
    Since the reformulation from \eqref{prob:approx} into \eqref{prob:approx1} is purely variable substitutions with an exact linearization from the bi-linear term $x_iz_{nk}$ into $s_{ink}$, the problem \eqref{prob:approx1} is equivalent to \eqref{prob:approx}.
	It follows that the results established in Lemma \ref{lemma:piecewise_function} for \eqref{prob:approx}  also hold for \eqref{prob:approx1}. In particular, we have
	\begin{equation}
	    F(\bx^*) \leq \widehat{F}(\bx^*) \quad \text{and} \quad \widehat{F}(\widehat{\bx}) \leq \frac{1 + \epsilon}{1 - \epsilon}F(\widehat{\bx}).
        \label{eq:temp}
	\end{equation}
    Since $\widehat{\bx}$ is an optimal solution to the approximate problem \eqref{prob:approx1}, we have
	$\widehat{F}(\bx^*) \leq \widehat{F}(\widehat{\bx})$. Combine with \eqref{eq:temp}, we get 
        \begin{equation}
	    F(\bx^*) \leq \widehat{F}(\widehat{\bx}) \leq \frac{1 + \epsilon}{1 - \epsilon}F(\widehat{\bx}).
     \label{eq:temp2}
	\end{equation}
	Moreover, since $x^*$ is optimal to 
	the problem \eqref{prob:CNL-assort}, we have $F(\widehat{\bx}) \leq F(\bx^*)$. Combine with \eqref{eq:temp2}, we get 
	\[
	F(\bx^*) \leq \widehat{F}(\widehat{\bx}) \leq \frac{1 + \epsilon}{1 - \epsilon}F(\widehat{\bx}) \leq \frac{1 + \epsilon}{1 - \epsilon}F(\bx^*),
	\]
	which directly implies the desired results.
	\endproof
	
	Theorem \ref{theo:approx_bound} indicates that one can obtain $\frac{1 - \epsilon}{1 + \epsilon}$-optimal solution to the assortment problem \eqref{prob:CNL-assort} by solving the linear-fractional program \eqref{prob:approx1}. In the next subsection, we will present the procedure to split $[L_n, U_n]$ into $K_n$ sub-intervals in such a way that  the properties in lemma \ref{lemma:piecewise_function} are satisfied and $K_n$ is minimized for any given accuracy level $\epsilon > 0$.
	
	\subsection{Discretization Procedure}
    \label{subsec:discretization}
    For any given $\epsilon > 0$, our goal is to partition $[L_n, U_n]$ into $K_n$ sub-intervals that satisfy the properties in Lemma \ref{lemma:piecewise_function}, i.e., on every sub-interval $[c^n_k, c^n_{k+1}]$,  $k\in[K_n]$, we have
    \begin{equation}
        \max_{c^n_k \leq t \leq c^n_{k+1}}\left\{\frac{\widehat{f^n}(t) - f^n(t)}{f^n(t)}\right\} \leq \epsilon
	\qquad
	\text{and}
	\qquad
	\max_{c^n_k \leq t \leq c^n_{k+1}}\left\{\frac{g^n(t) - \widehat{g^n}(t)}{g^n(t)}\right\} \leq \epsilon
    \label{eq:discrete_condition}
    \end{equation}
    where $\widehat{f^n}(\cdot)$ and $\widehat{g^n}(\cdot)$ take the forms as in \eqref{eq:approx_f} and \eqref{eq:approx_g}. Besides, we also want to minimize the number of sub-intervals $K_n$ to reduce the number of additional binary variables in the approximate problem. Without loss of generality, suppose that $k -1$ first sub-intervals ($k \geq 1$) have been constructed, i.e., $k$ points $c^n_1, c^n_2, \ldots c^n_{k}$ have been determined. We now present how to find the next point $c^n_{k+1} \in (c^n_{k}, U_n]$. 
 
	For the special case where $\sigma_n = 1$, we have $f^n(t) = 1$ and $g^n(t) = t$, which are linear functions. Thus we have $K_n = 1$, $c^n_1 = L_n$ and $c^n_{K_n + 1} = U_n$. For $\sigma_n < 1$, we first re-notate our piecewise approximation functions as follows. Let $\widehat{f^n}(c^n_k, u, t)$ be the approximation of $f^n(t)$ and $\widehat{g^n}(c^n_k, u, t)$ be the approximation of $g^n(t)$ at the point $t$ on a sub-interval $[c^n_k, u]$ ($c^n_k$ and $u$ are two breakpoints of the piece). 
    In the remaining parts of this paper, we will use $\widehat{f^n}(c^n_k, u, t)$ and $\widehat{g^n}(c^n_k, u, t)$ interchangeably with $\widehat{f^n}(t)$ and $\widehat{g^n}(t)$, i.e., we use the new notations when we care about the endpoints $c^n_k$ and $u$ of the linear piece, otherwise, the old notations are used with the same meanings.
    On $[c^n_k, u]$, let $\zeta^{fn}_k = (f^n(u) - f^n(c^n_k)) / (u - c^n_k)$ be the slope of $\widehat{f^n}(c^n_k, u, t)$ and $\zeta^{gn}_k = (g^n(u) - g^n(c^n_k)) / (u - c^n_k)$ be the slope of $\widehat{g^n}(c^n_k, u, t)$ in $[c^n_k,u]$.
    Then we have 
    \begin{align*}
        \widehat{f^n}(c^n_k, u, t) = f^n(c^n_k) + \zeta^{fn}_k(t - c^n_k),\\
        \widehat{g^n}(c^n_k, u, t) = g^n(c^n_k) + \zeta^{gn}_k(t - c^n_k).
    \end{align*}
    We have the following results for $\widehat{f^n}(c^n_k, u, t)$ and $\widehat{g^n}(c^n_k, u, t)$ with the proof given in Appendix \ref{appendix:proof_monotnone_of_approximators}.
    \begin{lemma}
    \label{lemma:monotnone_of_approximators}
        For $u \in (c^n_k, U_n]$,
        $\widehat{f^n}(c^n_k, u, t)$ is monotonically  increasing in $u$, $\widehat{g^n}(c^n_k, u, t)$ is monotonically  decreasing in $u$.
    \end{lemma}
    For every $t\in [c^n_k, u]$, let 
    \begin{align}
    \phi^{fn}(c^n_k, u, t) &= (\widehat{f^n}(c^n_k, u, t) - f^n(t)) / f^n(t)~ \text{ and }\nonumber\\
    \phi^{gn}(c^n_k, u, t) &= (g^n(t) - \widehat{g^n}(c^n_k, u, t)) / g^n(t) \nonumber
    \end{align}
     be the relative error of $\widehat{f^n}(c^n_k, u, t)$ w.r.t. $f^n(t)$, and  the relative error of $\widehat{g^n}(c^n_k, u, t)$ w.r.t. $g^n(t)$, respectively.  Let
    \begin{align}
     \Phi^{fn}(c^n_k, u) = \max_{c^n_k \leq t \leq u}\{\phi^{fn}(c^n_k, u, t)\}\text{ and }\nonumber \\
     \Phi^{gn}(c^n_k, u) = \max_{c^n_k \leq t \leq u}\{\phi^{gn}(c^n_k, u, t)\}\nonumber   
    \end{align}
    be the maximum relative errors of the approximators $\widehat{f^n}(c^n_k, u,t)$ and $\widehat{g^n}(c^n_k, u, t)$ on the sub-interval $[c^n_k, u]$, respectively. 
    Our goal is to find the point $u$ such that $\Phi^{fn}(c^n_k, u) \leq \epsilon$ and $\Phi^{gn}(c^n_k, u) \leq \epsilon$. 
    Since we also want to minimize the number of sub-intervals $K_n$, intuitively, the length of the sub-interval $\Delta_{nk} = |c^n_{k+1} - c^n_k|$ should be maximized. To this end,
	the value of $c^n_{k+1}$ can be determined by solving the following problem
	\begin{equation}
		\label{prob:bisection_discritization}
		c^n_{k+1} = 
		\max \left\{u\in (c^n_k, U_n] \ \Bigg|~ \Phi^{fn}(c^n_k, u)\leq \epsilon\ \text{and } \Phi^{gn}(c^n_k, u)\leq \epsilon\right\}.
	\end{equation}
	The following results indicate how we can calculate $\Phi^{fn}(c^n_k, u)$ and $\Phi^{gn}(c^n_k, u)$ for a given value of $u$ and solve  \eqref{prob:bisection_discritization} to find the next point $c^n_{k+1}$. Specifically, in Proposition \ref{proposition:unique_maximum_of_relative_error}, we show that $\phi^{fn}(c^n_k, u, t)$ and $\phi^{gn}(c^n_k, u, t)$ each have unique maximum on the sub-interval $[c^n_k, u]$, and thus we can directly calculate the values of $\Phi^{fn}(c^n_k, u)$ and $\Phi^{gn}(c^n_k, u)$. In Proposition \ref{lemma:monotone_of_max_error}, we show that $\Phi^{fn}(c^n_k, u)$ and $\Phi^{gn}(c^n_k, u)$ are monotonic in $u$, which allow us to efficiently solve  \eqref{prob:bisection_discritization} via binary search. 
 
	\begin{proposition}
        \label{proposition:unique_maximum_of_relative_error}
        On every sub-interval $[c^n_k, u] \subseteq [L_n, U_n]$, the relative error functions $\phi^{fn}(c^n_k, u, t)$ and $\phi^{gn}(c^n_k, u, t)$ have  unique maximums at
        \begin{align*}
        t^{fn}_k = \frac{(\sigma_n - 1)(f^n(c^n_k) - \zeta_k^{fn}c^n_k)}{\zeta_k^{fn}(2 - \sigma_n)}\nonumber \\ t^{gn}_k = \frac{\sigma_n(g^n(c^n_k) - \zeta_k^{gn}c^n_k)}{\zeta_k^{gn}(1 - \sigma_n)},    
        \end{align*}
        respectively.
	\end{proposition}
	\proof{Proof.}
        For $\sigma_n < 1$, the lemma can be verified simply using the derivatives of $\phi^{fn}(c^n_k, u, t)$ and $\phi^{gn}(c^n_k, u, t)$. Let take the first partial derivative of $\phi^{fn}(c^n_k, u, t)$ with respect to $t$ as
	\begin{align*}
		\frac{\partial}{\partial t}\phi^{f^n}(c^n_k, u, t) &= \frac{\partial}{\partial t}\left(\frac{\widehat{f^n}(c^n_k, u, t)}{f^n(t)} - 1\right) = \frac{\partial}{\partial t}\left(\widehat{f^n}(c^n_k, u, t) t^{1 - \sigma_n} - 1\right) \\ &= \zeta_k^{fn} t^{1 - \sigma_n} + (1 - \sigma_n)\widehat{f^n}(c^n_k, u, t) t ^{-\sigma_n}.
	\end{align*}
	Then, we can solve the equation $\frac{\partial}{\partial t}\phi^{f^n}(c^n_k, u, t) = 0$ and obtain a unique solution
    \begin{equation*}
        t^{fn}_k = \frac{(\sigma_n - 1)(f^n(c^n_k) - \zeta_k^{fn}c^n_k)}{\zeta_k^{fn}(2 - \sigma_n)}.
    \end{equation*}
    Since $c^n_k$ and $c^n_k$ are breakpoints in our approximation, we have $\phi^{fn}(c^n_k, u, c^n_k) = \phi^{fn}(c^n_k, u, c^n_k) = 0$. Furthermore, since $\widehat{f^n}(c^n_k, u, t) \geq f^n(t)\ \forall t\in [L_n, U_n]$, we have $\phi^{fn}(c^n_k, u, t) > 0 \ \forall t\in (c^n_k, u)$. It follows that $t^{fn}_k \in [c^n_k, u]$ and $\phi^{fn}(c^n_k, u, t)$ has a unique maximum at this point.


	For the case of $\phi^{gn}(c^n_k, u, t)$, we have
	\begin{align*}
		\frac{\partial}{\partial t}\phi^{g^n}(c^n_k, u, t) &= \frac{\partial}{\partial t}\left(1 - \frac{\widehat{g^n}(c^n_k, u, t)}{g^n(t)}\right) = \frac{\partial}{\partial t}\left(1 - \widehat{g^n}(c^n_k, u, t)t^{-\sigma_n}\right)\\ &= \sigma_n\widehat{g^n}(c^n_k, u, t)t^{-\sigma_n - 1} - \zeta_k^{gn}t^{-\sigma_n}.
	\end{align*}
	Solving the equation $\frac{\partial}{\partial t}\phi^{g^n}(c^n_k, u, t) = 0$, we obtain an unique solution $$t^{gn}_k = \frac{\sigma_n(g^n(c^n_k) - \zeta_k^{gn}c^n_k)}{\zeta_k^{gn}(1 - \sigma_n)}.$$ 
 Moreover, $\phi^{g^n}(c^n_k, u, c^n_k) = \phi^{g^n}(c^n_k, u, u) = 0$  and $\phi^{g^n}(c^n_k, u, t) >0$ for all $t\in (c^n_k,u)$. So,  $\phi^{gn}(c^n_k, u, t)$ has a unique maximum at $t^{fn}_k \in [c^n_k, u]$,
	which completes the proof.
	\endproof
	
	\begin{proposition}
        \label{lemma:monotone_of_max_error}
		The maximum relative errors of the approximators $\widehat{f^n}(c^n_k, u, t)$ and $\widehat{g^n}(c^n_k, u, t)$ on the sub-interval $[c^n_k, u] \subseteq [L_n, U_n]$ are monotonically non-decreasing in $u$, i.e., 
		 $\Phi^{fn}(c^n_k, u)$ and $\Phi^{gn}(c^n_k, u)$ are monotonically non-decreasing in $u$. 
	\end{proposition}

	With the above results,  Problem \ref{prob:bisection_discritization} can be solved via a bisection procedure. Starting from $c^n_1 = L_n$, we can determine all the breakpoints by repeatedly solving \eqref{prob:bisection_discritization} until reaching the upper bound $U_n$. We describe the bisection in Algorithm \ref{algo:DA} below, where $\delta > 0$ is a predefined tolerance used to stop the bisection. We assume that $\delta$ is significantly smaller than the accuracy level $\epsilon$, and thus can be removed from the overall performance guarantee. We further analyze the properties of the sub-intervals generated by Algorithm \ref{algo:DA} in the next section.


\begin{algorithm}[htb]
    \caption{Discretization} 
    \label{algo:DA}
    \SetKwRepeat{Do}{do}{until}
    For any $n\in [N]$:
     \begin{itemize}
		\item[\textit{Step 1.}] Start by $k := 1$ and $c^n_1 := L_n$.
		\item[\textit{Step 2.}] $a := c^n_k$, $b := U_n$.
		\item[\textit{Step 3.}] $u := (a + b) / 2$. If $\Phi^{fn}(a, u) \leq \epsilon$ and $\Phi^{gn}(a, u) \leq \epsilon$ then $a := u$, else $b := u$.
		\item[\textit{Step 4.}] If $b - a \leq \delta$ then $c^n_{k+1} := u$ and go to \textit{Step 5}, else repeat \textit{Step 3}.
		\item[\textit{Step 5.}] If $U_n - c^n_{k+1} \leq \delta$ then $c^n_{k+1} := U_n$ and stop, else $k := k+1$ and repeat \textit{Step 2}.
	\end{itemize}
\end{algorithm}

\subsection{Optimality and Bounds}
\label{subsec:K_opt_bounds}
Theorem \ref{theorem:min_K} bellow shows that the procedure described in the above section is, in fact, optimal, in the sense that there is no set of breakpoints that has a small number of breakpoints while achieving better or the same relative errors. 
\begin{theorem}
        \label{theorem:min_K}
        For any given $\epsilon > 0$, the above partitioning procedure generates the minimum number of sub-intervals required to achieve the performance guarantee of $\frac{1-\epsilon}{1 + \epsilon}$. In other words, there does not exist any other way to partition $[L_n, U_n]$ into $K'_n < K_n$ sub-intervals that offer the same performance guarantee.
	\end{theorem}
	\proof{Proof.}
	Assuming that  there exists a set of points $\{d^n_1, d^n_2, \ldots, d^n_{K'_n+1}\}$ that partitions the interval $[L_n, U_n]$ into $K'_n$ sub-intervals and satisfies the performance guarantee of $\frac{1 - \epsilon}{1 + \epsilon}$, i.e., on every sub-interval $k\in[K'_n]$, we have $\Phi^{fn}(d^n_k, d^n_{k+1}) \leq \epsilon$ and $\Phi^{gn}(d^n_k, d^n_{k+1}) \leq \epsilon$. 
 Let us denote $\{c^n_1, c^n_2, \ldots, c^n_{K_n+1}\}$ as the set of points generated by Algorithm \ref{algo:DA}.
 We will prove $K_n \leq K'_n$ by showing that, for any $k\in [K_n - 1]$,  if $d^n_k \leq c^n_k$ then $d^n_{k+1} \leq c^n_{k+1}$. 
	
	Suppose, on the contrary, that there exist $k \in [K_n - 1]$ such that $d^n_k \leq c^n_k$ and $d^n_{k+1} \geq c^n_{k+1}$, $k \leq K_n - 1$. Since $f^n(t)$ is convex and decreasing, $\widehat{f^n}(c^n_k, c^n_{k+1}, t)$ is decreasing in $c^n_k$ and increasing in $c^n_{k+1}$. This means, for every $t \in [c^n_k, c^n_{k+1}]$, $\widehat{f^n}(d^n_k, d^n_{k+1}, t) > \widehat{f^n}(c^n_k, c^n_{k+1}, t)$. Similarly, since $g^n(t)$ is concave and increasing, $\widehat{g^n}(c^n_k, c^n_{k+1}, t)$ is increasing in $c^n_k$ and decreasing in $c^n_{k+1}$, which implies that for every $t \in [c^n_k, c^n_{k+1}]$, $\widehat{g^n}(d^n_k, d^n_{k+1}, t) < \widehat{g^n}(c^n_k, c^n_{k+1}, t)$. Moreover, since $[c^n_k, c^n_{k+1}]$ is a sub-interval obtained from our procedure, there exists a point $t^*\in [c^n_k, c^n_{k+1}]$ such that 
	\[
	\frac{\widehat{f^n}(c^n_k, c^n_{k+1}, t^*) - f^n(t^*)}{f^n(t^*)} \ =\  \epsilon
	\qquad
	\text{or}
	\qquad
	\frac{g^n(t^*) - \widehat{g^n}(c^n_k, c^n_{k+1}, t^*)}{g^n(t^*)} \ =\  \epsilon.
	\]
	It follows that there exists a point $t^*\in [c^n_k, c^n_{k+1}]$ such that 
	\[
	\Phi^{fn}(d^n_k, d^n_{k+1}) \ \geq\ \frac{\widehat{f^n}(d^n_k, d^n_{k+1}, t^*) - f^n(t^*)}{f^n(t^*)} \ >\ 
	\frac{\widehat{f^n}(c^n_k, c^n_{k+1}, t^*) - f^n(t^*)}{f^n(t^*)} \ =\  \epsilon
	\]
	\noindent or
	\[
	\Phi^{gn}(d^n_k, d^n_{k+1}) \ \geq\ \frac{g^n(t^*) - \widehat{g^n}(d^n_k, d^n_{k+1}, t^*)}{g^n(t^*)} \ >\  \frac{g^n(t^*) - \widehat{g^n}(c^n_k, c^n_{k+1}, t^*)}{g^n(t^*)} \ =\  \epsilon,
	\]
	which contradicts to the assumption that $\Phi^{fn}(d^n_k, d^n_{k+1}) \leq \epsilon$ and $\Phi^{gn}(d^n_k, d^n_{k+1}) \leq \epsilon, \forall{k\in[K'_n]}$. Therefore, if $d^n_k \leq c^n_k$ then $d^n_{k+1} \leq c^n_{k+1}$ holds. Since $c^n_1 = d^n_1 = L_n$, it follows that $d^n_{K_n+1} \leq c^n_{K_n+1}$, or $d^n_{K_n+1} \leq U_n$, which indicates that $K_n \leq K'_n$, as desired.
	\endproof

To give a better sense of how the number of sub-intervals returned by Algorithm \ref{algo:DA} is affected by the model parameters such as $\sigma_n$ or $U_n,L_n$, we further establish lower and upper bounds for it. To this end, we first  bound $\Phi^{fn}(c^n_k, u)$ and $\Phi^{gn}(c^n_k, u)$ in the following lemma. 
\begin{lemma}\label{lm:bound-phi}
   Given $n\in [N]$ and scalars $c, u$ such that $L_n \leq c\leq u\leq U_n$, we can bound  $\Phi^{fn}(c, u)$ and $\Phi^{gn}(c, u)$ as
   \begin{align}
        \frac{2^{\sigma_n-2}(u^{\sigma_n-1} + c^{\sigma_n-1})}{(u+c)^{\sigma_n-1}}-1&\leq \Phi^{fn}(c, u) \leq \left(\frac{c}{u}\right)^{\sigma_n-1} - 1 \nonumber \\
        1-\frac{2^{\sigma_n-1}(u^{\sigma_n} + c^{\sigma_n})}{(u+c)^{\sigma_n}}&\leq \Phi^{gn}(c, u) \leq 1 - \left(\frac{c}{u}\right)^{\sigma_n}.  \nonumber
   \end{align}
\end{lemma}
It can be seen that the bounds presented in Lemma \ref{lm:bound-phi} are not trivial; they are all positive and will converge to 0 as $c$ approaches $u$. To bound $K_n$, let us consider the approximation of each function $f^n(\cdot)$ and $g^n(\cdot)$ separately. That is, let $K^{f}_n$ be the number of sub-intervals required to approximate the function $f^n(\cdot)$ separately such that $\Phi^{fn}(c^{fn}_{k}, c^{fn}_{k+1}) = \epsilon,\ \forall k < K^{f}_n$ and $\Phi^{fn}(c^{fn}_{k}, c^{fn}_{k+1}) \leq \epsilon$ for $k = K^{f}_n$. Similarly, let $K^{g}_n$ be the number of sub-intervals required to approximate the function $g^n(\cdot)$ separately such that $\Phi^{gn}(c^{gn}_{k}, c^{gn}_{k+1}) = \epsilon,\ \forall k < K^{g}_n$ and $\Phi^{gn}(c^{gn}_{k}, c^{gn}_{k+1}) \leq \epsilon$ for $k = K^{g}_n$. The following lemma provides lower and upper bounds for $K^{f}_n$ and $K^{g}_n$.  
\begin{lemma}\label{lm:lm5}
    For any $\epsilon \in (0,1-2^{\sigma_n-1}]$
    \begin{align}
      \frac{\ln(U_n/L_n)}{\ln t^f(\epsilon)}  &\leq K^{f}_n\leq \frac{(1-\sigma_n)\ln (U_n/L_n)}{\ln 
 (\epsilon+1)} +1   \nonumber \\
        \frac{\ln(U_n/L_n)}{\ln t^g(\epsilon)}  &\leq K^{g
    }_n\leq  \frac{-\sigma_n\ln (U_n/L_n)}{\ln(1 - \epsilon)} +1   \nonumber
    \end{align}
    where $t^f(\epsilon)$ and $t^g(\epsilon)$ are unique solutions to the following equations (respectively):
    \begin{align}
       \frac{t^{\sigma_n-1}+1}{(t+1)^{\sigma_n-1}} &= 2^{2-\sigma_n}(\epsilon+1)\label{eq:lm5-eq11} \\
        \frac{t^{\sigma_n}+1}{(t+1)^{\sigma_n}} &= 2^{1-\sigma_n}(1-\epsilon). \label{eq:lm5-eq12}     
    \end{align}
    Moreover,  $t^f(\epsilon)> 1$, $t^f(\epsilon)> 1$, and they are monotonically decreasing as $\epsilon$ decreases,  and $\lim_{\epsilon\rightarrow 0} \max\{t^f(\epsilon),t^g(\epsilon)\} = 1$.
\end{lemma}
It can be revealed from the proof of Lemma \ref{lm:lm5} that the functions $\delta(t) =  \ln({t^{\sigma_n-1}+1})-\ln ({(t+1)^{\sigma_n-1}})$ and  $\rho(t) = \ln({t^{\sigma_n}+1})-\ln({(t+1)^{\sigma_n}})$ are monotonic in $t$, thus $t^f(\epsilon)$ and $t^g(\epsilon)$ can be computed efficiently via bisection. The limit  $t^f(\epsilon)> 1$, $t^f(\epsilon)> 1$ implies that both the lower and upper bounds will go to infinity as $\epsilon\rightarrow 0$.

At this point, we are  ready to present bounds for the number of sub-intervals returned by Algorithm \ref{algo:DA}.

\begin{theorem}    \label{theorem:partitioning_convergence}
    For any given $\epsilon$ such that $0< \epsilon\leq \min_{n\in [N]}\{1- 2^{\sigma_n-1}\} $, the above partitioning procedure converges after a finite number of iterations and the number of generated sub-intervals $K_n$, for any $n\in [N]$, required to achieve the performance guarantee of $\frac{1 - \epsilon}{1 + \epsilon}$ is bounded as
    \begin{align}
       \ln \left(\frac{U_n}{L_n}\right) \max\left\{\frac{1}{\ln t^f(\epsilon)};~\frac{1}{\ln t^g(\epsilon)}\right\} \leq K_n \leq \ln \left(\frac{U_n}{L_n}\right) \left(\frac{1-\sigma_n}{\ln 
 (\epsilon+1)} - 
   \frac{\sigma_n}{\ln(1 - \epsilon)}\right) +1.\nonumber
    \end{align}
\end{theorem}
The proof of Theorem \ref{theorem:partitioning_convergence} is given in Appendix \ref{appendix:proof_convergence}. From Lemma \ref{lm:lm5}, it can be seen that the lower bounds are always positive and will go to infinity as $\epsilon$ approaches 0. We note that the condition $\epsilon \in (0,1-2^{\sigma_n-1}]$ is only required for the establishment of the lower bounds. This condition does not hold, for any $\epsilon>0$, if $\sigma_n=1$. In this case, $f^n(t)$ and $g^n(t)$ are linear in $t$ and we will only need one sub-interval to achieve the performance guarantee, i.e., $K_n= 1$. Moreover, since $\lim_{\epsilon \rightarrow 0} \{\epsilon/\ln (1+\epsilon)\} = \lim_{\epsilon \rightarrow 0} \{\epsilon/\ln (1-\epsilon)\} = -1$, the upper bound of $K_n$ stated in Theorem \ref{theorem:partitioning_convergence} implies that Algorithm \ref{algo:DA} stops after $\cO(1/\epsilon)$ iterations.


The lower bounds presented in Theorem \ref{theorem:partitioning_convergence} are somewhat ambiguous because they rely on $t^{f}(\epsilon)$ and $t^{g}(\epsilon)$, which are expressed as solutions to certain nonlinear equations. To provide a more precise but less strong lower bound for $K_n$, we introduce Corollary \ref{coro:bounds-Kn} below. The bounds can be deduced from the proof of Lemma \ref{lm:lm5}.  
\begin{corollary}\label{coro:bounds-Kn}
    The number of sub-intervals returned by Algorithm \ref{algo:DA} can be bounded from below as
    \[
    K_n\geq    \ln \left(\frac{U_n}{L_n}\right)\max\left\{\frac{1}{\ln t^{f*}};~\frac{1}{\ln t^{g*}}\right\}
    \]
    where 
    \begin{align*}
        t^{f*} &= \exp\left(\frac{\ln(1+\epsilon)+(2-\sigma_n)\ln 2}{1-\sigma_n}\right) - 1; ~ \text{ 
 and }
        t^{g*} = \left(\frac{1}{2^{1-\sigma_n}(1-\epsilon)-1}\right)^{1/\sigma_n}.
    \end{align*}
\end{corollary}

It can be seen that $t^{f*}> t^{f}(\epsilon)$ and $t^{f*}> t^{f}(\epsilon)$, so the lower bounds stated in Corollary \ref{coro:bounds-Kn} are lower (so weaker) than those  in Theorem \ref{theorem:partitioning_convergence}.  In addition, $t^{f*}$ and $t^{g*}$ will not approach 1 when $\epsilon$ tends to zero (except when $\sigma_n$ approaches 1), thus the lower bounds in Corollary \ref{coro:bounds-Kn} will not go to infinity as $\epsilon\rightarrow 0$. 


	In Tables \ref{tab:k_bisection} below, we examine the actual number of sub-intervals $K_n$ generated by Algorithm \ref{algo:DA} by varying the model parameters $U_n$, $\sigma_n$, and $\epsilon$. Here, $U_n$ is the upper bound of the total preference weights of nest $S_n$, which is related to the number of products in $S_n$. The value of $\sigma_n$ indicates the correlation level between the products, and $\epsilon$ decides the approximation accuracy. To see how these parameters affect the  number of generated sub-intervals, $\sigma_n$ is varied in $\{0.2, 0.3, 0.5, 0.7, 0.9\}$, $U_n$ is varied over $\{5, 10, 15\}$, $L_n$ is set to 1,  and $\epsilon$ takes value from \{0.1, 0.05, 0.01, 0.005, 0.001\}. These values of $\epsilon$ result in the performance guarantees of \{81.8\%, 90.5\%, 98.0\%, 99.0\%, 99.8\%\}, respectively. From the table, we can see that all three parameters have their influences on the number of generated sub-intervals. In particular, $K_n$ increases as $\epsilon$ or $\sigma_n$ decreases, while the opposite happens for $U_n$. These results indicate that smoother division is required for higher approximation accuracy, or if the correlation between products in the nest is high. 
    However, we can see that the value of $K_n$ is not necessarily too large to achieve a high accuracy level. For instance, 6 segments are sufficient for the approximation guarantee of 90\%, while the number is 17 segments for 99\%, which is almost the exact level.
	\input{tables/k_bisection.tex}
	
	\section{Solving the Approximate Problem}
	\label{sec:solution_methods}
        In the previous section, we showed that for any given accuracy level $\epsilon > 0$, we can obtain $\frac{1 - \epsilon}{1 + \epsilon}$-optimal solution of the assortment problem \eqref{prob:CNL-assort} by solving the approximate problem \eqref{prob:approx1}. This section presents two methods to solve the problem \eqref{prob:approx1}. The first method (Subsection \ref{subsec:bisection}) is based on a parametric approach, where a sequence of feasibility linear programs is solved to gradually refine the solution to a near-optimum. In the second method (Subsection \ref{subsec:milp_assort}), we make a transformation to convert the approximate problem \eqref{prob:approx1} into a MILP with that the problem \eqref{prob:approx1} can be solved by using such off-the-shelf linear solvers like IBM's CPLEX or Gurobi. The second approach also enables us to extend the approximation method to tackle assortment optimization under a mixture of CNL models.   
        

	\subsection{Bisection Method}
	\label{subsec:bisection}
 Since the approximate problem \eqref{prob:approx1} is a fractional program, it is convenient (and also popular) to use the Dinkelbach transform \citep{dinkelbach1967nonlinear} to simplify the fractional structure. That is, we can reformulate \eqref{prob:approx1} as
	\begin{align*}
		\max \quad& \lambda \label{prob:assort_bisection}\tag{Assort-BIS}\\
		\text{subject to}\quad & \exists (\bx,\by,\bz,\bs)\in \Theta~\text{ s.t. }{\sum_{\substack {i\in[m]\\n\in[N]}} a_{in} x_i + \sum_{\substack{i\in [m]\\n\in [N],k\in [K_n]}}b_{ink}s_{ink}
		}\geq \lambda({c+ \sum_{\substack{n\in [N]\\ k\in [K_n]}}d_{nk} {z_{nk}} }),
	\end{align*}
	where $\Theta$ is the feasible set of \eqref{prob:approx1}.
	Since $\lambda$ is just a scalar, we can perform a binary search to find a near-optimal solution to \eqref{prob:approx1}. Let $L_{\lambda}$ and $U_{\lambda}$ be the lower bound and upper bound of $\lambda$. We can obtain the value of $L_{\lambda}$ by running a greedy algorithm that iteratively adds the next highest-revenue product to the assortment as long as the constraints are satisfied. For $U_{\lambda}$, we can simply set $U_{\lambda} = \max_{i\in[m]}\{r_i\}$. Given a relative tolerance $\delta > 0$, the following  bisection method can be used to solve  \eqref{prob:assort_bisection}.

    \begin{framed}{\textbf{\underline{Bisection}}}        
	\begin{itemize}
		\item[\textit{Step 1.}] Set $\lambda := (U_{\lambda} + L_{\lambda}) / 2$.
		\item[\textit{Step 2.}] Solve the following feasibility program: $\lambda$ is feasible iff there exists {$ (\bx,\by,\bz,\bs)\in \Theta$ such that}  
        \begin{align*}
			{\sum_{\substack {i\in[m]\\n\in[N]}} a_{in} x_i + \sum_{\substack{i\in [m]\\n\in [N],k\in [K_n]}}b_{ink}s_{ink}
			} -  \lambda({c+ \sum_{\substack{n\in [N]\\ k\in [K_n]}}d_{nk} {z_{nk}} }) \geq 0.
		\end{align*}
		\item[\textit{Step 3.}] If feasible, then $L_{\lambda} \leftarrow \lambda$, else $U_{\lambda}\leftarrow\lambda$.
		\item[\textit{Step 4.}] If $(U_{\lambda} - L_{\lambda}) / L_{\lambda} \leq \delta$, then return the feasible assortment, else repeat \textit{Step 1}.
	\end{itemize}
    \end{framed} 
	This procedure can provide $\delta$-approximation solution to the problem \eqref{prob:assort_bisection} in $\cO(\tau\log (1/\delta))$, where $\tau$ is the computing time to solve the feasibility problem in \textit{Step 2}.

Theorem \ref{theorem:bisection_guarantee} below gives a performance guarantee of $\frac{(1 - \delta)(1 - \epsilon)}{1 + \epsilon}$ for any solution returned by the above bisection procedure. 
 
	\begin{theorem}
        \label{theorem:bisection_guarantee}
		Let $\widetilde{\bx}$ be a feasible solution to problem \eqref{prob:assort_bisection} returned by the bisection method with a relative tolerance of $\delta$. Then, the expected revenue obtained by offering the assortment $\widetilde{\bx}$ differs from the optimal expected revenue $F(\bx^*)$ by no more than a factor of $\frac{(1 - \delta)(1 - \epsilon)}{1 + \epsilon}$.
	\end{theorem}
	\proof{Proof.}
	Let $\widehat{L}_{\lambda}$ and $\widehat{U}_{\lambda}$ be the bounds of parameter $\lambda$ in the last iteration of the bisection method, i.e. $(\widehat{U}_{\lambda} - \widehat{L}_{\lambda}) / \widehat{L}_{\lambda} \leq \delta$. Let $\widehat{\bx}$ be the optimal solution of the approximate problem \eqref{prob:approx1}. Then we have
	\[
	\widehat{L}_{\lambda} \leq \widehat{F}(\widetilde{\bx}) \leq \widehat{F}(\widehat{\bx}) \leq \widehat{U}_{\lambda}.
	\]
	It follows that
	\[
	0 \leq \frac{\widehat{F}(\widehat{\bx}) - \widehat{F}(\widetilde{\bx})}{\widehat{F}(\widehat{\bx})} \leq \frac{\widehat{U}_{\lambda} - \widehat{L}_{\lambda}}{\widehat{L}_{\lambda}} \leq \delta,
	\]
	which equivalent to 
	\begin{equation}
	    (1 - \delta)\widehat{F}(\widehat{\bx}) \leq \widehat{F}(\widetilde{\bx}) \leq \widehat{F}(\widehat{\bx}).
        \label{eq:tmp1_theo4_proof}
	\end{equation}
        From Theorem \ref{theo:approx_bound}, we have 
        \begin{equation}
            F(\bx^*) \leq \widehat{F}(\widehat{\bx}) \leq \frac{1 + \epsilon}{1 - \epsilon}F(\bx^*).
            \label{eq:tmp2_theo4_proof}
        \end{equation}
        Combining \eqref{eq:tmp1_theo4_proof} and \eqref{eq:tmp2_theo4_proof}, we get
        \begin{equation}
            (1 - \delta)F(\bx^*) \leq \widehat{F}(\widetilde{\bx}) \leq \frac{1 + \epsilon}{1 - \epsilon}F(\bx^*).
            \label{eq:tmp3_theo4_proof}
        \end{equation}
	From Lemma \ref{lemma:piecewise_function}, we have
        \begin{equation}
           \widehat{F}(\widetilde{\bx}) \leq \frac{1 + \epsilon}{1 - \epsilon}F(\widetilde{\bx}).
            \label{eq:tmp4_theo4_proof}
        \end{equation}
    Combining \eqref{eq:tmp3_theo4_proof} and \eqref{eq:tmp4_theo4_proof}, we get
	\[
	(1 - \delta)F(\bx^*) \leq \frac{1 + \epsilon}{1 - \epsilon}F(\widetilde{\bx}) \leq \frac{1 + \epsilon}{1 - \epsilon}F(\bx^*),
	\]
	or we can rewrite as $\frac{(1 - \delta)(1 - \epsilon)}{1 + \epsilon}F(\bx^*) \leq F(\widetilde{\bx}) \leq F(\bx^*),$
	which implies the desired result. 
	\endproof
	
	\subsection{Mixed-integer Linear Program Reformulation}
	\label{subsec:milp_assort}
	In this section, we transform  \eqref{prob:approx1} in Section \ref{sec:approximation} into a MILP using the Charnes-Cooper transformation \citep{charnes1973explicit} and the Glover’s linearization scheme \citep{glover1975improved}. First, let $w = 1 / (c+\sum_{n\in [N]}\sum_{k\in [K_n]} d_{nk}z_{nk})$. Since $c$ and $d_{nk}$ are all positive, $w$ is positive. Let $L^w$ and $U^w$ be the lower and upper bound of $w$. These bounds can simply be calculated by $L^w = 1 / (c+\sum_{n\in [N]}\sum_{k\in [K_n]} d_{nk})$ and $U^w = 1 / c$. Then, we can rewrite the program as follows
	\begin{align}
		\max_{\bx \in \cX} \ &  \left\{\widehat{F}(\bx)  = \sum_{i\in[m]}\sum_{n\in[N]} \left(a_{in} wx_i + \sum_{k\in[K_n]}b_{ink}ws_{ink}\right) \right\}\nonumber\\
		\text{subject to}\qquad & cw+ \sum_{n\in[N]}\sum_{k\in[K_n]}d_{nk}wz_{nk} = 1 \nonumber\\
		& \sum_{i\in S_n} \alpha_{in}v_iwx_i = (L_n-v_{0n})w+ \sum_{k\in [K_n]}\Delta_{nk}wz_{nk},\;\forall n\in [N]\nonumber \\ 
		&wy_{n,k+1} \leq wy_{nk} \leq wz_{nk},\;\forall n\in [N], k\in [K_n-1] \nonumber \\
		&wz_{n,k+1}\leq wy_{nk}, \forall n\in [N], k\in [K_n-2] \nonumber \\
		&ws_{ink} \leq wx_i, ws_{ink} \leq wz_{nk},\;\forall i\in [m], n\in [N], k \in [K_n] \nonumber\\
		&ws_{ink} \geq wx_i+wz_{nk}-w,\;\forall i\in [m], n\in [N], k \in [K_n] \nonumber\\
		& y_{nk} \in \{0,1\},z_{nk}\in[0,1],  s_{ink}\in [0, 1], \forall i\in [m], n\in [N], k \in [K_n].\nonumber
	\end{align}
	Now let $w^x_i = w x_i$, $w^y_{nk} = w y_{nk}$, $w^s_{ink} =  w s_{ink}$ and $w^z_{nk} = w z_{nk}$, for all $i\in [m], n\in [N], k\in [K_n]$. Note that only continuous variables are transformed automatically by the Charnes-Cooper transformation, i.e., $w s_{ink}$ into $w^s_{ink}$ and $w z_{nk}$ into $w^z_{nk}$. Therefore, to linearize the bi-linear terms $w x_i$ and $w y_{nk}$, we integrate the transformation with the Glover’s linearization scheme and obtain a MILP as follows:
	\begin{align}
		\max_{\bx \in \cX} \ &  \left\{\widehat{F}(\bx)  = \sum_{i\in[m]}\sum_{n\in[N]} \left(a_{in} w^x_i + \sum_{k\in[K_n]}b_{ink}w^s_{ink} \right)\right\}\label{prob:milp}\tag{Assort-MILP}\\
		\text{subject to}\qquad & cw+ \sum_{n\in[N]}\sum_{k\in[K_n]}d_{nk}w^z_{nk} = 1 \nonumber\\
		& \sum_{i\in S_n} \alpha_{in}v_iw^x_i = (L_n-v_{0n})w+ \sum_{k\in [K_n]}\Delta_{nk}w^z_{nk},\;\forall n\in [N]\nonumber \\ 
		&w^y_{n,k+1} \leq w^y_{nk} \leq w^z_{nk},\;\forall n\in [N], k\in [K_n-1] \nonumber \\
		&w^z_{n,k+1}\leq w^y_{nk}, \forall n\in [N], k\in [K_n-2] \nonumber \\
		&w^s_{ink} \leq w^x_i, w^s_{ink} \leq w^z_{nk},\;\forall i\in [m], n\in [N], k \in [K_n] \nonumber\\
		&w^s_{ink} \geq w^x_i+w^z_{nk}-w,\;\forall i\in [m], n\in [N], k \in [K_n] \nonumber\\
		&w - U^w (1-x_i) \leq w^x_i \leq w,\ w^x_i \leq U^w x_i, \forall{i\in[m]}\label{ctr:milp_glover_x} \\
		&w - U^w (1-y_{nk}) \leq w^y_{nk} \leq w,\ w^y_{nk} \leq U^w y_{nk}, \forall{n\in[N], k\in[K_n]}\label{ctr:milp_glover_y} \\
		&w\in \bbR_+, y_{nk} \in \{0,1\}, w^x_i\in \bbR, w^y_{nk} \in \bbR, w^z_{nk} \in \bbR, w^s_{ink}\in \bbR,\nonumber\\
		&\qquad\qquad\qquad\qquad\qquad\qquad\qquad \forall{i\in[m], n\in[N], k\in[K_n]}.\nonumber
	\end{align}
    In the above program, auxiliary constraints \eqref{ctr:milp_glover_x} and \eqref{ctr:milp_glover_y} result from the linearization of the bilinear terms $w x_i$ and $w y_{nk}$. Compared to the program \eqref{prob:approx1}, the transformation and linearization increase the model size by additional $\sum_{n\in[N]}K_n + m + 1$ real variables and $\sum_{n\in[N]}3K_n + 3m + 1$ constraints. However, since the difficulty of an integer program often depends much more on the number of integer variables rather than the number of continuous variables, we expect that the above transformation does not significantly increase the complexity of the problem compared to \eqref{prob:approx1}. Program 
     \eqref{prob:milp} can be handled by an available commercial MILP solver such as CPLEX or Gurobi. In a practical setting, the solver is often carried out with a time limit. Suppose that, within a given time budget $\tau$, the solver takes \eqref{prob:milp} and returns a solution $\widetilde{\bx}$ with an optimality gap of $\delta$, which is qualified as the relative gap between the solution $\widetilde{\bx}$ and the best bound. We then have the following performance guarantee for $\widetilde{\bx}$ that can be verified in the same way as for Theorem \ref{theorem:bisection_guarantee}.
    \begin{corollary}
    \label{corol:milp_guarantee}
        Let $\widetilde{\bx}$ be a solution returned by a solver after a given time budget $\tau$ and $\delta$ be the optimality gap from $\widetilde{\bx}$ to the best bound. Then, the obtained expected revenue $F(\widetilde{\bx})$ differs from the optimal expected revenue $F(\bx^*)$ by no more than a factor of $\frac{(1 - \delta)(1 - \epsilon)}{1 + \epsilon}$.
    \end{corollary}
 
	\section{Joint Assortment and Price Optimization}
	\label{sec:formulation-jap}
	This section considers the joint assortment and price optimization problem under the CNL model, where the retailer has to select a subset of products and set their prices to offer to the customer so that the expected revenue is maximized. There are $m$ products that has been assigned to $n$ nests $S_1, S_2, \ldots, S_n$. We assume that for each product $i \in [m]$, the decision maker can offer a price selected from a finite set $P_i = \{p_{i1}, p_{i2}, \ldots, p_{iL}\}$ of $L$ price levels. Here, we assume that the price decisions are discrete. This setting is not restrictive in practice, as a seller usually only takes into account  a finite number of possible prices, such as those that are expressed in the form of ``\textit{x.99}''.
 
 Let $v_{il}$ be the preference weight of the product $i$ when its price is set to $p_{il}$. Here we assume that the relationship between $p_{il}$ and $v_{il}$ can be completely arbitrary. Our approach is to convert the joint assortment-price optimization into an assortment optimization problem with an extended set of products. In particular, we consider each product $i$ and its $L$ possible prices as $L$ alternatives where each of them is associated with one price level. Let $\bx \in \{0, 1\}^{m\times L}$ be a binary vector representing  a combination of assortment and price, where $x_{il} = 1$ if and only if the product $i$ is offered at price $p_{il}$. Similar to the assortment optimization problem considered in Section \ref{sec:formulation-assort}, the probability that the customer choose nest $S_n$ is given by $P(S_n|\bx)=W_n^{\sigma_n} / \sum_{n'\in[N]}W_{n'}^{\sigma_{n'}}$, where $W_n = v_{0n} + \sum_{i\in [S_n]}\sum_{l\in[L]}\alpha_{in}x_{il}v_{il}$ is the total preference weights of all alternatives in nest $S_n$. In the case the customer selects $S_n$, the probability that she purchases a product $i\in S_n$ at the offered price $p_{il}$ is determined by $P(i|S_n) = \alpha_{in}x_{il}v_{il}/W_n$. The probability that product $i$ is purchased at price $p_{il}$ can be calculated as $$
	P(i | \bx) =\frac{\sum_{n\in [N]}W^{\sigma_n-1}_n(\alpha_{in}x_{il}v_{il})}{ \sum_{n'\in[N]}W^{\sigma_{n'}}_{n'}}.$$ The objective function  can now be written as
	\[
	F(\bx) = \frac{\sum_{i\in [m]}\sum_{l\in [L]}\sum_{n\in[N]}W^{\sigma_n-1}_n \left(\alpha_{in}x_{il}v_{il}p_{il}\right)}{\sum_{n\in[N]}W^{\sigma_{n}}_{n}}.
	\]
	Under the CNL model, the joint assortment and price optimization problem can now be formulated as follows:
	\begin{align}
		\max_{\bx\in \cX}&\left\{ \ F(\bx) = \frac{\sum_{i\in [m]}\sum_{l\in [L]}\sum_{n\in[N]}W^{\sigma_n-1}_n \left(\alpha_{in}x_{il}v_{il}p_{il}\right)}{\sum_{n\in[N]}W^{\sigma_{n}}_{n}}\right\}\label{prob:assort-price}\tag{\sf A\&P}\\
		\text{subject to}\qquad & \sum_{l\in[L]}x_{il} \leq 1, \forall{i \in [m]}\label{ctr:jap}\\
		& W_n = v_{n0} + \sum_{i\in [S_n]}\sum_{l\in[L]}\alpha_{in}x_{il}v_{il}, \forall{n\in [N]}.\nonumber
	\end{align}
	The constraint \eqref{ctr:jap} is to ensure that each selected product is offered with only one price. The problem \eqref{prob:assort-price} has the same form as \eqref{prob:CNL-assort}, and thus, can be reformulated using the proposed approach and obtain an approximation with a performance guarantee of $\frac{1 - \epsilon}{1 + \epsilon}$ as follows:

\begin{align}
	\max_{\bx \in \cX} \ &  \left\{\widehat{F}(\bx)  = \frac{\sum_{i\in[m], l\in [L],n\in[N]} a_{iln} x_{il} + \sum_{i\in [m], l\in [L],n\in [N],k\in [K_n]}b_{ilnk}s_{ilnk}
	}{c+ \sum_{n\in [N], k\in [K_n]}d_{nk} {z_{nk}} } \right\}\label{prob:assort-price-approx}\tag{\sf A\&P-LFP}\\
	\text{s.t.}\qquad & \sum_{l\in[L]}x_{il} \leq 1, \forall{i \in [m]}\nonumber \\ & v_{0n} + \sum_{i\in [S_n]}\sum_{l\in[L]}\alpha_{in}x_{il}v_{il} = L_n+ \sum_{k\in [K_n]}\Delta_{nk}z_{nk},\;\forall n\in [N] \nonumber \\
	&y_{n,k+1} \leq y_{nk} \leq z_{nk},\;\forall n\in [N], k\in [K_n-1] \nonumber \\
	&z_{n,k+1}\leq y_{nk}, \forall n\in [N], k\in [K_n-2] \nonumber \\
	&s_{ilnk} \leq x_{il}, s_{ilnk} \leq z_{nk},\;\forall i\in [m], l\in[L], n\in [N], k \in [K_n] \nonumber\\
	&s_{ilnk} \geq x_{il}+z_{nk}-1,\;\forall i\in [m], l\in [L], n\in [N], k \in [K_n] \nonumber\\
	& y_{nk} \in \{0,1\},z_{nk}\in[0,1],  s_{ilnk}\in [0, 1], \forall i\in [m], l\in[L], n\in [N], k \in [K_n]\nonumber
\end{align}
where $a_{iln} = \alpha_{in}v_{il}p_{il}f^n(L_n)$, $b_{ilnk}=\alpha_{in}v_{il}p_{il}\gamma^{fn}_k\Delta_{nk}$, $c = \sum_{n\in[N]}g^n(L_n)$, and $d_{nk} = \Delta_{nk}\gamma^{gn}_k$.
Similar to the previous section, the following corollary establishes a performance bound from solving 	\eqref{prob:assort-price-approx}.
	\begin{corollary}
		Let ${\bx^*}$ be the optimal solution of the problem \eqref{prob:assort-price} and $\widehat{\bx}$ be the optimal solution of the approximate problem \eqref{prob:assort-price-approx}. Then, we have $\frac{1-\epsilon}{1 + \epsilon}F(\bx^*) \leq F(\widehat{\bx}) \leq F(\bx^*)$.
	\end{corollary}

    To solve \eqref{prob:assort-price-approx}, we can use the bisection  presented in Subsection \ref{subsec:bisection}, or we can convert \eqref{prob:assort-price-approx} into a MILP using the Charnes-Cooper transformation \citep{charnes1973explicit} and the Glover’s linearization scheme \citep{glover1975improved} as follows:
 \begin{align}
	\max_{\bx \in \cX} \ &  \left\{\widehat{F}(\bx)  = \sum_{i\in[m]}\sum_{l\in [L]}\sum_{n\in[N]} \left(a_{iln} x_{il} + \sum_{k\in [K_n]}b_{ilnk}s_{ilnk}\right)
	 \right\}\label{prob:assort-price-MILP}\tag{\sf A\&P-MILP}\\
	\text{subject to}\qquad & \sum_{l\in[L]}x_{il} \leq 1, \forall{i \in [m]}\nonumber \\
    &cw+ \sum_{n\in [N]}\sum_{k\in [K_n]}d_{nk} {w^z_{nk}} = 1 \nonumber\\ 
    &\sum_{i\in [S_n]}\sum_{l\in[L]}\alpha_{in}v_{il}w^x_{il} = (L_n - v_{0n})w+ \sum_{k\in [K_n]}\Delta_{nk}w^z_{nk},\;\forall n\in [N] \nonumber \\ 
	&w^y_{n,k+1} \leq w^y_{nk} \leq w^z_{nk},\;\forall n\in [N], k\in [K_n-1] \nonumber \\
	&w^z_{n,k+1}\leq w^y_{nk}, \forall n\in [N], k\in [K_n-2] \nonumber \\
	&w^s_{ilnk} \leq x_{il}, w^s_{ilnk} \leq w^z_{nk},\;\forall i\in [m], l\in[L], n\in [N], k \in [K_n] \nonumber\\
	&w^s_{ilnk} \geq w^x_{il}+w^z_{nk}-1,\;\forall i\in [m], l\in [L], n\in [N], k \in [K_n] \nonumber\\
	&w - U^w (1-x_{il}) \leq w^x_{il} \leq w,\ w^x_{il} \leq U^w x_{il}, \forall{i\in[m], l\in[L]} \nonumber \\
		&w - U^w (1-y_{nk}) \leq w^y_{nk} \leq w,\ w^y_{nk} \leq U^w y_{nk}, \forall{n\in[N], k\in[K_n]} \nonumber \\
		&w\in \bbR_+, y_{nk} \in \{0,1\}, w^x_{il}\in \bbR, w^y_{nk} \in \bbR, w^z_{nk} \in \bbR, w^s_{ink}\in \bbR,\nonumber\\
		&\qquad\qquad\qquad\qquad\qquad\qquad\qquad \forall{i\in[m], l\in[L], n\in[N], k\in[K_n]}.\nonumber
\end{align}
	
	\section{Assortment Optimization under Mixture of CNL models}
	\label{sec:mixture}
	In this section, we discuss  the assortment optimization problem under a mixture of CNL models. This is motivated by practical situations where  there can be
multiple customer segments served by the firm and the customers belonging to different segments may
have different preferences and underlying correlation structures. 
 To this end, we consider an assortment optimization problem with $T$ types of customers indexed by $\{1, 2, \ldots, T\}$. We use $\theta_t$ to denote the probability that a customer of type $t\in[T]$ is interested in purchasing from our system, where $\sum_{t\in[T]}\theta_t = 1$. For each product $i \in [m]$, let $v_{ti}$ be its preference weight associated with a customer of type $t$, and $W_{tn}=v_{0n} + \sum_{i\in S_n}\alpha_{in}x_iv_{ti}$ captures the total preference weights of all options in nest $S_n$. Then, the probability that a customer of type $t$ selects a product in nest $S_n$ is given by $P_t(S_n|\bx) = W_{tn}^{\sigma_n}/ \sum_{n'\in[N]}W_{tn'}^{\sigma_{n'}}$. Suppose that the customer has decided to make a purchase from nest $S_n$, the probability that she chooses a product $i\in S_n$ is then given by $P_t(i|S_n) = \alpha_{in}x_iv_{ti}/W_{tn}$. Putting all together, the probability that a product $i\in [m]$ is purchased given an assortment $\bx$ has the following form:
	\begin{align*}
		P(i|x) &= \sum_{t\in[T]}\theta_t\cdot\frac{\sum_{n\in[N]}W_{tn}^{\sigma_n}}{\sum_{n'\in[N]}W_{tn'}^{\sigma_{n'}}}\cdot \frac{\alpha_{in}x_iv_{ti}}{W_{tn}}\\
		&= \sum_{t\in[T]}\theta_t\cdot\frac{\sum_{n\in[N]}W_{tn}^{\sigma_n-1}(\alpha_{in}x_iv_{ti})}{\sum_{n\in[N]}W_{tn}^{\sigma_{n}}}.
	\end{align*}
	The expected revenue across all the customer types can now be calculated as
	\[
	F(\bx) = \sum_{t\in[T]}\theta_t\cdot\frac{\sum_{i\in[m]}\sum_{n\in[N]}W_{tn}^{\sigma_n-1}(\alpha_{in}x_ir_iv_{ti})}{\sum_{n\in[N]}W_{tn}^{\sigma_{n}}}.
	\]
	With this objective function, the assortment optimization problem under the mixture of CNL models can be formulated as follows:
	\begin{align}
		\max_{\bx \in \cX } \qquad &\left\{F(\bx) = \sum_{t\in[T]}\theta_t\cdot\frac{\sum_{i\in[m]}\sum_{n\in[N]}W_{tn}^{\sigma_n-1}(\alpha_{in}x_ir_iv_{ti})}{\sum_{n\in[N]}W_{tn}^{\sigma_{n}}}\right\}\label{prob:MCNL-assort}\tag{\sf MA}\\
        \text{subject to} \quad &W_{tn}=v_{0n} + \sum_{i\in S_n}\alpha_{in}x_iv_{ti}, \forall{t\in[T], n\in [N]}.\nonumber
	\end{align}
	

    The assortment problem \eqref{prob:CNL-assort} can be considered as a special case of problem \eqref{prob:MCNL-assort} with $T = 1$. Since \eqref{prob:CNL-assort} is NP-hard, the problem \eqref{prob:MCNL-assort} is also NP-hard.
    We can solve \eqref{prob:MCNL-assort} by applying our previous approximation approach to every fraction of the objective function. With similar transformations as in Section \ref{sec:approximation}, we can obtain an approximate problem in the form of maximizing a sum of linear fractional functions as follows:
	\begin{align}
		\max_{\bx \in \cX} \ &  \left\{\widehat{F}(\bx)  = \sum_{t\in[T]}\theta_t\cdot\frac{\sum_{i\in[m],n\in[N]} a_{tin} x_i + \sum_{i\in [m],n\in [N],k\in [K_n]}b_{tink}s_{tink}
		}{c_t+ \sum_{n\in [N], k\in [K_n]}d_{tnk} {z_{tnk}}} \right\}\label{prob:MCNL-approx}\tag{\sf MA-LFP}\\
		\text{subject to} \qquad & v_{0n} + \sum_{i\in S_n} \alpha_{in}x_iv_{ti} = L_{tn}+ \sum_{k\in [K_n]}\Delta_{tnk}z_{tnk},\;\forall t\in[T], n\in [N] \nonumber \\ 
		&y_{tn,k+1} \leq y_{tnk} \leq z_{tnk},\;\forall t\in[T], n\in [N], k\in [K_n-1] \nonumber \\
		&z_{tn,k+1}\leq y_{tnk}, \forall t\in[T], n\in [N], k\in [K_n-2] \nonumber \\
		&s_{tink} \leq x_i, s_{tink} \leq z_{tnk},\;\forall t\in[T], i\in [m], n\in [N], k \in [K_n] \nonumber\\
		&s_{tink} \geq x_i+z_{tnk}-1,\;\forall t\in[T], i\in [m], n\in [N], k \in [K_n] \nonumber\\
		& y_{tnk} \in \{0,1\},z_{tnk}\in[0,1],  s_{tink}\in [0, 1], \forall t\in[T], i\in [m], n\in [N], k \in [K_n]\nonumber
	\end{align}
where $a_{tin} = \alpha_{in}v_{ti}r_{i}f^n(L_n)$, $b_{tink}=\alpha_{in}v_{ti}r_{i}\gamma^{fn}_k\Delta_{nk}$, $c = \sum_{n\in[N]}g^n(L_n)$, and $d_{nk} = \Delta_{nk}\gamma^{gn}_k$.
 Since each ratio in the objective function of \eqref{prob:MCNL-assort} is approximated in a similar way as in Section \ref{sec:approximation}, we also has the approximation bounds for each ratio given by $F_t(\bx) \leq \widehat{F_t}(\bx) \leq \frac{1 + \epsilon}{1 - \epsilon}F_t(\bx)$, resulting in an overall approximation bounds for the objective function of \eqref{prob:MCNL-approx} as $F(\bx) \leq \widehat{F}(\bx) \leq \frac{1 + \epsilon}{1 - \epsilon}F(\bx)$. We have the performance guarantee for the problem \eqref{prob:MCNL-approx} given in the following corollary, which can be derived in a similar way as in Theorem \ref{theo:approx_bound}. 
	\begin{corollary}
		Let $\bx^*$ be the optimal solution of the problem \eqref{prob:MCNL-assort} and $\widehat{\bx}$ be the optimal solution to the approximate problem \eqref{prob:MCNL-approx}, then we have $\frac{1 - \epsilon}{1 + \epsilon}F(\bx^*) \leq F(\widehat{\bx}) \leq F(\bx^*)$.
	\end{corollary}

	Due to the nature of the objective function, the bisection method in Subsection \ref{subsec:bisection} can not be applied to solve  \eqref{prob:MCNL-approx}. Nevertheless, since $c_t+ \sum_{n\in [N], k\in [K_n]}d_{tnk} {z_{tnk}} > 0\ \forall{t\in[T]}$, we can use the Charnes-Cooper transformation and the Glover’s linearization scheme to transform \eqref{prob:MCNL-approx} into a MILP as follows
	\begin{align}
		\max_{\bx \in \cX} \ &  \left\{\widehat{F}(\bx)  = \sum_{t\in[T]}\sum_{i\in[m]}\sum_{n\in[N]}\left( a_{tin} w^x_{ti} + \sum_{k\in [K_n]}b_{tink}w^s_{tink}
		\right)\right\}\label{prob:MCNL-MILP}\tag{\sf MA-MILP}\nonumber\\
		\text{subject to}\qquad & c_tw_t+ \sum_{n\in [N]}\sum_{k\in [K_n]}d_{tnk} {w^z_{tnk}} = \theta_t, \forall{t\in[T]}\nonumber\\
		& \sum_{i\in S_n} \alpha_{in}v_{ti}w^x_{ti} = (L_{tn}-v_{0n})w_t+ \sum_{k\in [K_n]}\Delta_{tnk}w^z_{tnk},\;\forall t\in[T], n\in [N] \nonumber \\ 
		&w^y_{tn,k+1} \leq w^y_{tnk} \leq w^z_{tnk},\;\forall t\in[T], n\in [N], k\in [K_n-1] \nonumber \\
		&w^z_{tn,k+1}\leq w^y_{tnk}, \forall t\in[T], n\in [N], k\in [K_n-2] \nonumber \\
		&w^s_{tink} \leq w^x_{ti}, w^s_{tink} \leq w^z_{tnk},\;\forall t\in[T], i\in [m], n\in [N], k \in [K_n] \nonumber\\
		&w^s_{tink} \geq w^x_{ti}+w^z_{tnk}-w_t,\;\forall t\in[T], i\in [m], n\in [N], k \in [K_n] \nonumber\\
		&w_t - U_t^w (1-x_i) \leq w^x_{ti} \leq w_t,\ w^x_{ti} \leq U_t^w x_i, \forall{t\in[t], i\in[m]} \nonumber \\
		&w_t - U_t^w (1-y_{tnk}) \leq w^y_{tnk} \leq w_t,\ w^y_{tnk} \leq U_t^w y_{tnk}, \forall{t\in[t], n\in[N], k\in[K_n]} \nonumber \\
		&w_t\in \bbR_+, y_{tnk} \in \{0,1\}, w^x_{ti}\in \bbR, w^y_{tnk} \in \bbR, w^z_{tnk} \in \bbR, w^s_{tink}\in \bbR,\nonumber\\
		& \qquad\qquad\qquad\qquad\qquad\qquad\forall{t\in[t], i\in[m], n\in[N], k\in[K_n]}.\nonumber
	\end{align}
	
	\section{Numerical Experiments}
	\label{sec:experiment}
    In this section, we evaluate the performance of our proposed methods through extensive numerical experiments.
    The experimental settings and our datasets are described in Subsection \ref{subsec:exp_description}. The computational results for the assortment problem are presented in Subsection \ref{subsec:exp_assort} and  the results for the joint assortment and price optimization problem are reported in Subsection \ref{subsec:exp_jap}. Finally, we provide numerical results for the assortment problem under the mixture of CNL models in Subsection \ref{subsec:exp_mixture}.
	
	\subsection{Experimental Settings}
        \label{subsec:exp_description}
	We conduct experiments on a large number of test instances to evaluate the performance of the bisection method (denoted as BIS) and the MILP approach (denoted as MILP). The tolerance $\delta$ in the halting condition of the bisection method is set to 0.001, which ensures that the optimality gap does not exceed 0.1\%. The methods are tested with three different levels of performance guarantee, i.e., 90\%, 95\%, and 99\%. To achieve this, the parameter $\epsilon$ in the bisection method is set to 0.0521, 0.0251, and 0.0045, while in the MILP method, $\epsilon$ is set to 0.0526, 0.0256, and 0.005, respectively. Note that a higher approximation accuracy requires a larger number of linear approximation segments,  resulting in a larger model in size. We use two criteria to evaluate the performance of our methods, including the required computing time and the solution quality. Since we do not know optimal solutions to our test problems, the solution quality is evaluated based on the best-obtained solution for each instance (denoted by $\textsf{BOS}$), which is ensured to be at least 99\% of the optimum.
	
	Our test instances are randomly generated where the number of products $m$ varies from 20 to 200 and the number of nests $N$ takes a value from $\{5, 10\}$. Since the nests may overlap with others, we define a parameter $\gamma \geq 1$ as the average number of nests to which one product belongs in order to control the overlapping rate. In this study, the value of $\gamma$ is set to 1.2. With this parameter, the cross-nested correlation structure over $m$ products and $N$ nests is constructed as follows. We first randomly generate $m$ products. Next, $\ceil{(\gamma - 1) \times m}$ products are randomly sampled (with replacement) from $m$ products generated in the previous step. All these $\ceil{\gamma \times m}$ elements are then randomly assigned to $N$ nests in such a way that a product can appear in multiple nests, but no more than once in each nest. Finally, the value of allocation parameter $\alpha_{in}$ is generated and normalized such that $\alpha_{in} = 0$ if the product $i$ is not a member of the nest $S_n$ and $\sum_{n\in[N]}\alpha_{in} = 1$, for all product $i \in [m]$.
	In all experiments, we apply two types of cardinality constraints. The first constraint limits the total number of offered products over the whole assortment, i.e., $\sum_{i\in [m]}x_i \leq c$, where $c = \ceil{0.5 \times m}$ is the maximum number of offered products. The second type corresponds to a cardinality constraint that limits the total number of offered products in each nest, i.e., $\sum_{i\in S_n}x_i \leq c_n$, where $c_n = \ceil{0.8 \times |S_n|}$ is the maximum number of offered products in nest $S_n$. Noting that our approximation and solution methods are able to work with any set of linear constraints.
	
	All the experiments are implemented using C++ and run on Intel(R) Xeon(R) CPU E5-2698 v3 @ 2.30GHz. The linear programs are carried out by IBM ILOG CPLEX 22.1 where the number of CPUs is set to 8 cores. 
 
	\subsection{Assortment Optimization under CNL Model}
	\label{subsec:exp_assort}
	We generate two datasets to test the performance of our methods. The first dataset includes small instances where the number of products $m$ varies over \{20, 30\}, while the second dataset includes large instances with $m \in \{50, 100, 150, 200\}$. The dissimilarity parameter $\sigma_n$ of each nest $S_n$ is uniformly generated from $[0.25, \Bar{\sigma}]$, where $\Bar{\sigma}$ takes value from \{0.5, 0.75, 1.00\} for the small dataset, and is fixed to 0.75 for the large instances. We label our small test instances as $(N, m, \Bar{\sigma})\in \{5, 10\}\times\{20, 30\}\times\{0.5, 0.75, 1.00\}$ and the large instances as $(N, m)\in \{5, 10\}\times\{50, 100, 150, 200\}$. For the small dataset, we generate 1000 test instances for each combination $(N, m, \Bar{\sigma})$, which results in 12,000 instances in total. On the other hand, since the model size for the large dataset may prevent us from solving a large number of tests, we only generate 10 instances for each parameter combination $(N, m)$, which results in 80 instances in total. For each instance, to generate the preference $v_i$ and the revenue $r_i$ of product $i$, we first randomly and uniformly generate  $u_i \in (0, 1]$, $X_i \in [0.1, 10]$ and $Y_i \in [0.1, 10]$. These numbers are then used to calculate $r_i$ and $v_i$ as  $r_i = u^2_i \times X_i$ and $v_i = (1 - u_i) \times Y_i$. By doing this, products with higher revenues (i.e. high costs) are more likely to be less preferred. For the alternative of leaving the system without purchasing, we set the preference weight $v_{0n} = 1, \forall n \in [N]$.
	
	\input{tables/assortment_small_new.tex}
	\begin{figure}
		\centering
		\includegraphics[scale=0.8]{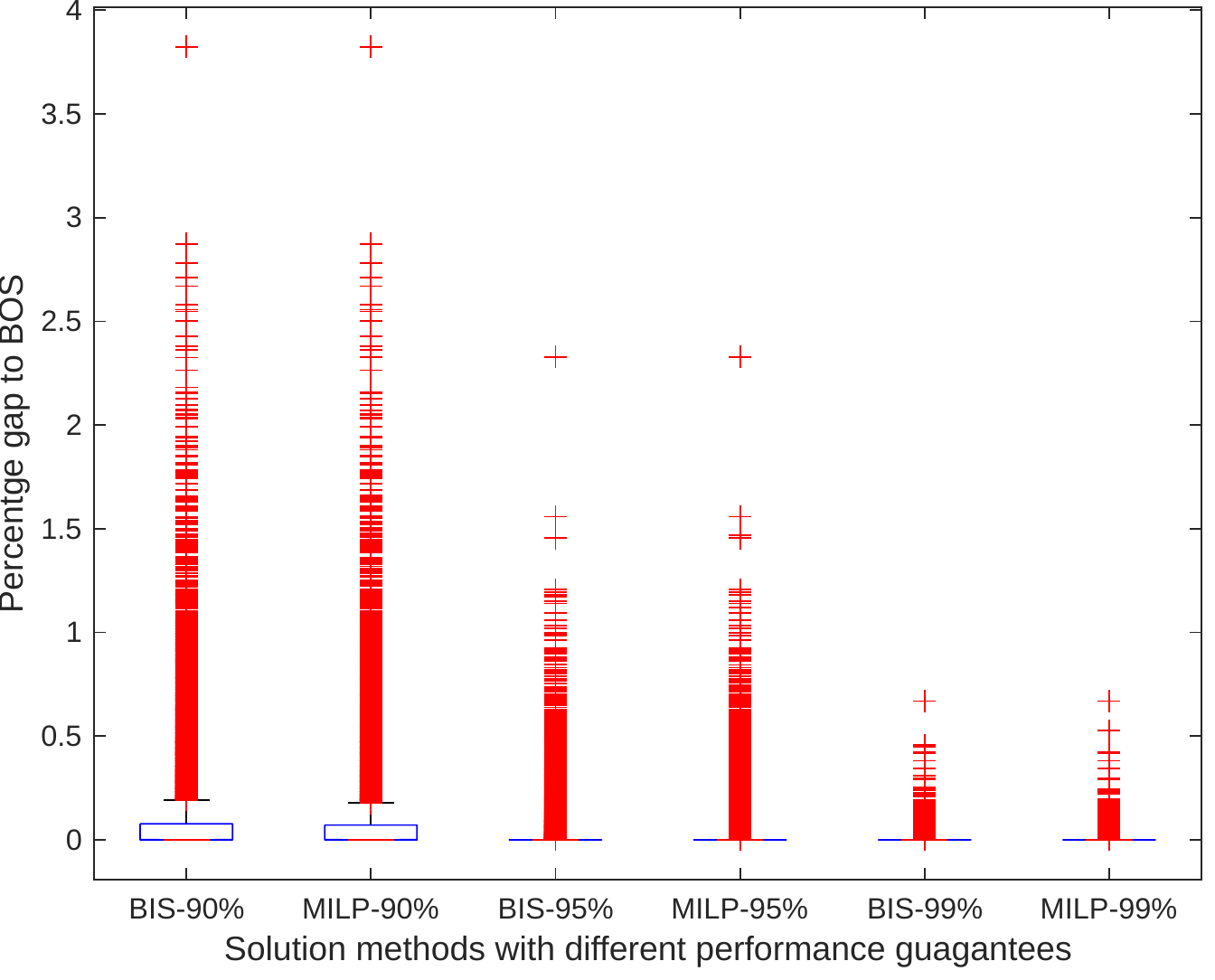}
		\caption{Percentage gaps for small instances of the assortment optimization problem.}
		\label{fig:assort_small_gap}
	\end{figure}
	
	For the small dataset, we run the bisection method (Section \ref{subsec:bisection}) and the MILP \eqref{prob:milp} without setting a time limit. The computational results are reported in Table \ref{table:assort_small}, where better entries are marked in bold. In each row of this table, the first column contains the label ($N, m, \Bar{\sigma}$), and the following columns provide the results under different performance guarantees, including the average runtime (calculated over 1000 instances) and the number of instances (out of 1000) on which the solver yields the \textsf{BOS}. Recall that, for each instance, the \textsf{BOS} is guaranteed to be at least 99\% of the optimal assortment. In terms of solution quality, Table \ref{table:assort_small} shows that our methods perform impressively well. It is notable that the highest approximation guarantee does not always yield the best-obtained solution (\textsf{BOS}), indicating that the methods work very well even with lower guarantees. In particular, the \textsf{BOS} of about two-thirds of the dataset can be found under the guarantee of 90\%, and about three-fourths when the guarantee is increased to 95\%. Regarding the average runtime, the table shows that larger models (caused by a larger number of products, larger 
    number of nests, higher correlation between products, or higher performance guarantee) require more computing time, as expected. However, almost all instances of the dataset can be solved within 1 second under the guarantees of 90\% and 95\%, and 3 seconds under the guarantee of 99\%. In comparison between the two solution methods, we can see that the MILP method outperforms the bisection method in finding the $\textsf{BOS}$. The MILP method is faster than bisection for small instances, as the bisection method requires repeatedly solving MILPs for different values of $\lambda$. In contrast, the MILP approach requires only a single program to be executed to find optimal solutions. However, for large instances, checking the feasibility of $\lambda$ in the bisection method is much faster than solving the MILP, making the bisection method more effective as the model size grows.
	
	In addition, we also evaluate the methods by calculating the percentage gap between the expected revenue generated by these methods and the expected revenue obtained from offering the $\textsf{BOS}$. The gap is calculated as $$100 \times \frac{F(\textsf{BOS}) - F(\bx)}{F(\textsf{BOS})},$$ where $\bx$ is the assortment returned by our method and $\textsf{BOS}$ in this formula is the best solution obtained for the same instance. Figure \ref{fig:assort_small_gap} displays the gaps through box plots, where the outliers, represented by red plus symbols, are 1.5 times larger than the interquartile range. As we can see from the figure, if we ignore outliers, the maximum gaps to $\textsf{BOS}$ are approximately only 0.2\% when the performance guarantee is set to 90\%, and almost negligible when the guarantee is raised to higher levels. From these results, we can indirectly calculate the gap between the obtained solution $\bx$ and the optimal assortment $\bx^*$ as $$\frac{F(\bx^*) - F(\bx)}{F(\bx^*)} = 1 -\frac{F(\bx)}{F(\bx^*)} = 1 - \frac{F(\bx)}{F(\textsf{BOS})}\times \frac{F(\textsf{BOS})}{F(\bx^*)} \leq 1 - (1 - 0.002) \times 0.99 \approx 1.2\%,$$ which indicates that under the performance guarantee of 90\%, the gaps between the objective values yielded by solutions obtained from our methods and the optimal expected revenue are no larger than 1.2\%. 
	
	\input{tables/assortment_large.tex}
	For the large dataset, the solver's time limit is set to 3600 seconds for each instance. Table \ref{table:assort_large} reports the number of instances (out of 10 instances of each $(N, m)$ combination) that are solved to their optimality by our methods, as well as the average computing time required by the solver to confirm  optimality. The results demonstrate that our methods scale very well. Specifically, the methods are able to solve almost all instances with a performance guarantee of 95\% or lower in less than 5 minutes each, except for only the largest instances of 200 products and 10 nests, which require an average running time of 10 minutes. When the guarantee is increased to 99\%, only several of the largest instances are left unsolved within the time budget. These results indicate that the methods are efficient and effective for practical use in real-world situations.
	
	\subsection{Joint Assortment and Price Optimization under CNL Model}
	\label{subsec:exp_jap}
	To test our methods on the joint assortment optimization and pricing problem, we generate a dataset by varying the number of products $m$ over \{20, 30\}, the number of nests $N$ over \{5, 10\}, and the number of price levels $L$ over \{2, 5\}. This creates 8 parameter combinations labeled as $(N, m, L) \in \{5, 10\} \times \{20, 30\} \times \{2, 5\}$. By randomly generating 100 instances for each combination, we obtain a dataset of 800 instances in total. For each instance, the dissimilarity parameters are uniformly generated from [0.25, 0.75].
	In this experiment, we assume that the price $p_{il}$ and the preference weight (or attractiveness) $v_{il}$ follow a parametric relationship suggested by \citet{li2011pricing} and \citet{gallego2014multiproduct}, i.e., $v_{il} = \exp(\mu_{i} - \eta_{i}p_{il})$. This is a decreasing function of the price $p_{il}$, where $\mu_{i}$ is the price-independent utility component (or ideal utility when the price is set to 0) and $\eta_{i} > 0$ is a price-sensitive parameter. For each product $i\in[m]$, we randomly and uniformly generate $\mu_i$ from [-1, 1] and $\eta_i$ from (0, 1). To generate the price set $P_i$, we first generate a random number $\alpha_i$ from (0, 1) and then calculate the prices as $p_{il} =  l\times \alpha_i + 0.5$, $\forall l\in [L]$.
	
\begin{figure}
		\centering
		\includegraphics[scale=0.8]{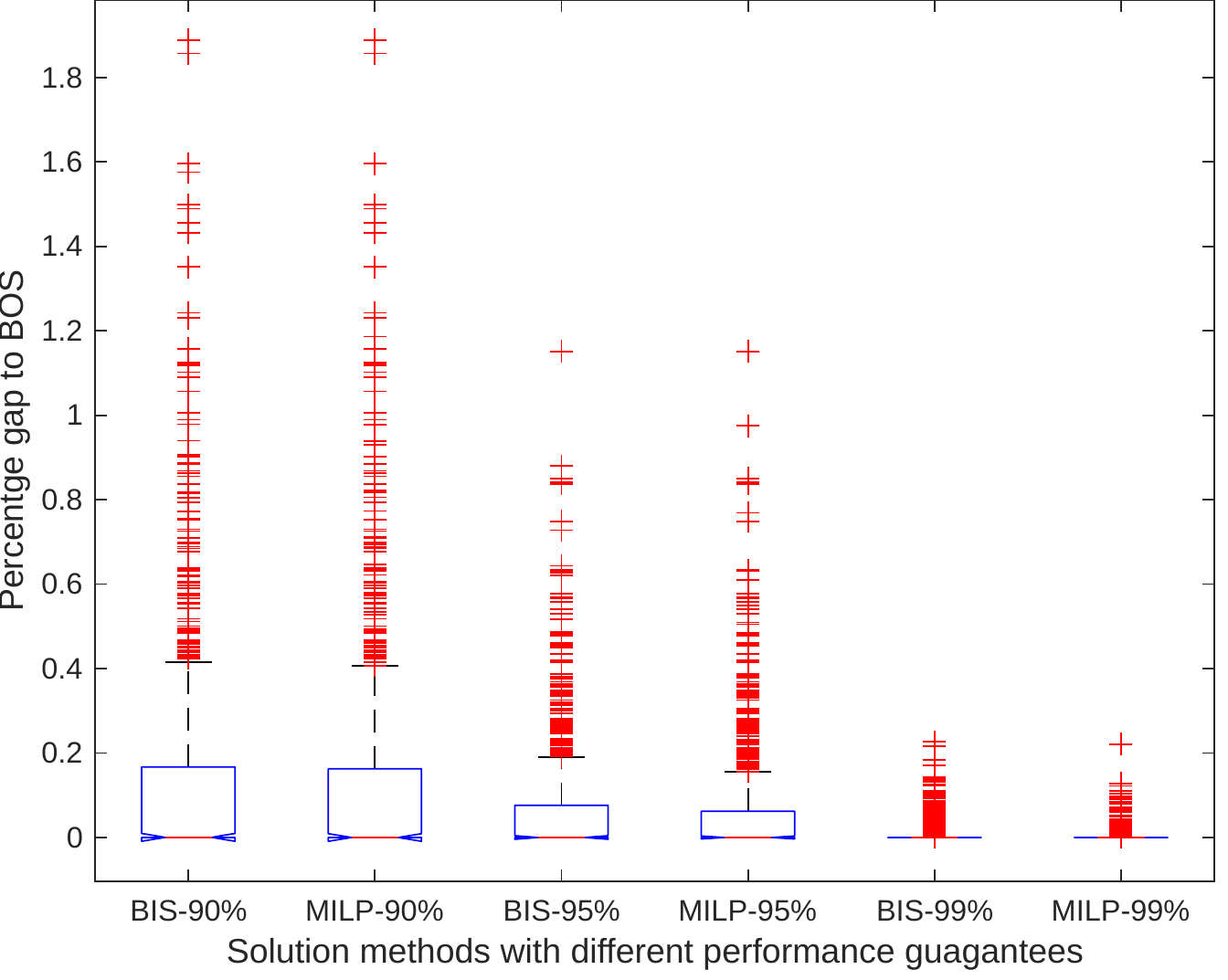}
		\caption{Percentage gap for the joint assortment optimization and pricing problem.}
		\label{fig:jap_small_gap}
	\end{figure}

	We run the bisection method and solve the MILP in \eqref{prob:assort-price-MILP} by CPLEX without setting the time limit. The results are summarized in Table \ref{table:AnP} and Figure \ref{fig:jap_small_gap}, with the same formats as reporting the results of small instances in the previous experiment. It shows that the $\textsf{BOS}$ for a half dataset is attained at the 90\% approximation guarantee. While the MILP method still outperforms bisection in terms of solution quality, bisection is superior in computing time. In particular, while the MILP method takes about 5 minutes to solve the largest instances ($m=30$) under the performance guarantee of 99\%, the bisection method can complete the same task in less than a half minute. Concerning the percentage gap, the boxplots in Figure \ref{fig:jap_small_gap} show that the maximum gaps between the assortments returned by our methods and the corresponding $\textsf{BOS}$ are around 0.4\% at the performance guarantee of 90\%, 0.2\% at the guarantee of 95\%, and very close to zero when the guarantee is set to 99\% (excluding the outliers). These outcomes indicate that our methods are also highly effective in optimizing both assortment and price simultaneously.
	
	\input{tables/assortment_price_small.tex}
	
	\subsection{Assortment Optimization under Mixture of CNL Models}
	\label{subsec:exp_mixture}
	\input{tables/mixture}
	\begin{figure}
		\centering
		\includegraphics[scale=0.8]{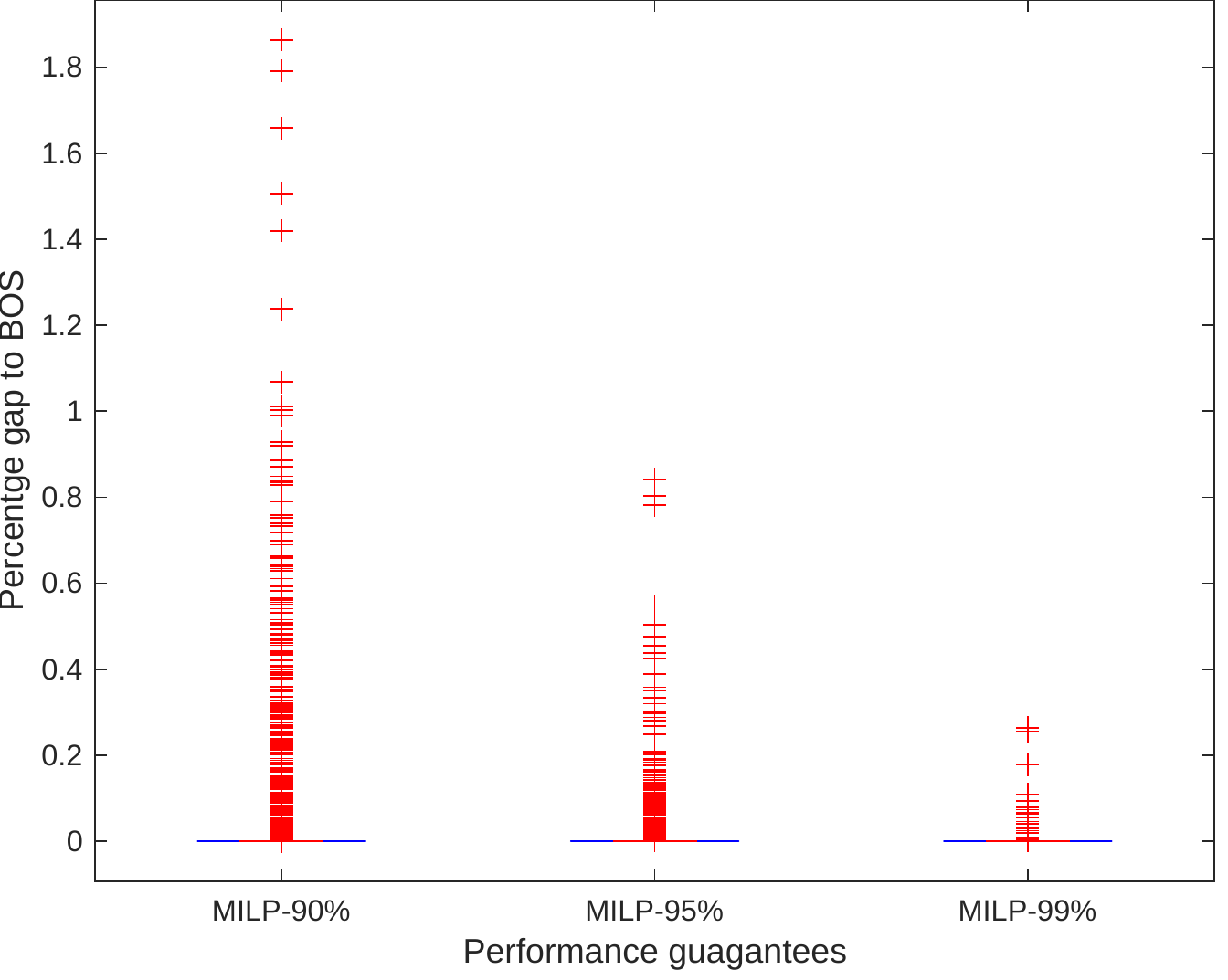}
		\caption{The percentage gap on the assortment problem under mixture of CNL models.}
		\label{fig:mixed_assort_small_gap}
	\end{figure}
	Finally, we test the performance of our approximation method and the MILP in \eqref{prob:MCNL-MILP} when dealing with the assortment problem under a mixture of CNL models. Noting that the bisection method in Section \ref{subsec:bisection} is not applicable for the mixed-CNL model due to the presence of multiple ratios in the objective function. To generate the dataset for this experiment, we vary the number of products $m$ over \{20, 30\}, the number of nests $N$ over \{5, 10\}, and the number of customer types $T$ over \{2, 5, 10\}. By this way, we obtain 12 parameter combinations $(N, m, T) \in \{5, 10\}\times \{20, 30\}\times \{2, 5, 10\}$. For each combination, we generate 100 test instances, resulting in 1200 instances in total. The revenues and preference weights of these test problems are generated using the same way as in Section \ref{subsec:exp_assort}. In particular, for a product $i\in [m]$, we first randomly and uniformly generate  $u_i \in (0, 1]$ and $X_i \in [0.1, 10]$. The revenue of this product is then  calculated as $r_i = u^2_i \times X_i$. For each customer types $t\in [T]$, we uniformly generate  $Y_{ti} \in [0.1, 10]$ and then calculate the preference weights $v_{ti} = (1 - u_i) \times Y_{ti}$, $\forall t\in[T], i\in [m]$. All the parameters $\theta_t$ are generated randomly in [0, 1] and then are normalized so that different customer types have different arriving probabilities and the total probability $\sum_{t\in[T]}\theta_t$ = 1.
	
	The computational results are reported in Table \ref{table:assort_mixture} and Figure \ref{fig:mixed_assort_small_gap}, showing that our method performs remarkably well under the mixture of multiple CNL models. Specifically, at the performance guarantee of 90\%, the outputs of 924/1200 instances are at least as good as 99\% of the optimal assortments. This number increases to 1071/1200 instances when the guarantee is set to 95\%. For the largest instances, it takes 1 minute to confirm the optimality of a solution with the guarantee of 90\%, 2 minutes for 95\%, and 10 minutes for 99\%. Notably, Figure \ref{fig:mixed_assort_small_gap} shows that if we exclude the red outliers, the percentage gap between the objective value yielded by the obtained assortment and the one from the $\textsf{BOS}$ is almost zero for all the levels of accuracy. Including the outliers, the maximum gap w.r.t.  $\textsf{BOS}$ is less than 2\% at the performance guarantee of 90\%, less than 1\% at 95\%, and approximately 0.3\% at 99\%.
	
	\section{Conclusion}
	\label{sec:conclu}
	The question of how businesses can effectively optimize their product assortments to meet the changing needs and preferences of their customers while maximizing revenue is a topic of critical importance in the field of revenue management. This paper, for the first time, studies the assortment optimization problem with general linear constraints under the CNL model. Compared to the traditional MNL and NL models, the CNL model is more general and does not suffer from the limitations of the IIA property or the transitive correlation assumption, and thus, allows more sophisticated substitution patterns between products. However, the resulting optimization problem is also more complex and harder to handle. To overcome the computational challenge, we developed an approximation method with a worst-case performance guarantee of $\frac{1 - \epsilon}{1 + \epsilon}$ for any given $\epsilon>0$. By employing a discretization technique, we have approximated the initial problem by a linear-fractional program. We have demonstrated that optimal solutions to this approximate problem can be derived by solving a single or a series of MILPs. Additionally, we have expanded our analysis to include a joint assortment optimization and pricing problem, as well as an assortment problem involving multiple customer classes. Our experimental results based on various test instances have demonstrated that the proposed methods are highly practical in returning near-optimal solutions. Future research directions would be to consider an assortment problem under a more general class of choice models such as the network Generative Extreme Value model \citep{daly2006general,mai2017dynamic}, or constrained assortment optimization under choice parameter uncertainties.

 \section*{Acknowledgments}
The project is  supported by  the National Research
Foundation Singapore and DSO National Laboratories under the AI Singapore Programme (AISG Award No: AISG2-
RP-2020-017).
	\bibliographystyle{plainnat_custom}
	\bibliography{refs}
	
	\newpage
	\begin{appendices}
 \section{Missing Proofs}
		\subsection{Proof of Theorem \ref{theorem:np_hard}}
        \label{appendix:proof_np_hard}
		To prove the theorem, we show that we can reduce any instance of the set partition problem to an instance of the assortment feasibility problem under the CNL model. Set partition problem is a well-known NP-complete problem \citep{karp2010reducibility} which was also used by \citet{rusmevichientong2014assortment} to establish the NP-hardness of the assortment problem under the mixture of two multinomial logit models.
		
		The set partition problem is defined as follows. Given a set of integers $(c_1, c_2, \ldots, c_m)$ such that $\sum_{i \in [m]}c_i = 2T$. The question is, whether there exists a subset $S$ such that $\sum_{i \in S} c_i = T$. We will construct an instance of the assortment feasibility problem in a such way that there exists a feasible subset $S$ of the partition problem if and only if there exists an assortment $\bx$ such that the expected revenue $F(\bx) \geq K$.
		
		Let's construct the instance for the assortment feasibility problem as follows. There are $m+1$ products that are organized into two nests. For the first $m$ products, we have the revenues $r_i = 6$ and the preference weights $v_i = c_i$ where $c_i$ is the $i^{th}$ element in the partition problem. For the last product, we have $r_{m+1} = 9$ and $v_{m+1}=4T$. The first nest contains all $m+1$ products with the allocation parameters $\alpha_{i1} = 1$ for $i \in [m]$ and $\alpha_{i1} = \frac{1}{4}$ for $i = m+1$. The second nest includes only the last product $m+1$, i.e., $\alpha_{i2} = 0$ for $i \in [m]$ and $\alpha_{i2} = \frac{3}{4}$ for $i = m+1$. We set the preference weights of the no purchase option as $v_{01} = T$ and $v_{02} = 0$. The dissimilarity parameters are set as $\sigma_1 = \sigma_2 = 0.5$. The value of the feasible threshold $K$ is set to 7.
		
		We see that if we offer only the last product $m+1$, then we obtain the expected revenue as
		\[
		F(\bx) = \frac{(\frac{1}{4}\cdot4T + T)^{-0.5}(\frac{1}{4}\cdot4T\cdot9) + (\frac{3}{4}\cdot4T)^{-0.5}(\frac{3}{4}\cdot4T\cdot9)}{(\frac{1}{4}\cdot4T + T)^{0.5} + (\frac{3}{4}\cdot4T)^{0.5}} \ \approx\  6.98
		\]
		which is larger than the revenues of all other products. Therefore, the last product is always included in the optimal assortment since we want to maximize the expected revenue. Thus, the remaining question is whether there exists a subset of first $m$ products such that the obtained expected revenue $F(x) \geq K$. Now suppose that we offer a subset $S$ of the first $m$ products together with the last product $m+1$. Then we have the expected revenue as follows
		\[
		F(\bx) = \frac{(\sum_{i\in S}v_i + \frac{1}{4}\cdot4T + T)^{-0.5}(\sum_{i\in S}v_i\cdot 6+ \frac{1}{4}\cdot4T\cdot9) + (\frac{3}{4}\cdot4T)^{-0.5}(\frac{3}{4}\cdot4T\cdot9)}{(\sum_{i\in S}v_i + \frac{1}{4}\cdot4T + T)^{0.5} + (\frac{3}{4}\cdot4T)^{0.5}}
		\]
		Let denote $C = \sum_{i\in S}v_i$. Since $r_i = c_i$ for $i\in [m]$, we have $C = \sum_{i\in S}c_i$. Then, $F(\bx) \geq K$ is equivalent to
		\[
		\frac{(C + 2T)^{-0.5}(6C+ 9T) + (3T)^{-0.5}(27T)}{(C + 2T)^{0.5} + (3T)^{0.5}} \geq 7
		\]
		or we can rewrite as
		\[
		\frac{6C + 9T}{\sqrt{C + 2T}} + 9\sqrt{3T} \geq 7 \left(\sqrt{C + 2T} + \sqrt{3T}\right).
		\]
		Multiply both sides to $\sqrt{C + 2T}$, the inequality is equivalent to
		\[
		6C + 9T + 9\sqrt{3T}\sqrt{C + 2T} \geq 7 (C + 2T) + 7\sqrt{3T}\sqrt{C + 2T}.
		\]
		We rearrange the terms and rewrite as
		\[
		(C + 2T) - 2\sqrt{3T}\sqrt{C + 2T} + 3T \leq 0,
		\]
		which is equivalent to 
		\[
		\left(\sqrt{3T} -\sqrt{C + 2T}\right)^2 \leq 0.
		\]
		The last inequality holds if and only if $\sqrt{3T} =\sqrt{C + 2T}$, which equivalent to $C = T$, or $\sum_{i\in S}c_i = T$ according to the definition of $C$. This implies that there exists a subset $S$ in the partition problem such that $\sum_{i\in S}c_i = T$ if and only if there exists an assortment $\bx$ of products that satisfies $F(\bx) \geq K$, or finding the subset $S$ is equivalent to finding the assortment $\bx$, which completes the proof.

        \subsection{Proof of Lemma \ref{lemma:monotnone_of_approximators}}
        \label{appendix:proof_monotnone_of_approximators}
        Let's show that
        $\widehat{f^n}(c^n_k, u, t)$ is monotonically  increasing in $u$ by taking its derivative as follows
        \begin{align*}
            \frac{\partial}{\partial u}\widehat{f^n}(c^n_k, u, t) &= 
            \frac{\partial}{\partial u}\left(f^n(c^n_k) + \frac{f^n(u) - f^n(c^n_k)}{u - c^n_k}(t - c^n_k)\right)  = (t - c^n_k) \frac{\partial}{\partial u} \frac{u^{\sigma_n - 1} - (c^n_k)^{\sigma_n - 1}}{u - c^n_k}\\
            &= (t - c^n_k) \frac{(\sigma_n - 2)u^{\sigma_n - 1} - (\sigma_n - 1)c^n_ku^{\sigma_n - 2} + (c^n_k)^{\sigma_n - 1}}{(u - c^n_k)^2}\\
            &= (t - c^n_k) \frac{h(u)}{(u - c^n_k)^2},
        \end{align*}
        where $h(u) = (\sigma_n - 2)u^{\sigma_n - 1} - (\sigma_n - 1)c^n_ku^{\sigma_n - 2} + (c^n_k)^{\sigma_n - 1}$. Let take the derivative of $h(u)$ as
        \begin{align*}
            h'(u) &= (\sigma_n - 2)(\sigma_n - 1)u^{\sigma_n - 2} - (\sigma_n - 1)(\sigma_n - 2)c^n_ku^{\sigma_n - 3} \\
            &= (\sigma_n - 1)(\sigma_n - 2)u^{\sigma_n - 3}(u - c^n_k)\\
            & > 0, 
        \end{align*}
        which indicates that $h(u)$ is increasing in $u$. Moreover, by the L'Hospital's rule, we have
        \begin{align*}
            \lim_{u \downarrow c^n_k}\frac{\partial}{\partial u}\widehat{f^n}(c^n_k, u, t) &=  \lim_{u \downarrow c^n_k} (t - c^n_k)\frac{h(u)}{(u - c^n_k)^2} = (t - c^n_k) \lim_{u \downarrow c^n_k} \frac{h'(u)}{2(u - c^n_k)} \\
            &= (t - c^n_k)\frac{(\sigma_n - 1)(\sigma_n - 2)u^{\sigma_n - 3}}{2}\\
            & > 0,
        \end{align*}
        implying that $\lim_{u \downarrow c^n_k}h(u) > 0$. 
        Combine all together, we have $h(u) > 0\ \forall u \in (c^n_k, U_n]$. It follows that $\frac{\partial}{\partial u}\widehat{f^n}(c^n_k, u, t) > 0\ \forall u \in (c^n_k, U_n]$, or $\widehat{f^n}(c^n_k, u, t)$ is monotonically increasing in $u$, as desired.

        Similarly, we can prove $\widehat{g^n}(c^n_k, u, t)$ is monotonically  decreasing in $u$, which completes the proof.

\subsection{Proof of Proposition \ref{lemma:monotone_of_max_error}}

	Using contradiction, assume that there exists $u' > u$ such that $\Phi^{fn}(c^n_k, u) > \Phi^{fn}(c^n_k, u')$, or $\max_{c^n_k \leq t \leq u}\{\phi^{fn}(c^n_k, u, t)\}\  > \ \max_{c^n_k \leq t \leq u'}\{\phi^{fn}(c^n_k, u', t)\}$, according the the definition of $\Phi^{fn}(c^n_k, u)$. This is equivalent to
	\begin{equation}
		\max_{c^n_k \leq t \leq u}\left\{\frac{\widehat{f^n}(c^n_k, u, t) - f^n(t)}{f^n(t)}\right\} \quad > \quad \max_{c^n_k \leq t \leq u'}\left\{\frac{\widehat{f^n}(c^n_k, u', t) - f^n(t)}{f^n(t)}\right\}.
		\label{eq:Phi_monotonic}
	\end{equation}
	\noindent It follows from Lemma \ref{lemma:monotnone_of_approximators} that $\widehat{f^n}(c^n_k, u', t) > \widehat{f^n}(c^n_k, u, t)$ for every $u' > u$. Thus, \eqref{eq:Phi_monotonic} is equivalent to
	\begin{equation*}
		\max_{c^n_k \leq t \leq u}\left\{\frac{\widehat{f^n}(c^n_k, u, t) - f^n(t)}{f^n(t)}\right\} \quad > \quad \max_{c^n_k \leq t \leq u'}\left\{\frac{\widehat{f^n}(c^n_k, u, t) - f^n(t)}{f^n(t)}\right\},
	\end{equation*}
	which contradicts the fact that $[c^n_k, u] \subset [c^n_k, u']$. It follows that for any $u' > u$, we have $\Phi^{fn}(c^n_k, u') \geq \Phi^{fn}(c^n_k, u)$, indicating that the function $\Phi^{fn}(c^n_k, u)$ is monotonically non-decreasing in $u$. We can also prove that $\Phi^{gn}(c^n_k, u)$ is monotonically non-decreasing in $u$ using the same contradiction and complete the proof.

\subsection{Proof of Lemma \ref{lm:bound-phi}}
 We first recall that
\[
\Phi^{fn}(c, u)  = \max_{t|~c\leq t\leq u} \left\{ \phi^{fn} (c,u,t)\right\} = \max_{t|~c\leq t\leq u} \left\{ \frac{\widehat{f^n}(c,u,t)}{f^n(t)}-1\right\}.
\]
Since both $\widehat{f^n}(c,u,t)$ and $f^n(t)$ are monotonically decreasing on $[c, t]$, we  have $\frac{\widehat{f^n}(c,u,t)}{f^n(t)} \leq \frac{f^n(c)}{f(u)} = (c/u)^{\sigma_n-1}$, implying that $\Phi^{fn}(c, u) \leq (c/u)^{\sigma_n-1} -1$. Similarly, we  also see that both  $\widehat{g^n}(c,u,t)$ and $g^n(t)$ are increasing on $[c,u]$, thus $\widehat{g^n}(c,u,t)/g^n(t)\geq (c/u)^{\sigma_n}$. It then follows that 
\[
\Phi^{gn}(c, u)  = \max_{t|~c\leq t\leq u} \left\{ \phi^{gn} (c,u,t)\right\} = \max_{t|~c\leq t\leq u} \left\{ 1 - \frac{\widehat{g^n}(c_n^k,u,t)}{f^n(t)}\right\} \leq 1-\left(\frac{c}{u}\right)^{\sigma_n}.
\]
Moreover, to have lower bounds for $
\Phi^{fn}(c, u)$  and $\Phi^{gn}(c, u)$, we just evaluate  $\phi^{fn} (c,u,t)$ and $\phi^{gn} (c,u,t)$ at the center point of the interval, i.e., $t = (c+u)/2$, as follows:
\begin{align}
    \Phi^{fn}(c, u)  &= \max_{t|~c\leq t\leq u} \left\{ \phi^{fn} (c,u,t)\right\}  \geq
\phi^{fn} \left(c,u,\frac{c+u}{2}\right) \nonumber \\
&= \frac{f^n(c)+f^n(u)}{2 f^n(\frac{c+u}{2})}-1\nonumber \\
&= \frac{2^{\sigma_n-2} ((c)^{\sigma_n-1} +u^{\sigma_n-1} ) }{\left(c+u\right)^{\sigma_n-1}}
-1,\nonumber
\end{align}
and
\begin{align}
    \Phi^{gn}(c, u)  &= \max_{t|~c\leq t\leq u} \left\{ \phi^{gn} (c,u,t)\right\}  \geq
\phi^{gn} \left(c,u,\frac{c+u}{2}\right) \nonumber \\
&= 1 - \frac{g^n(c)+g^n(u)}{2 g^n(\frac{c+u}{2})}\nonumber \\
&= 1 - \frac{2^{\sigma_n-1} ((c_k^n)^{\sigma_n} +u^{\sigma_n} ) }{\left(c+u\right)^{\sigma_n}},\nonumber
\end{align}
which gives us the desired bounds.

\subsection{Proof of Lemma \ref{lm:lm5}}

Since $\Phi^{fn}(c^{fn}_k, c^{fn}_{k+1}) = \epsilon$ for all $k \in[K^{f}_n-1]$ and $\Phi^{gn}(c^{gn}_k, c^{gn}_{k+1}) = \epsilon$ for all $k \in[K^{g}_n-1]$, using the inequalities in Lemma \ref{lm:bound-phi} we have 
\begin{align}
    \left(\frac{c^{fn}_k}{c^{fn}_{k+1}}\right)^{\sigma_n-1} \geq \epsilon +1, ~ k \in [K^{f}_n-1] \nonumber \\
     \left(\frac{c^{fn}_k}{c^{gn}_{k+1}}\right)^{\sigma_n} \leq 1 - \epsilon , ~ k \in [K^{g}_n-1], \nonumber
\end{align}
which implies
\begin{align}
   c^{fn}_{k+1}\geq c^{fn}_k \left(\epsilon+1\right)^{1/(1-\sigma_n)}, ~ k \in [K^{f}_n-1] \label{eq:proof-lm5-eq1} \\
    {c^{gn}_{k+1}}  \geq {c^{gn}_{k}}(1 - \epsilon)^{-1/\sigma_n} , ~ k \in [K^{g}_n-1]. \label{eq:proof-lm5-eq2}
\end{align}
For ease of notation, let  $\rho^f = \left(\epsilon+1\right)^{1/(1-\sigma_n)}$ and $\rho^g = (1 - \epsilon)^{-1/\sigma_n}$. 
From \eqref{eq:proof-lm5-eq1} and \eqref{eq:proof-lm5-eq2} we will have  
\begin{align}
 c^{fn}_{K^f_n} \geq (\rho^f)^{K^f_n-1}  c^{fn}_{1} =  L_n (\rho^f)^{K^f_n-1} \nonumber \\  
  c^{gn}_{K^f_n} \geq (\rho^g)^{K^g_n-1}  c^{gn}_{1} =  L_n (\rho^g)^{K^g_n-1}. \nonumber
\end{align}
In addition, we have $U_n \geq c^{fn}_{K^f_n}$ and $U_n \geq c^{gn}_{K^f_n}$. Putting all the above together, we have
\begin{align}
 (\rho^f)^{K^f_n-1} \leq \frac{U_n}{L_n}~\text{ or } K^f_n \leq \frac{\ln (U_n/L_n)}{\ln \rho^f} +1  \nonumber \\
 (\rho^g)^{K^g_n-1} \leq \frac{U_n}{L_n}~\text{ or } K^g_n \leq \frac{\ln (U_n/L_n)}{\ln \rho^g} +1. \nonumber
\end{align}
We obtain the desired upper bounds for $K^f_n$ and $K^g_n$. 

To validate the lower bounds, we will identify two scalars $\alpha, \beta>1$ such that  $c^{fn}_{k+1}\leq c^{fn}_k \alpha$ for all $k\in [K^f_n]$  and $c^{gn}_{k+1}\leq c^{gn}_k \beta$ for all $k\in [K^g_n]$. To this end,  we will utilize the lower bound functions established in Lemma \ref{lm:proof-lm5-lm1}. To facilitate the later exposition, let us consider the following equations, for any given $\epsilon>0$:
 \begin{align}
        \frac{t^{\sigma_n-1}+1}{(t+1)^{\sigma_n-1}} = 2^{2-\sigma_n}(\epsilon+1) \label{eq:lmlm5-eq1}\\
        \frac{t^{\sigma_n}+1}{(t+1)^{\sigma_n}} = 2^{1-\sigma_n}(1-\epsilon) \label{eq:lmlm5-eq2}
    \end{align}
And let  $t^f(\epsilon)$  and $t^g(\epsilon)$ be solutions to \eqref{eq:lmlm5-eq1} and \eqref{eq:lmlm5-eq2}, respectively.
\begin{lemma}\label{lm:proof-lm5-lm1}
Given any $\epsilon \in (0,1-2^{\sigma_n-1}]$, the following results hold
\begin{itemize}
    \item[(i)] $t^f(\epsilon)$  and $t^g(\epsilon)$ are always uniquely determined  in $[1, \infty]$.
   
    \item[(ii)]  The following holds
      \begin{itemize}
          \item[(ii-1)] $t\leq t^f(\epsilon)$ if and only if 
          $
               \frac{t^{\sigma_n-1}+1}{(t+1)^{\sigma_n-1}} \leq 2^{2-\sigma_n}(\epsilon+1) 
           $
           \item[(ii-2)] $t\leq t^
           g(\epsilon)$ if and only if 
          $
               \frac{t^{\sigma_n}+1}{(t+1)^{\sigma_n}} \geq 2^{1-\sigma_n}(1-\epsilon) 
           $
      \end{itemize}
    \item[(iii)] $t^f(\epsilon)> 1$, $t^f(\epsilon)> 1$, and they are monotonically decreasing as $\epsilon$ decreases,  and $\lim_{\epsilon\rightarrow 0} \max\{t^f(\epsilon),t^g(\epsilon)\} = 1$.
\end{itemize}  
\end{lemma}
\proof{Proof.}
For \eqref{eq:lmlm5-eq1}, we see that it is equivalent to 
\[
\begin{aligned}
&t^{\sigma_n-1}+1 = 2^{2-\sigma_n} (\epsilon+1) (t+1)^{\sigma_n-1}\\
\text{ or } &\ln (t^{\sigma_n-1}+1) + (1-\sigma_n)\ln (t+1) = \ln (2^{2-\sigma_n} (\epsilon+1)).
\end{aligned}
\] 
Let us now consider the function $\delta(t) = \ln (t^{\sigma_n-1}+1) + (1-\sigma_n)\ln (t+1)$. Taking the first-order derivative of $\delta(t)$ w.r.t $t$, for any $t\geq 1$,  we have 
\begin{align}
\frac{\partial \delta(t)}{\partial t} &= (\sigma_n-1)t^{\sigma_n-2}\frac{1}{t^{\sigma_n-1}+1} + (1-\sigma_n)\frac{1}{t+1}    \nonumber \\
&= \frac{1-\sigma_n}{(t^{\sigma_n-1}+1)(t+1)}(t^{\sigma_n}+1 - t^{\sigma_n-1} - t^{\sigma_n-2}) \nonumber\\
&= \frac{1-\sigma_n}{(t^{\sigma_n-1}+1)(t+1)}( 1 - t^{\sigma_n-2})\stackrel{(a)}{\geq}0, 
\end{align}
where $(a)$ is because $t\geq 1$, thus $1 - t^{\sigma_n-2}\geq 0$. So, $\delta(t)$ is  monotonically increasing on $[1,0]$. On the other hand, $\delta(1) = (2-\sigma_n)\ln 2 < \ln (2^{2-\sigma_n} (\epsilon+1))$ and $\lim_{t\rightarrow +\infty}\delta(t)  = +\infty$. All these imply that the equation $\delta(t) = \ln (2^{2-\sigma_n} (\epsilon+1))$ always has a unique solution at $t = t^f(\epsilon)>1$. 

In a similar way, \eqref{eq:lmlm5-eq2} is equivalent to 
\[
\ln(t^{\sigma_n}+1) - \sigma_n \ln (t+1) = (1-\sigma_n)\ln 2 + \ln (1-\epsilon).
\]
We then consider function $\rho(t) =\ln(t^{\sigma_n}+1) - \sigma_n \ln (t+1) $. Taking the first-order derivative of $\rho(t)$ w.r.t. $t\geq 1$ we have 
\begin{align}
   \frac{\partial \rho(t)}{\partial t} &= \frac{\sigma_n t^{\sigma_n-1}}{t^{\sigma_n}+1} - \frac{\sigma_n}{t+1} \nonumber \\ 
   &= \sigma_n \frac{t^{\sigma_n }+t^{\sigma_n-1} - t^{\sigma_n} - 1}{(t^{\sigma_n}+1)(t+1)} = \sigma_n \frac{t^{\sigma_n-1} -1}{(t^{\sigma_n}+1)(t+1)}\stackrel{(b)}{\leq} 0, 
\end{align}
where $(b)$ is due to the fact that $t\geq 1$, thus $t^{\sigma_n-1}\leq 1$. So, $\rho(t)$ is monotonically decreasing in  $t$ on $[1,\infty)$. Moreover,  $\rho(1) =  (1-\sigma_n)\ln 2 > (1-\sigma_n) \ln 2 + \ln (1-\epsilon)$ for any $\epsilon\geq 0$, and $\lim_{t\rightarrow \infty}\rho(t) = 0$. In addition,  since $\epsilon \in (0, 1-2^{\sigma_n-1}]$, we have $2^{1-\sigma_n}(1-\epsilon)\geq 1$ and $(1-\sigma_n)\ln 2 + \ln (1-\epsilon) \geq 0$, implying that $\rho(1) > (1-\sigma_n)\ln 2 + \ln (1-\epsilon) \geq \lim_{t\rightarrow \infty }\rho(t)$.  As a result, the equation $\rho(t) = (1-\sigma_n) \ln 2 + \ln (1-\epsilon)$ always  has a unique solution at $ t^g(\epsilon)>1$.

The two claims \textit{(ii-1)} and \textit{(ii-2)} can be directly verified based on the results that $\delta(t)$ is monotonically increasing and $\rho(t)$ is monotonically decreasing over the interval $[1,\infty]$. 
Recall that $t^f(\epsilon)$ and $t^g(\epsilon)$ are unique solutions to the equations $\delta(t) = \ln (2^{2-\sigma_n} (\epsilon+1))$ and $\rho(t) = (1-\sigma_n) \ln 2 + \ln (1-\epsilon)$, respectively. So, 
\textit{(iii)} can be verified based on the observation that $\lim_{\epsilon\rightarrow 0} |\delta(1)- \ln (2^{2-\sigma_n} (\epsilon+1))|=0$ and $\lim_{\epsilon\rightarrow 0} |\rho(1)- (1-\sigma_n) \ln 2 + \ln (1-\epsilon)|=0$.  We complete the proof of Lemma \ref{lm:proof-lm5-lm1}.
\endproof

Let us now return to the proof of Lemma \ref{lm:lm5}. From Lemma \ref{lm:bound-phi}, we see that, for any $k\in [K^f_n]$
\begin{align}
     &\Phi^{fn}(c^{fn}_k, c^{fn}_{k+1}) =  \epsilon \geq   \frac{2^{\sigma_n-2}((c^{fn}_{k+1}))^{\sigma_n-1} + (c^n_k)^{\sigma_n-1})}{(c^{fn}_{k+1}+c^n_k)^{\sigma_n-1}}-1.\nonumber \\
     \text{Thus, }~& (1+\epsilon)2^{2-\sigma_n}\geq \frac{t^{\sigma_n-1} + 1}{(t+1)^{\sigma_n-1}}, 
\end{align}
where $t = c^{fn}_{k+1}/c^{fn}_{k}$ for notational simplicity. Lemma \ref{lm:proof-lm5-lm1} then tells us that $c^{fn}_{k+1}/c^{fn}_{k} \leq t^f(\epsilon)$. In a similar way, we also have $c^{gn}_{k+1}/c^{gn}_{k} \leq  t^g(\epsilon)$ for any $k\in [K^g_n]$. We now note that if we continue to find a next point $c^{fn*}_{K^f_n+1}$ after $c^{fn}_{K^f_n}$ such that $\Phi^{fn}(c^{fn}_{K^f_n}, c^{fn*}_{K^f_n+1}) = \epsilon$, then  $c^{fn*}_{K^f_n+1}) \geq U_n$. From the inequality $c^{fn}_{k+1}/c^{fn}_{k} \leq t^f(\epsilon)$, $\forall k\in [K^f_n]$, we see that $U_n\leq c^{fn*}_{K^f_n+1}) \leq t^f(\epsilon) c^{fn}_{K^f_n} \leq t^f(\epsilon)^{K^f_n} c^{fn}_1$, implying 
\[
(t^f(\epsilon))^{K^f_n} \geq \frac{U_n}{L_n}~\text{ or } K^f_n \geq \frac{\ln (U_n/L_n)}{\ln t^f(\epsilon)}.
\]
In a similar way, we also have $$K^g_n\geq\frac{\ln (U_n/L_n)}{\ln t^g(\epsilon)},$$ which completes the proof.

    \subsection{Proof of Theorem \ref{theorem:partitioning_convergence}}
    \label{appendix:proof_convergence}
    In this proof, given any $n\in [N]$, we first show that, for a given accuracy level $\epsilon > 0$, the number of sub-intervals required by our approximation to achieve the performance guarantee of $\frac{1 - \epsilon}{1 + \epsilon}$ is at most $K^{f}_n + K^{g}_n - 1$ and is at least $\max \{K^f_n;~ K^f_n\}$. The lower and upper bounds then are just direct results from Lemma \ref{lm:lm5}.

We first let $\{c^{fn}_1,\ldots,c^{fn}_{K^f_n+1} \}$  and  $\{c^{gn}_1,\ldots,c^{gn}_{K^g_n+1}\}$ be the set of breakpoints generated for $f^n(\cdot)$ and $g^n(\cdot)$, respectively, to optimally achieve the approximation guarantee. That is,  $\Phi^{fn}(c^{fn}_{k}, c^{fn}_{k+1}) = \epsilon,\ \forall k < K^{f}_n$ and $\Phi^{fn}(c^{fn}_{k}, c^{fn}_{k+1}) \leq \epsilon$ for $k = K^{f}_n$, and $\Phi^{gn}(c^{gn}_{k}, c^{gn}_{k+1}) = \epsilon,\ \forall k < K^{g}_n$ and $\Phi^{gn}(c^{gn}_{k}, c^{gn}_{k+1}) \leq \epsilon$ for $k = K^{g}_n$. We also let $\{c^n_1,\ldots,c^n_{K_n+1}\}$ be the set of the breakpoints returned by Algorithm \ref{algo:DA}. According to Theorem \ref{theorem:bisection_guarantee}, each of $K_n, K^f_n, K^g_n$ is optimal for the respective approximation.
    
    Without the loss of generality, we assume that $K^f_n\geq K^g_n$, noting that the other case can be handled in the same way. By contradiction, assume that $K_n < \max\{K^f_n;~ K^g_n\} = K^f_n$. Since the points $\{c^n_1,\ldots,c^n_{K_n+1}\}$ are chosen in such a way that $\phi^{fn}(c^n_{k}, c^n_{k+1})\leq \epsilon$ and $\phi^{gn}(c^n_{k}, c^n_{k+1})\leq \epsilon$, the set $\{c^n_1,\ldots,c^n_{K_n+1}\}$ offers the same (or even better) performance guarantee as that given by $\{c^{fn}_1,\ldots,c^{fn}_{K^f_n+1} \}$. However, our contradiction assumption says that $K_n < K^f_n$, which contradicts to the fact that the set $\{c^{fn}_1,\ldots,c^{fn}_{K^f_n+1} \}$ is optimal for the approximation of $f^n(\cdot)$. We then can conclude that $K_n\geq K^f_n = \max\{K^f_n;~ K^g_n\}$.  
    
   For the upper bounds, we unite the two sets of breakpoints $\{c^{fn}_1,\ldots,c^{fn}_{K^f_n+1} \}$  and  $\{c^{gn}_1,\ldots,c^{gn}_{K^g_n+1}\}$ to form a new set of $K^f_n+ K^g_n$ points $ \{c^{n'}_1,\ldots,c^{n'}_{K^f_n+ K^g_n}\}$, where $c^{n'}_k \leq c^{n'}_{k+1}$ for all $k\in [K^f_n+ K^g_n]$. Here we note that there are $K^f_n+ K^g_n+2$ points in total, but the two first points and the two last points of the two series are identical, i.e., $c^{fn}_1 = c^{gn}_1 =L_n$ and $c^{fn}_{K^f_n+1} = c^{gn}_{K^g_n+1} =U_n$. It can be seen that, for any $k\in [K^f_n+ K^g_n]$, there are always two sub-intervals -- one from the breakpoints $\{c^{fn}_1,\ldots,c^{fn}_{K^f_n+1}\}$ and another from $\{c^{fn}_1,\ldots,c^{fn}_{K^f_n+1}\}$ -- that contains $[c^{n'}_k,c^{n'}_{k+1}]$. In other words, there exist $k_1 \in [K^f_n]$ and $k_2 \in [K^g_n]$ such that $[c^{n'}_k,c^{n'}_{k+1}] \subseteq [c^{fn}_{k_1}, c^{fn}_{k_1+1}]$ and $[c^{n'}_k,c^{n'}_{k+1}] \subseteq [c^{gn}_{k_2}, c^{gn}_{k_2+1}]$. This implies, for any $k\in [K^f_n+ K^g_n]$, 
   \begin{align}
       \phi^{fn}(c^{n'}_k,c^{n'}_{k+1}) &\leq \phi^{fn}(c^{fn}_{k_1}, c^{fn}_{k_1+1}) \leq \epsilon\nonumber \\
       \phi^{gn}(c^{n'}_k,c^{n'}_{k+1}) &\leq \phi^{gn}(c^{gn}_{k_2}, c^{gn}_{k_2+1}) \leq \epsilon.\nonumber 
   \end{align}
 As a result, the new set of breakpoints  $ \{c^{n'}_1,\ldots,c^{n'}_{K^f_n+ K^g_n}\}$ can offer a similar or better approximation guarantee as the set $\{c^n_1,\ldots,c^n_{K_n+1}\}$. The optimality property shown in Theorem \ref{theorem:bisection_guarantee} then implies that  the number of sub-intervals from $\{c^{n'}_1,\ldots,c^{n'}_{K^f_n+ K^g_n}\}$ cannot be smaller than $K_n$, i.e., 
 $K_n \leq K^f_n+ K^g_n-1$. We complete the proof.

\subsection{Proof of Corollary \ref{coro:bounds-Kn}}

For any $k\in [K^f_n]$, let $t = c^{fn}_{k+1}/c^{fn}_k$,  we have
\[
(1+\epsilon)2^{2-\sigma_n} \geq \frac{t^{\sigma_n-1} +1}{(t+1)^{\sigma_n-1}} \geq \frac{1}{(t+1)^{\sigma_n-1}},
\]
which implies that 
\[
\begin{aligned}
&(1+\epsilon)2^{2-\sigma_n} \geq (t+1)^{1-\sigma_n}\\
\text{ or } &t \leq \exp\left(\frac{\ln(1+\epsilon)+(2-\sigma_n)\ln 2}{1-\sigma_n}\right) - 1 \stackrel{\mathrm{def}}{=} t^{f*}.  
\end{aligned}
\]
Thus, $c^{fn}_{k+1}/c^{fn}_{k} \leq t^{f*}$. Using the same arguments as in the proof of Lemma \ref{lm:lm5}, we can show that $$K^f_n \geq \frac{\ln(U_n/L_n)}{\ln t^{f*}}.$$
It can also be seen that 
\[
\frac{\ln(1+\epsilon)+(2-\sigma_n)\ln 2}{1-\sigma_n} > \frac{(2-\sigma_n)\ln 2}{1-\sigma_n} \geq \ln 2,
\]
implying $t^{f*}> 1$. 

We now consider $K^g_n$. In a similar way, we denote  $t = c^{gn}_{k+1}/c^{gn}_k$ and write 
\[
 2^{1-\sigma_n}(1-\epsilon) 
 \leq \frac{t^{\sigma_n}+1}{(t+1)^{\sigma_n}}\leq \frac{t^{\sigma_n}+1}{t^{\sigma_n}}, 
\]
which implies
\[
t\leq \left(\frac{1}{2^{1-\sigma_n}(1-\epsilon)-1}\right)^{1/\sigma_n}\stackrel{\mathrm{def}}{=} t^{g*},
\]
noting that $2^{1-\sigma_n}(1-\epsilon)-1>0$ due to our assumption that $\epsilon \in (0,1-2^{\sigma_n-1}]$. It then follows that $c^{gn}_{k+1}/c^{gn}_{k} \leq t^{g*}$ for all $k\in [K^g_n]$. Similarly, we have
\[
K^g_n \geq \frac{\ln(U_n/L_n)}{\ln t^{g*}}.
\]
We complete the proof.

	 \section{General Dissimilarity Parameters}
    In the main paper, we investigate the assortment optimization problem under the CNL model, assuming that the dissimilarity parameters of all nests varied within the unit interval. We now focus on  cases where there are no restrictions on these parameters, i.e., $\sigma_n$ can be greater than 1 for any nest $S_n$. 
    This setting would be of relevance in some choice modeling contexts \citep{train2009discrete}. 
    While this relaxation has no effect on our previous formulations, it does result in significant changes in the properties of the approximation problem \eqref{prob:approx} due to the changes in the non-linear components $f^n(W_n)$ and $g^n(W_n)$ of the objective function \eqref{eq:obj_assort}. In this appendix, we show that the performance guarantee for the approximation problem \eqref{prob:approx1} still holds for $\sigma_n > 1$, and that the approximation method presented in Subsection \ref{subsec:discretization} can still be used for this case.

    \subsection{General Performance Guarantee}
    For a given accuracy level $\epsilon > 0$ and dissimilarity parameter $\sigma_n > 0$, $n\in [N]$, we first establish the approximation bounds for the approximate problem \eqref{prob:approx} in the following lemma.
    \begin{lemma}
        \label{lemma:general_approx_bounds}
        For any given $\epsilon > 0$, the assortment problem \eqref{prob:CNL-assort} with dissimilarity parameters $\sigma_n > 0,\ n\in[N]$ can be approximated by problem \eqref{prob:approx} where the approximate objective value is bounded by $\frac{1 - \epsilon}{1 + \epsilon}F(\bx) \leq \widehat{F}(\bx) \leq \frac{1 + \epsilon}{1 - \epsilon}F(\bx)$.
    \end{lemma}
    \proof{Proof.}
        For a given $\epsilon > 0$, the piecewise linear approximators $\widehat{f^n}(W_n)$ and $\widehat{g^n}(W_n)$ need to be constructed in such a way that, for every sub-interval $[c^n_k, c^n_{k+1}]$, the following conditions hold
        \begin{equation}
            \max_{c^n_k \leq W_n \leq c^n_{k+1}} \left\{\frac{\Big|f^n(W_n) - \widehat{f^n}(W_n)\Big|}{f^n(W_n)}\right\} \leq \epsilon
            \label{eq:appendix_phi_f_12}
        \end{equation}
        and 
        \begin{equation}
           \max_{c^n_k \leq W_n \leq c^n_{k+1}}\left\{\frac{\Big|  g^n(W_n)-\widehat{g^n}(W_n)\Big|}{g^n(W_n)}\right\} \leq \epsilon.
           \label{eq:appendix_phi_g_12}
        \end{equation}
        
        In the case that $f^n(W_n) \leq \widehat{f^n}(W_n)$ (i.e., $f^n(W_n) = W_n^{\sigma_n - 1}$ is convex, or $\sigma_n \in [0, 1] \bigcup [2, +\infty]$), the condition \eqref{eq:appendix_phi_f_12} is equivalent to $f^n(W_n) \leq \widehat{f^n}(W_n) \leq (1 + \epsilon)f^n(W_n)$. In the case that $f^n(W_n) \geq \widehat{f^n}(W_n)$ (i.e., $f^n(W_n) = W_n^{\sigma_n - 1}$ is concave, or $\sigma_n \in [1, 2]$), \eqref{eq:appendix_phi_f_12} can be rewritten as $(1 -\epsilon)f^n(W_n) \leq \widehat{f^n}(W_n) \leq f^n(W_n)$. Overall, for any $\sigma_n \geq 0$, we have the worst-case approximation bound for the approximator $\widehat{f^n}(W_n)$ given by $(1 -\epsilon)f^n(W_n) \leq \widehat{f^n}(W_n) \leq (1 + \epsilon)f^n(W_n)$. Similarly, we also have the approximation bound for $\widehat{g^n}(W_n)$ given by $(1 -\epsilon)g^n(W_n) \leq \widehat{g^n}(W_n) \leq (1 + \epsilon)g^n(W_n)$. Substitute the approximators to the objective function \eqref{eq:obj_assort} of the assortment problem \eqref{prob:CNL-assort}, we have:
        \begin{align*}
            \frac{\sum_{i\in [m]} \sum_{n\in [N]}(1 - \epsilon){f^n}(W_n) (\alpha_{in}x_ir_iv_i)}{\sum_{n\in [N]} (1+\epsilon){g^n}(W_n)} &\leq \frac{\sum_{i\in [m]} \sum_{n\in [N]}\widehat{f^n}(W_n) (\alpha_{in}x_ir_iv_i)}{\sum_{n\in [N]} \widehat{g^n}(W_n)} \\ &\leq \frac{\sum_{i\in [m]} \sum_{n\in [N]}(1+\epsilon){f^n}(W_n) (\alpha_{in}x_ir_iv_i)}{\sum_{n\in [N]} (1-\epsilon){g^n}(W_n)}
        \end{align*}
        for we can rewrite and obtain the approximation bounds for the objective function as
        \[
        \frac{1 - \epsilon}{1 + \epsilon}F(\bx) \leq \widehat{F}(\bx) \leq \frac{1 + \epsilon}{1 - \epsilon}F(\bx),
        \]
        which completes the proof.
    \endproof

    From the above approximation bounds, we have the worst-case performance guarantee for the approximate problem \eqref{prob:approx} as in following theorem. 
    \begin{theorem}
        For dissimilarity parameter $\sigma_n \geq 0,\ n \in [N]$,
        let $\bx^*$ be the optimal solution of the assortment problem \eqref{prob:CNL-assort} and $\widehat{\bx}$ be the optimal solution of the approximate problem \eqref{prob:approx1}. Then, we have $\frac{(1 - \epsilon)^2}{(1 + \epsilon)^2}F(\bx^*) \leq F(\widehat{\bx}) \leq F(\bx^*)$.
    \end{theorem}
    \proof{Proof.}
        As we have mentioned, relaxing the restrictions on the dissimilarity parameters make no effect on our formulations, and thus we obtain an approximate problem that is exactly the same as the linear-fractional program given in \eqref{prob:approx1}, which is equivalent to \eqref{prob:approx}.
	Therefore, it follows from Lemma \ref{lemma:general_approx_bounds} above that
	\begin{equation}
	    \frac{1 - \epsilon}{1 + \epsilon}F(\bx^*) \leq \widehat{F}(\bx^*) \quad \text{and} \quad \widehat{F}(\widehat{\bx}) \leq \frac{1 + \epsilon}{1 - \epsilon}F(\widehat{\bx}).
        \label{eq:tmp1_from_approx_bound_sigma_12}
	\end{equation}
        Since $\widehat{\bx}$ is optimal solution of the approximate problem \eqref{prob:approx1}, we have
    	$\widehat{F}(\bx^*) \leq \widehat{F}(\widehat{\bx})$. Combine with \eqref{eq:tmp1_from_approx_bound_sigma_12}, we get 
            \begin{equation}
    	    \frac{1 - \epsilon}{1 + \epsilon}F(\bx^*) \leq \frac{1 + \epsilon}{1 - \epsilon}F(\widehat{\bx}).
         \label{eq:tmp2_from_approx_bound_sigma_12}
    	\end{equation}
    	Since $x^*$ is optimal solution of 
    	the problem \eqref{prob:CNL-assort}, we have $F(\widehat{\bx}) \leq F(\bx^*)$. Combine with \eqref{eq:tmp2_from_approx_bound_sigma_12}, we get 
    	\[
    	\frac{1 - \epsilon}{1 + \epsilon}F(\bx^*) \leq \frac{1 + \epsilon}{1 - \epsilon}F(\widehat{\bx}) \leq \frac{1 + \epsilon}{1 - \epsilon}F(\bx^*).
    	\]
    	By dividing all the terms by $\frac{1 + \epsilon}{1 - \epsilon}$, we can obtain the desired result and complete the proof.
    \endproof
    The above theorem implies that for any dissimilarity parameters $\sigma_n \geq 0\ \forall n\in[N]$, we can obtain a $\frac{(1 - \epsilon)^2}{(1 + \epsilon)^2}$-optimal solution to the assortment problem \eqref{prob:CNL-assort} by solving the linear-fractional program \eqref{prob:approx1}.

        \subsection{The Applicability of the Discretization Procedure}
        We now show that the discretization procedure in Subsection \ref{subsec:discretization} can still be used for the assortment problem where the dissimilarity parameters $\sigma_n > 1, n\in[N]$, and thus the procedure is applicable for any values of dissimilarity parameters. This procedure relies on two facts that $i)$ the relative approximation error functions $\phi^{fn}(c^n_k, u, t)$ and $\phi^{gn}(c^n_k, u, t)$  have unique maximums in the sub-interval $[c^n_k, u] \subseteq [L_n, U_n]$ (Proposition \ref{proposition:unique_maximum_of_relative_error}), and $ii)$ the maximum relative error functions $\Phi^{fn}(c^n_k, u)$ and $\Phi^{gn}(c^n_k, u)$ are monotonically non-decreasing in the second endpoint $u$ of the linear piece (Lemma \ref{lemma:monotone_of_max_error}). Then, if these results hold for $\sigma_n > 1$, then we can use the same procedure for this case. 
        We have the following results, which are generalizations of the results in Subsection \ref{subsec:discretization}:

        \begin{proposition}
            For any dissimilarity parameters $\sigma_n \geq 0\ \forall n\in [N]$, the relative errors $\phi^{fn}(c^n_k, u, t)$ has a unique maximum at $$t^{fn} = \frac{(\sigma_n - 1)(f^n(c^n_k) - \gamma_k^{fn}c^n_k)}{\gamma_k^{fn}(2 - \sigma_n)}$$ and $\phi^{gn}(c^n_k, u, t)$ has a unique maximum at $$t^{gn} = \frac{\sigma_n(g^n(c^n_k) - \gamma_k^{gn}c^n_k)}{\gamma_k^{gn}(1 - \sigma_n)}$$ on the sub-interval $[c^n_k, u]$.
        \end{proposition}
        \proof{Proof.}
        For special cases where $\sigma_n = 1$ or $\sigma_n = 2$, $f^n(t) = t^{\sigma_n - 1}$ becomes linear. Therefore, we have $\phi^{fn}(c^n_k, u, t) = 0,\ \forall  t\in [L_n, U_n]$, and only the error of $\widehat{g^n}(t)$ needs to be considered.
        For $\sigma_n \neq 1$ and $\sigma_n \neq 2$, the relative error function of the approximator $\widehat{f^n}(c^n_k, u, t)$ at the point $t$ on the sub-interval $[c^n_k, u]$ is given by
        \[
        \phi^{fn}(c^n_k, u, t) = \frac{\Big|\widehat{f^n}(c^n_k, u, t) - f^n(t)\Big|}{f^n(t)} = \Bigg|\frac{\widehat{f^n}(c^n_k, u, t)}{f^n(t)} - 1\Bigg|.
        \]
        If the function $f^n(t) = t^{\sigma_n-1}$ is convex (i.e., $0 \leq \sigma_n < 1 \text{ or } \sigma_n > 2$), then we have $\widehat{f^n}(c^n_k, u, t) \geq f^n(t),\ \forall t\in [c^n_k, u]$, and thus $\phi^{fn}(c^n_k, u, t) = \frac{\widehat{f^n}(c^n_k, u, t)}{f^n(t)} - 1$. This function has a unique maximum at $$t^{fn} = \frac{(\sigma_n - 1)(f^n(c^n_k) - \gamma_k^{fn}c^n_k)}{\gamma_k^{fn}(2 - \sigma_n)},$$ as shown in Proposition \ref{proposition:unique_maximum_of_relative_error}. On the other hand, if $f^n(t)$ is concave (i.e., $1 < \sigma_n < 2$), we have $\widehat{f^n}(c^n_k, u, t) \leq f^n(t),\ \forall t\in [c^n_k, u]$, and thus $$\phi^{fn}(c^n_k, u, t) = 1 -  \frac{\widehat{f^n}(c^n_k, u, t)}{f^n(t)}.$$ In this case, we have $$\frac{\partial}{\partial t}\phi^{f^n}(c^n_k, u, t) = -\gamma_k^{fn} t^{1 - \sigma_n} - (1 - \sigma_n)\widehat{f^n}(c^n_k, u, t) t ^{-\sigma_n}.$$ It is then straightforward to solve the equation $\frac{\partial}{\partial t}\phi^{f^n}(c^n_k, u, t) = 0$ to show that $\phi^{f^n}(c^n_k, u, t)$ has a unique maximum at the same point $$t^{fn} = \frac{(\sigma_n - 1)(f^n(c^n_k) - \gamma_k^{fn}c^n_k)}{\gamma_k^{fn}(2 - \sigma_n)}.$$

        Similarly, for the cases where $\sigma_n = 0$ or $\sigma_n = 1$, $g^n(t) = t^{\sigma_n}$ is linear. Thus, $\phi^{gn}(c^n_k, u, t) = 0,\ \forall t \in [L_n, U_n]$. For $\sigma_n > 0$ and $\sigma_n \neq 1$, we can do the same things as with $\phi^{fn}(c^n_k, u, t)$ to show that $\phi^{gn}(c^n_k, u, t)$ also has a unique maximum at $$t^{gn} = \frac{\sigma_n(g^n(c^n_k) - \gamma_k^{gn}c^n_k)}{\gamma_k^{gn}(1 - \sigma_n)}$$ and complete the proof.
        \endproof

        \begin{proposition}
            For dissimilarity parameters $\sigma_n \geq 0,\ n\in[N]$, the maximum relative error functions $\Phi^{fn}(c^n_k, u) = \max_{c^n_k \leq t \leq u}\{\phi^{fn}(c^n_k, u, t)\}$ and $\Phi^{gn}(c^n_k, u) = \max_{c^n_k \leq t \leq u}\{\phi^{gn}(c^n_k, u, t)\}$ are monotonically non-decreasing in $u$.
        \end{proposition}
        \proof{Proof.}
            The proof is similar, but more general to the proof of Lemma \ref{lemma:monotone_of_max_error}. We first prove that the function $\Phi^{fn}(c^n_k, u)$ is monotonically non-decreasing by supposing to the contrary that there exists $u' > u$ such that $\Phi^{fn}(c^n_k, u) > \Phi^{fn}(c^n_k, u')$. This is equivalent to
        	\begin{equation}
        		\max_{c^n_k \leq t \leq u}\left\{\frac{\Big|\widehat{f^n}(c^n_k, u, t) - f^n(t)\Big|}{f^n(t)}\right\} \quad > \quad \max_{c^n_k \leq t \leq u'}\left\{\frac{\Big|\widehat{f^n}(c^n_k, u', t) - f^n(t)\Big|}{f^n(t)}\right\}.
        		\label{eq:abs_Phi_monotonic}
        	\end{equation}

            In the cases that $0\leq \sigma_n \leq 1$ or $\sigma_n \geq 2$, $f^n(t) = t^{\sigma_n - 1}$ is convex and monotonic in $t \in [L_n, U_n]$. It follows that $\widehat{f^n}(c^n_k, u, t) \geq f^n(t)$ and $\widehat{f^n}(c^n_k, u, t)$ is increasing in $u$, leading to the same case stated in Lemma \ref{lemma:monotone_of_max_error}.

            For $\sigma_n \in [1, 2]$, the function $f^n(t)$ is now concave and monotonic in $t\in [L_n, U_n]$. It follows that $\widehat{f^n}(c^n_k, u, t) \leq f^n(t)$ and $\widehat{f^n}(c^n_k, u, t)$ is decreasing in $u$. The Inequality \eqref{eq:abs_Phi_monotonic} is thus equivalent to
            \begin{equation}
        		\max_{c^n_k \leq t \leq u}\left\{\frac{f^n(t) - \widehat{f^n}(c^n_k, u, t)}{f^n(t)}\right\} \quad > \quad \max_{c^n_k \leq t \leq u'}\left\{\frac{f^n(t) - \widehat{f^n}(c^n_k, u', t)}{f^n(t)}\right\}.
                \label{eq:Phi_f_concave_monotonic}
            \end{equation}
            Now, since $\widehat{f^n}(c^n_k, u, t)$ is decreasing in $u$, we have $\widehat{f^n}(c^n_k, u, t) > \widehat{f^n}(c^n_k, u', t)$ for $u < u'$. Therefore,  \eqref{eq:Phi_f_concave_monotonic} can be rewritten as
            \begin{equation*}
        		\max_{c^n_k \leq t \leq u}\left\{\frac{f^n(t) - \widehat{f^n}(c^n_k, u, t)}{f^n(t)}\right\} \quad > \quad \max_{c^n_k \leq t \leq u'}\left\{\frac{f^n(t) - \widehat{f^n}(c^n_k, u, t)}{f^n(t)}\right\},
            \end{equation*}
            which contradicts the fact that $[c^n_k, u] \subset [c^n_k, u']$. Therefore, we have $\Phi^{fn}(c^n_k, u') \geq \Phi^{fn}(c^n_k, u)$ for any $u' > u$, indicating that $\Phi^{fn}(c^n_k, u)$ is monotonically non-decreasing in $u$ for any dissimilarity parameter $\sigma_n \geq 0, n\in[N]$. 
            
            Using the same contradiction arguments, we can also prove that $\Phi^{gn}(c^n_k, u)$ is monotonically non-decreasing in $u$ and complete the proof. 
        \endproof

        These above results indicate that the discretization procedure presented in Subsection \ref{subsec:discretization} can be applied to the assortment problem with any dissimilarity parameter $\sigma_n \geq 0,\ n\in [N]$. Moreover, since the proof of Theorem \ref{theorem:min_K} depends solely on the convexity and monotonicity of the functions $f^n(t)$ and $g^n(t)$, the minimum number of sub-intervals established in this theorem remains valid for $\sigma_n > 1$. In other words, for any dissimilarity parameter $\sigma_n \geq 0,\ n\in [N]$, the discretization procedure presented in Subsection \ref{subsec:discretization} partitions the interval $[L_n, U_n]$ into $K_n$ sub-intervals, where the number of sub-intervals $K_n$ is minimized for any given accuracy level $\epsilon > 0$.
    \end{appendices}
\end{document}

%% file: tables/k_bisection.tex
\begin{table}[htb]
\centering
\caption{Number of generated sub-intervals}
\label{tab:k_bisection}
\begin{tabular}{@{}lllllllllllllllll@{}}
\toprule
$\sigma_n$ &  & 0.2 & 0.2 & 0.2 & 0.3 & 0.3 & 0.3 & 0.5 & 0.5 & 0.5 & 0.7 & 0.7 & 0.7 & 0.9 & 0.9 & 0.9 \\ \midrule
$U_n$     &  & 5   & 10  & 15  & 5   & 10  & 15  & 5   & 10  & 15  & 5   & 10  & 15  & 5   & 10  & 15  \\ \midrule
$\epsilon = 0.1$   &  & 3   & 4   & 4   & 2   & 3   & 4   & 2   & 3   & 3   & 2   & 2   & 2   & 1   & 1   & 1   \\
$\epsilon = 0.05$  &  & 4   & 5   & 6   & 3   & 5   & 5   & 3   & 4   & 4   & 2   & 3   & 3   & 1   & 2   & 2   \\
$\epsilon = 0.01$  &  & 7   & 10  & 12  & 7   & 9   & 11  & 5   & 8   & 9   & 4   & 6   & 6   & 2   & 3   & 4   \\
$\epsilon = 0.005$ &  & 10  & 14  & 17  & 9   & 13  & 15  & 7   & 10  & 12  & 6   & 8   & 9   & 3   & 4   & 5   \\
$\epsilon = 0.001$ &  & 22  & 31  & 37  & 20  & 29  & 34  & 16  & 23  & 27  & 12  & 17  & 19  & 6   & 9   & 11  \\ \bottomrule
\end{tabular}
\end{table}

%% file: tables/assortment_small_new.tex
\begin{table}
\centering
\caption{Results for the small dataset of the assortment problem.}
\label{table:assort_small}
\resizebox{\textwidth}{!}{%
\begin{tabular}{@{}rllllllllllllllllll@{}}
\toprule
\multirow{4}{*}{($N,   m, \Bar{\sigma}$)} &  & \multicolumn{5}{c}{90\%}                                                                                                                                                             &  & \multicolumn{5}{c}{95\%}                                                                                                                                                          &  & \multicolumn{5}{c}{99\%}                                                                                                                                                         \\ \cmidrule(lr){3-7} \cmidrule(lr){9-13} \cmidrule(l){15-19} 
                                    &  & \multicolumn{2}{c}{\begin{tabular}[c]{@{}c@{}}No.   of inst.  \\ with \textsf{BOS}\end{tabular}} &  & \multicolumn{2}{c}{\begin{tabular}[c]{@{}c@{}}Average\\ runtime (s)\end{tabular}} &  & \multicolumn{2}{c}{\begin{tabular}[c]{@{}c@{}}No. of inst.\\ with \textsf{BOS}\end{tabular}} &  & \multicolumn{2}{c}{\begin{tabular}[c]{@{}c@{}}Average \\ runtime (s)\end{tabular}} &  & \multicolumn{2}{c}{\begin{tabular}[c]{@{}c@{}}No. of inst.\\ with \textsf{BOS}\end{tabular}} &  & \multicolumn{2}{c}{\begin{tabular}[c]{@{}c@{}}Average\\ runtime (s)\end{tabular}} \\ \cmidrule(lr){3-4} \cmidrule(lr){6-7} \cmidrule(lr){9-10} \cmidrule(lr){12-13} \cmidrule(lr){15-16} \cmidrule(l){18-19} 
                                    &  & \multicolumn{1}{c}{BIS}\quad                       & \multicolumn{1}{c}{MILP}                      &  & \multicolumn{1}{c}{BIS}                 & \multicolumn{1}{c}{MILP}                &  & \multicolumn{1}{c}{BIS}\quad                     & \multicolumn{1}{c}{MILP}                    &  & \multicolumn{1}{c}{BIS}                 & \multicolumn{1}{c}{MILP}                 &  & \multicolumn{1}{c}{BIS}\quad                     & \multicolumn{1}{c}{MILP}                    &  & \multicolumn{1}{c}{BIS}                 & \multicolumn{1}{c}{MILP}                \\ \midrule
(5, 20, 0.50)                       &  & 708                                           & \textbf{715}                                           &  & 0.67                                    & \textbf{0.21}                                    &  & 800                                         & \textbf{822}                                         &  & 0.83                                    & \textbf{0.38}                                     &  & 936                                         & \textbf{963}                                         &  & 1.60                                    & \textbf{1.09}                                    \\
(5, 20, 0.75)                       &  & 703                                           & \textbf{719}                                           &  & 0.55                                    & \textbf{0.15}                                    &  & 812                                         & \textbf{839}                                         &  & 0.69                                    & \textbf{0.23}                                     &  & 931                                         & \textbf{957}                                         &  & 1.27                                    & \textbf{0.79}                                    \\
(5, 20, 1.00)                       &  & 697                                           & \textbf{729}                                           &  & 0.42                                    & \textbf{0.10}                                    &  & 793                                         & \textbf{817}                                         &  & 0.52                                    & \textbf{0.15}                                     &  & 911                                         & \textbf{953}                                         &  & 0.95                                    & \textbf{0.46}                                    \\
(5, 30, 0.50)                       &  & 526                                           & \textbf{565}                                           &  & 1.46                                    & \textbf{0.52}                                    &  & 675                                         & \textbf{712}                                         &  & 1.99                                    & \textbf{0.89}                                     &  & 864                                         & \textbf{942}                                         &  & 6.85                                    & \textbf{3.32}                                    \\
(5, 30, 0.75)                       &  & 530                                           & \textbf{566}                                           &  & 1.15                                    & \textbf{0.36}                                    &  & 692                                         & \textbf{741}                                         &  & 1.49                                    & \textbf{0.59}                                     &  & 855                                         & \textbf{924}                                         &  & 3.86                                    & \textbf{2.01}                                    \\
(5, 30, 1.00)                       &  & 540                                           & \textbf{572}                                           &  & 0.92                                    & \textbf{0.21}                                    &  & 694                                         & \textbf{752}                                         &  & 1.17                                    & \textbf{0.36}                                     &  & 864                                         & \textbf{949}                                         &  & 2.61                                    & \textbf{1.19}                                    \\
(10, 20, 0.50)                      &  & 766                                           & \textbf{774}                                           &  & 0.30                                    & \textbf{0.25}                                    &  & 860                                         & \textbf{879}                                         &  & \textbf{0.34}                                    & 0.39                                     &  & 947                                         & \textbf{976}                                         &  & \textbf{0.53}                                    & 1.01                                    \\
(10, 20, 0.75)                      &  & 765                                           & \textbf{779}                                           &  & 0.24                                    & \textbf{0.17}                                    &  & 865                                         & \textbf{871}                                         &  & \textbf{0.28}                                    & 0.30                                     &  & 928                                         & \textbf{971}                                         &  & \textbf{0.43}                                    & 0.85                                    \\
(10, 20, 1.00)                      &  & 782                                           & \textbf{805}                                           &  & 0.18                                    & \textbf{0.11}                                    &  & 860                                         & \textbf{887}                                         &  & 0.22                                    & \textbf{0.17}                                     &  & 934                                         & \textbf{978}                                         &  & \textbf{0.34}                                    & 0.53                                    \\
(10, 30, 0.50)                      &  & 590                                           & \textbf{628}                                           &  & \textbf{0.96}                                    & 1.27                                    &  & 728                                         & \textbf{778}                                         &  & \textbf{1.16}                                    & 2.81                                     &  & 886                                         & \textbf{949}                                         &  & \textbf{1.99}                                    & 10.26                                   \\
(10, 30, 0.75)                      &  & 586                                           & \textbf{624}                                           &  & 0.79                                    & \textbf{0.66}                                    &  & 704                                         & \textbf{762}                                         &  & \textbf{0.94}                                    & 1.51                                     &  & 860                                         & \textbf{945}                                         &  & \textbf{1.61}                                    & 7.56                                    \\
(10, 30, 1.00)                      &  & 610                                           & \textbf{641}                                           &  & 0.60                                    & \textbf{0.28}                                    &  & 730                                         & \textbf{769}                                         &  & 0.74                                    & \textbf{0.59}                                     &  & 857                                         & \textbf{948}                                         &  & \textbf{1.21}                                    & 3.02                                    \\ \midrule
Total                          &  & 7803                                          & \textbf{8117}                                          &  & -                                    & -                                    &  & 9213                                        & \textbf{9629}                                        &  & -                                    & -                                     &  & 10773                                       & \textbf{11455}                                       &  & -                                    & -                                    \\ \bottomrule
\end{tabular}}
\end{table}

%% file: tables/assortment_large.tex
\begin{table}
\centering
\caption{Results on the large dataset of the assortment problem.}
\label{table:assort_large}
\resizebox{\textwidth}{!}{%
\begin{tabular}{@{}rllllllllllllllllll@{}}
\toprule
\multirow{4}{*}{($N,   m$)} &  & \multicolumn{5}{c}{90\%}                                                             &  & \multicolumn{5}{c}{95\%}                                                             &  & \multicolumn{5}{c}{99\%}                                                             \\ \cmidrule(lr){3-7} \cmidrule(lr){9-13} \cmidrule(l){15-19} 
                          &  & \multicolumn{2}{c}{{\begin{tabular}[c]{@{}c@{}}No. of \\solved inst.\end{tabular}}} &  & \multicolumn{2}{c}{{\begin{tabular}[c]{@{}c@{}}Average \\ runtime (s)\end{tabular}}} &  & \multicolumn{2}{c}{{\begin{tabular}[c]{@{}c@{}}No. of \\solved inst.\end{tabular}}} &  & \multicolumn{2}{c}{{\begin{tabular}[c]{@{}c@{}}Average \\ runtime (s)\end{tabular}}} &  & \multicolumn{2}{c}{{\begin{tabular}[c]{@{}c@{}}No. of \\solved inst.\end{tabular}}} &  & \multicolumn{2}{c}{{\begin{tabular}[c]{@{}c@{}}Average \\ runtime (s)\tablefootnote{The average runtime are  calculated over the solved instances.}\end{tabular}}} \\ \cmidrule(lr){3-4} \cmidrule(lr){6-7} \cmidrule(lr){9-10} \cmidrule(lr){12-13} \cmidrule(lr){15-16} \cmidrule(l){18-19} 
                          &  & BIS                & MILP               &  & BIS                & MILP               &  & BIS                & MILP               &  & BIS                & MILP               &  & BIS                & MILP               &  & BIS                & MILP               \\ \midrule
(5, 50)                   &  & \textbf{10}                 & \textbf{10}                 &  & 1.84               & \textbf{0.57}               &  & \textbf{10}                 & \textbf{10}                 &  & 3.15               & \textbf{1.04}               &  & \textbf{10}                 & \textbf{10}                 &  & 13.01              & \textbf{4.89}               \\
(5, 100)                  &  & \textbf{10}                 & \textbf{10}                 &  & 9.47               & \textbf{2.16}               &  & \textbf{10}                 & \textbf{10}                 &  & 16.53              & \textbf{4.93}               &  & \textbf{10}                 & \textbf{10}                 &  & 174.75             & \textbf{23.45}              \\
(5, 150)                  &  & \textbf{10}                 & \textbf{10}                 &  & 35.74              & \textbf{6.15}               &  & \textbf{10}                 & \textbf{10}                 &  & 88.76              & \textbf{13.66}              &  & \textbf{10}                 & \textbf{10}                 &  & 680.19             & \textbf{71.18}              \\
(5, 200)                  &  & \textbf{10}                 & \textbf{10}                 &  & 81.75              & \textbf{15.87}              &  & \textbf{10}                 & \textbf{10}                 &  & 198.50             &\textbf{ 26.85}              &  & 9                  & \textbf{10}                 &  & 1374.16            & \textbf{139.16}             \\
(10, 50)                  &  & \textbf{10}                 & \textbf{10}                 &  & \textbf{1.97}               & 2.30               &  & \textbf{10}                 & \textbf{10}                 &  & \textbf{2.95}               & 7.35               &  & \textbf{10}                 & \textbf{10}                 &  & \textbf{6.80}               & 66.92              \\
(10, 100)                 &  & \textbf{10}                 & \textbf{10}                 &  & 16.12              & \textbf{15.33}              &  & \textbf{10}                 & \textbf{10}                 &  & \textbf{28.11}              & 51.37              &  & \textbf{10}                 & \textbf{10}                 &  & \textbf{222.98}             & 809.01             \\
(10, 150)                 &  & \textbf{10}                 & \textbf{10}                 &  & 71.74              & \textbf{58.21}              &  & \textbf{10}                 & \textbf{10}                 &  & \textbf{180.75}             & 206.73             &  & \textbf{10}                 & 6                  &  & \textbf{1099.25}            & 1523.83            \\
(10, 200)                 &  & \textbf{10}                 & \textbf{10}                 &  & 281.18             & \textbf{131.30}             &  & \textbf{10}                 & \textbf{10}                 &  & 580.81             & \textbf{523.83}             &  & \textbf{7}                  & 3                  &  & \textbf{1243.30}            & 2302.54            \\ \midrule
Total                &  & \textbf{80}                 & \textbf{80}                 &  & -              & -              &  & \textbf{80}                 & \textbf{80}                 &  & -             & -             &  & \textbf{76}                 & 69                 &  & -             & -
             \\ \bottomrule
\end{tabular}}
\end{table}

%% file: tables/assortment_price_small.tex
\begin{table}
\centering
\caption{Results for the joint assortment optimization and pricing problem.}
\label{table:AnP}
\resizebox{\textwidth}{!}{%
\begin{tabular}{@{}rllllllllllllllllll@{}}
\toprule
\multirow{4}{*}{($N,   m, L$)} &  & \multicolumn{5}{c}{90\%}                                                                                                                                                                       &                      & \multicolumn{5}{c}{95\%}                                                                                                                                                                         &                      & \multicolumn{5}{c}{99\%}                                                                                                                                                                         \\ \cmidrule(lr){3-7} \cmidrule(lr){9-13} \cmidrule(l){15-19} 
                             &  & \multicolumn{2}{c}{\begin{tabular}[c]{@{}c@{}}No. of inst.\\ with \textsf{BOS}\end{tabular}} & \multicolumn{1}{c}{} & \multicolumn{2}{c}{\begin{tabular}[c]{@{}c@{}}Average\\ runtime (s)\end{tabular}} & \multicolumn{1}{c}{} & \multicolumn{2}{c}{\begin{tabular}[c]{@{}c@{}}No. of   inst.\\ with \textsf{BOS}\end{tabular}} & \multicolumn{1}{c}{} & \multicolumn{2}{c}{\begin{tabular}[c]{@{}c@{}}Average\\ runtime (s)\end{tabular}} & \multicolumn{1}{c}{} & \multicolumn{2}{c}{\begin{tabular}[c]{@{}c@{}}No. of   inst.\\ with \textsf{BOS}\end{tabular}} & \multicolumn{1}{c}{} & \multicolumn{2}{c}{\begin{tabular}[c]{@{}c@{}}Average\\ runtime (s)\end{tabular}} \\ \cmidrule(lr){3-4} \cmidrule(lr){6-7} \cmidrule(lr){9-10} \cmidrule(lr){12-13} \cmidrule(lr){15-16} \cmidrule(l){18-19} 
                             &  & BIS                                      & MILP                                     &                      & BIS                                     & MILP                                    &                      & BIS                                       & MILP                                      &                      & BIS                                     & MILP                                    &                      & BIS                                       & MILP                                      &                      & BIS                                     & MILP                                    \\ \midrule
(5, 20, 2)                   &  & \textbf{61}                                       & \textbf{61}                                       &                      & 0.38                                    & \textbf{0.21}                                    &                      & \textbf{70}                                        & \textbf{70}                                        &                      & 0.62                                    & \textbf{0.43}                                    &                      & 95                                        & \textbf{96}                                        &                      & \textbf{1.27}                                    & 1.76                                    \\
(5, 20, 5)                   &  & 55                                       & \textbf{58}                                       &                      & \textbf{1.10}                                    & 2.97                                    &                      & 65                                        & \textbf{69}                                        &                      & \textbf{1.53}                                    & 5.45                                    &                      & 83                                        & \textbf{84}                                        &                      & \textbf{4.89}                                    & 21.01                                   \\
(5, 30, 2)                   &  & \textbf{23}                                       & \textbf{23}                                       &                      & 1.14                                    & \textbf{0.84}                                    &                      & \textbf{52}                                        & 51                                        &                      & 1.81                                    & \textbf{1.39}                                    &                      & 84                                        & \textbf{89}                                        &                      & \textbf{6.46}                                    & 8.81                                    \\
(5, 30, 5)                   &  & \textbf{30}                                       & 29                                       &                      & \textbf{15.35}                                   & 110.36                                  &                      & \textbf{35}                                        & \textbf{35}                                        &                      & \textbf{11.62}                                   & 92.90                                   &                      & 74                                        & \textbf{78}                                        &                      & \textbf{25.91}                                   & 162.84                                  \\
(10, 20, 2)                  &  & 75                                       & \textbf{77}                                       &                      & 0.20                                    & \textbf{0.13}                                    &                      & 88                                        & \textbf{89}                                        &                      & 0.28                                    & \textbf{0.19}                                    &                      & 94                                        & \textbf{96}                                        &                      & 0.53                                    & \textbf{0.50}                                    \\
(10, 20, 5)                  &  & 61                                       & \textbf{64}                                       &                      & \textbf{0.43}                                    & 1.08                                    &                      & 68                                        & \textbf{76}                                        &                      & \textbf{0.68}                                    & 2.45                                    &                      & 80                                        & \textbf{88}                                        &                      & \textbf{2.49}                                    & 8.58                                    \\
(10, 30, 2)                  &  & 54                                       & \textbf{55}                                       &                      & \textbf{0.44}                                    & 0.74                                    &                      & \textbf{62}                                        & 61                                        &                      & \textbf{0.77}                                    & 2.01                                    &                      & 84                                        & \textbf{89}                                        &                      & \textbf{1.56}                                    & 7.55                                    \\
(10, 30, 5)                  &  & 40                                       & \textbf{46}                                       &                      & \textbf{2.07}                                    & 23.24                                   &                      & 49                                        & \textbf{59}                                        &                      & \textbf{5.07}                                    & 56.19                                   &                      & 64                                        & \textbf{72}                                        &                      & \textbf{14.14}                                   & 299.82                                  \\ \midrule
Total                   &  & 399                                      & \textbf{413}                                      &                      & -                                    & -                                   &                      & 489                                       & \textbf{510}                                       &                      & -                                    & -                                   &                      & 658                                       & \textbf{692}                                       &                      & -                                    & -                                   \\ \bottomrule
\end{tabular}
}
\end{table}

%% file: tables/mixture.tex
\begin{table}
\centering
\caption{Results of the MILP method for the assortment problem under the mixture of CNL models.}
\label{table:assort_mixture}
\resizebox{\textwidth}{!}{%
\begin{tabular}{@{}rlllllllllllllll@{}}
\toprule
\multirow{3}{*}{($N,   m, T$)} &  & \multicolumn{4}{c}{90\%}                                                                                                                                                & \multicolumn{1}{c}{} & \multicolumn{4}{c}{95\%}                                                                                                                                                & \multicolumn{1}{c}{} & \multicolumn{4}{c}{99\%}                                                                                                                                                \\ \cmidrule(lr){3-6} \cmidrule(lr){8-11} \cmidrule(l){13-16} 
                             &  & \multicolumn{2}{l}{\begin{tabular}[c]{@{}l@{}}No. of inst.\\ with \textsf{BOS}\end{tabular}} & \multicolumn{2}{l}{\begin{tabular}[c]{@{}l@{}}Avarage\\ runtime (s)\end{tabular}} &                      & \multicolumn{2}{l}{\begin{tabular}[c]{@{}l@{}}No. of inst.\\ with \textsf{BOS}\end{tabular}} & \multicolumn{2}{l}{\begin{tabular}[c]{@{}l@{}}Avarage\\ runtime (s)\end{tabular}} &                      & \multicolumn{2}{l}{\begin{tabular}[c]{@{}l@{}}No. of inst.\\ with \textsf{BOS}\end{tabular}} & \multicolumn{2}{l}{\begin{tabular}[c]{@{}l@{}}Avarage\\ runtime (s)\end{tabular}} \\ \midrule
(5, 20, 2)                   &  & \quad\quad                                        & 70                                        & \quad                                      & 0.26                                      &                      &  \quad\quad                                       & 95                                        & \quad                                     & 0.54                                       &                      &  \quad\quad                                       & 100                                       &\quad                                      & 2.60                                       \\
(5, 20, 5)                   &  &                                         & 81                                        &                                       & 1.61                                      &                      &                                         & 85                                        &                                      & 3.47                                       &                      &                                         & 98                                        &                                      & 16.82                                      \\
(5, 20,10)                  &  &                                         & 83                                        &                                       & 6.06                                      &                      &                                         & 91                                        &                                      & 12.38                                      &                      &                                         & 98                                        &                                      & 54.44                                      \\
(5, 30, 2)                   &  &                                         & 63                                        &                                       & 0.78                                      &                      &                                         & 78                                        &                                      & 2.16                                       &                      &                                         & 94                                        &                                      & 14.85                                      \\
(5, 30, 5)                   &  &                                         & 73                                        &                                       & 15.98                                     &                      &                                         & 86                                        &                                      & 34.55                                      &                      &                                         & 99                                        &                                      & 143.61                                     \\
(5, 30,10)                  &  &                                         & 74                                        &                                       & 60.72                                     &                      &                                         & 91                                        &                                      & 124.01                                     &                      &                                         & 100                                       &                                      & 519.82                                     \\
(10, 20, 2)                  &  &                                         & 81                                        &                                       & 0.19                                      &                      &                                         & 93                                        &                                      & 0.37                                       &                      &                                         & 96                                        &                                      & 1.41                                       \\
(10, 20, 5)                  &  &                                         & 85                                        &                                       & 0.90                                      &                      &                                         & 94                                        &                                      & 1.92                                       &                      &                                         & 98                                        &                                      & 8.59                                       \\
(10, 20,10)                 &  &                                         & 92                                        &                                       & 2.96                                      &                      &                                         & 96                                        &                                      & 6.48                                       &                      &                                         & 99                                        &                                      & 31.04                                      \\
(10, 30, 2)                  &  &                                         & 69                                        &                                       & 0.93                                      &                      &                                         & 85                                        &                                      & 2.22                                       &                      &                                         & 98                                        &                                      & 12.56                                      \\
(10, 30, 5)                  &  &                                         & 72                                        &                                       & 7.15                                      &                      &                                         & 83                                        &                                      & 15.28                                      &                      &                                         & 96                                        &                                      & 78.93                                      \\
(10, 30,10)                 &  &                                         & 81                                        &                                       & 21.04                                     &                      &                                         & 94                                        &                                      & 48.75                                      &                      &                                         & 99                                        &                                      & 260.73                                     \\ \midrule
Total                   &  &                                         & 924                                       &                                       &  -                                     &                      &                                         & 1071                                      &                                      &    -                                   &                      &                                         & 1175                                      &                                      &  -                                    \\ \bottomrule
\end{tabular}}
\end{table}